\documentclass[preprint]{informs3}

\DoubleSpacedXI 
\usepackage{bm}
\usepackage{natbib}
\usepackage{placeins}
\usepackage[dvipsnames]{xcolor}
\usepackage{subcaption}
\usepackage{flexisym}
\usepackage{dsfont}
\usepackage{diagbox}
\usepackage{tikz}
\usepackage{pgfplots}
\pgfplotsset{compat=1.17}
\usetikzlibrary{backgrounds,matrix}
\usepackage{mathtools}

\def\E{{\mathbb E}}

\newcommand{\abs}[1]{\left\lvert#1\right\rvert}
\newcommand{\norm}[1]{\left\lVert#1\right\rVert}
\newcommand{\ex}{\text{exp}}
\newcommand{\expemph}{\text{\emph{exp}}}
\newcommand{\control}{\text{con}}
\newcommand{\controlemph}{\text{\emph{con}}}
\newcommand{\treatment}{\text{tre}}
\newcommand{\treatmentemph}{\text{\emph{tre}}}

\makeatletter
\newcommand{\leqnomode}{\tagsleft@true\let\veqno\@@leqno}

\usepackage{endnotes}
\let\footnote=\endnote
\let\enotesize=\normalsize
\def\notesname{Endnotes}%
\def\makeenmark{\hbox to1.275em{\theenmark.\enskip\hss}}
\def\enoteformat{\rightskip0pt\leftskip0pt\parindent=1.275em
  \leavevmode\llap{\makeenmark}}

\usepackage{natbib}
 \bibpunct[, ]{(}{)}{,}{a}{}{,}%
 \def\bibfont{\small}%
 \def\BIBand{and}%

 \usepackage{hyperref}

\TheoremsNumberedThrough     
\ECRepeatTheorems

\EquationsNumberedThrough   

\begin{document}

\RUNAUTHOR{Delarue and Karakolios}

\RUNTITLE{Pricing Experiments in Matching Marketplaces under Interference}

\TITLE{Pricing Experiments in Matching Marketplaces under Interference: Designs and Estimators}

\ARTICLEAUTHORS{
\AUTHOR{Arthur Delarue}
\AFF{H. Milton Stewart School of Industrial and Systems Engineering, Georgia Institute of Technology, Atlanta, GA 30332, \EMAIL{arthur.delarue@isye.gatech.edu}} 
\AUTHOR{Kleanthis Karakolios}
\AFF{Georgia Institute of Technology, Atlanta, GA 30332, \EMAIL{kleanthis@gatech.edu}}
}

\ABSTRACT{%
Interference between treated and untreated units is a source of bias in marketplace experiments. In this paper, we specifically consider pricing interventions, in which a platform seeks to adjust base pricing levels at the marketplace level in order to increase demand. In a matching marketplace, this type of experiment leads to a crucial design question: should the platform match treated and untreated units differently because they paid different prices? We find that standard estimation techniques are biased, but the sign of this bias depends strongly on this design choice. Bias can be reduced by using the ``shadow price estimator'', which relies on the optimal dual solution of the platform's supply-demand matching problem --- especially when the platform chooses to ignore pricing differences at matching time. We validate our findings both theoretically in a fluid limit setting, and numerically in a finite-sample setting.
}%

\maketitle

\section{Introduction}

Experimentation is a fundamental decision-making tool for online platforms. Randomized controlled trials (or A/B tests) are often used to evaluate the potential gains from any proposed improvement to the platform. Unfortunately, experiments in online settings can suffer from \emph{interference}. Unlike, say, clinical trials, where each patient's outcome depends only on their own treatment, platform users can interact in such a way that treating one user may affect the outcome of another. Interference biases estimation --- sometimes so much that the bias exceeds the effect size \citep{holtz2023reducingInterferenceBiasClusterRandomization}.

Interference can be difficult to quantify in a general platform setting, but it can sometimes be modeled by leveraging the specific structure of a marketplace. For example, in matching marketplaces (such as ride-hailing marketplaces) where the platform uses a linear programming algorithm to match demand units with supply units, the platform can use information about the matching duals to quantify and reduce interference bias. \cite{bright2024reducing} explored this idea in the context of simple demand-boosting interventions to propose a new estimator which provably reduces bias as compared to standard techniques.

In this paper, we revisit the approach of \cite{bright2024reducing} in the more general context of pricing interventions, where the goal of the platform is to decide whether to adjust base pricing levels. To achieve this, the platform must accurately quantify the tradeoff between boosting total demand but reducing revenue per request. Base pricing adjustments apply to the entire marketplace, so the platform seeks to measure the global treatment effect (the total system value when everyone is treated minus the total system value when no one is).

In the experiment state, however, treated and untreated demand units, distinguishable by the revenue they bring to the platform, may coexist. This poses a core design question: should the matching algorithm treat these units as indistinguishable for the sake of the experiment, or should the platform take into account the fact that units paying a lower price bring in less revenue? We find that this seemingly innocuous question can have enormous practical implications.

Our work first establishes that in both settings, the standard estimator can be de-biased under certain conditions by using shadow-price based techniques. We then compare the conditions for the bias-reducing shadow price estimator to be unbiased for any choice of treatment fraction. Surprisingly, the two designs differ substantially: ignoring revenue differences in matching means that estimation is unbiased whenever the treatment effect is small enough relative to the total demand in the marketplace --- a condition often met in practice. In contrast, taking into account revenue differences in matching leads to unbiased estimation only when there is little or no interference in the platform to begin with --- a much stronger requirement.

The difference in these designs is especially surprising because the platform may not even be consciously choosing one design over another --- indeed, if the pricing and matching systems are run separately (a common practice in large ride-hailing firms), the choice of whether to match treated and untreated units differently may come down to an arbitrary software engineering decision. The results of this paper advocate for a more principled approach.

Our analysis relies on a simple bipartite matching model of a marketplace, with Poisson arrivals of demand and supply over a matching cycle. To derive generalizable insights, we focus on the fluid limit of the system, corresponding to a ``thick'' marketplace with many units of all types. We also analyze the performance of our proposed designs and estimators numerically in the finite-sample setting.

\section{Related Work}

Because technology companies increasingly rely on A/B testing to evaluate new strategies \citep{thomke2020building}, experimentation has received significant attention in the operations management literature. Historically, randomized controlled trials' perhaps most successful application is to inform medical decisions, for example in developing new drugs \citep{rosenberger2015randomization}. In this setting, it is reasonable to assume that treating a patient cannot improve or worsen the outcome of another patient: this Stable Unit Treatment Value Assumption (SUTVA) turns out to be a fundamental condition for the validity of many inference techniques \citep{imbensRubin2015causalInference}.

In marketplace experiments, however, SUTVA may not hold. The classic example is ``cannibalization,'' whereby treated demand units consume supply units that would otherwise have been available for control demand units, thus changing their outcome. In recent years, many approaches have been proposed to reduce bias, involving both new designs and new estimators. On the design side, a popular approach is to change the unit of randomization to reduce interactions between treated and untreated units. Clustering approaches ensure that units in the same neighborhood or region are assigned to the same treatment group \citep{chamandy2016experimentation,holtz2023reducingInterferenceBiasClusterRandomization}; meanwhile, switchback designs assign units arriving to the platform in the same time interval to the same treatment group \citep{bojinov2023design,xiong2023data}. \cite{jia2023clustered} propose a spatio-temporal clustering approach which combines both ideas. More complex designs may also reduce bias, such as two-sided designs that affect both sides of the marketplace simultaneously \citep{johari2022experimentalDesign2sidedPlatformsBiasAnalysis,bajari2023experimentalDesignMarketplaces}. Instead of proposing new designs, some works propose new estimators, for example based on Q-learning \citep{farias2022markovian} or linear programming duality \citep{bright2024reducing}.

Several papers focus on pricing interventions (understanding whether and how to adjust prices to increase total profit). As in \cite{wager2021experimenting}, the goal is usually to learn a \emph{marginal} treatment effect, i.e., the effect of an infinitesimal change in prices, which can then inform decision-making in a gradient-descent-like approach. \cite{hu2022average} and \cite{munro2021treatmentEffectsMarketEquilibrium} show that such marginal effects can be perfectly decomposed into a direct term (the effect of treating a unit on that unit's outcome) and an indirect term (the effect of treating a unit on other units' outcomes). \cite{li2023experimenting} analyze a switchback design to estimate marginal effects in systems where interference arises from queueing congestion. Our work differs from these approaches by considering general (not necessarily marginal) price interventions --- an analysis made tractable by the structure of the downstream matching model. Though our work is methodologically most similar to \cite{bright2024reducing}, we reach similar conclusions to \cite{dhaouadi2023price}. This work highlights that in pricing experiments, cannibalization can lead not just to incorrect estimation but also to incorrect decisions (i.e., suggesting that prices should be increased when in fact the opposite is true) --- a result that mirrors our own finding that the sign of the bias of the standard estimator is dependent on the exact experimental design selected. In a longer-term experimental setting, \cite{le2023price} also find that pricing interventions can lead to more complex interference patterns.

In this paper, we focus on matching marketplaces, where supply and demand units are assigned to each other by a centralized algorithm. Matching marketplaces are common in some domains such as ride-hailing \citep[see][for other examples]{shi2023optimal}. A key feature of matching-based marketplaces is that they must make both pricing and matching decisions. In practice, these decisions are often made by different systems owned by different teams \citep{qin2020ride,azagirre2024better}. Similarly, pricing and matching decisions are usually considered separately in the literature \citep[see][for a recent review]{lobel2021revenue}, with some works focusing on dynamic pricing \citep{banerjee2015pricing,castillo2017surge,bimpikis2019spatial,freund2021pricing} while others focus on matching \citep{wang2024demand,ozkan2020dynamic}. Some joint approaches have been proposed \citep[e.g.,][]{yan2020dynamic,ma2020spatio}, but tractability often requires significant simplifications, usually in the matching problem \citep{yan2023pricing}. In contrast, this paper considers the impact of simple pricing interventions (choosing between two base pricing levels) on a general bipartite matching model.

\section{Model Description}\label{sec:model-description}

We now describe our matching marketplace model and the associated experimental setting. We first describe the supply and demand arrival processes and how they are affected by the proposed treatment. We then describe the matching algorithm and the experimental setting.

\subsection{Supply and Demand Arrival Processes}

As in~\citet{bright2024reducing}, we consider a two-sided matching marketplace where demand and supply units arrive on the platform independently and randomly over time, observe relevant information provided by the platform, then decide whether or not to submit a match \emph{request}. The platform aggregates these requests over a period of time called a matching cycle, then matches supply and demand units at the end of the cycle using a matching algorithm. We model the matching algorithm as a linear optimization problem maximizing the total value generated by the matched units. 

On the demand side, we consider a finite number of demand types indexed by $i=1, \hdots, n_d$. We assume that demand intent units of type $i$ arrive to the platform independently according to a Poisson process with parameter $\tilde \lambda_i$. Upon arrival, a unit receives information about the state of the marketplace, and decides whether or not to submit a request. The goal of the platform is to understand how to affect the intent unit's request probability. The demand intent unit may be exposed to \emph{control} (the default experience on the marketplace). In this case, it will request a match with probability $0 < p_i < 1$, or irrevocably exits the marketplace with probability $1-p_i$. All requests are assumed independent of each other. Alternatively, the demand intent unit may be exposed to \emph{treatment} (in the form of a price reduction). In this case, it will request a match with probability $p_i+q_i>p_i$. Demand intent units requesting a match are called demand units, and their arrival process to the platform is effectively a thinned Poisson process with arrival rate $\lambda_i=\tilde \lambda_i p_i$ under control, and $\lambda_i + \beta_i = \tilde \lambda_i (p_i + q_i)$ under treatment, for each demand type $i \in [n_d]$. In total, the demand arrival rates are denoted by $\bm{\lambda} = (\lambda_1, \hdots, \lambda_{n_d})$ under control, and $\bm{\lambda} + \bm{\beta} = (\lambda_1+\beta_1, \hdots, \lambda_{n_d}+\beta_{n_d})$ under treatment. For a given demand arrival rate $\bm{\lambda}$, we denote by $\bm{D}^{\tau,\bm{\lambda}}$ the random vector corresponding to the number of demand units of each type; $D_i^{\tau, \bm{\lambda}}$ follows a Poisson distribution with parameter $\lambda_i\tau$, where $\tau$ is a scaling parameter which controls the relative abundance of demand and supply units.

Similarly, on the supply side each unit belongs to one of a finite number of supply types indexed by $j=1, \hdots, n_s$. The units arrive to the platform  independently of each other according to a Poisson process with parameter $\gamma_j$ for each type of supply $j\in[n_s]$. The vector of supply arrival rates is denoted by $\bm{\gamma} = (\gamma_1, \hdots, \gamma_{n_s})$. Supply units may also choose to observe the state of the platform and leave before they are matched --- for simplicity, the arrival rate $\gamma_j$ already accounts for these possible defections. We denote by $S_j^{\tau, \bm{\gamma}}\sim \text{Poisson}(\gamma_j\tau)$ the random variable corresponding to the number of supply units arriving to the platform in a matching cycle.

\subsection{Matching Mechanism}

In the control state, we assume that matching a demand unit of type $i$ to a supply unit of type $j$ yields value $v_{i,j}$ to the platform.
A particular feature of our model is that demand units may bring a different value to the platform under treatment, because exposing a unit to treatment has a cost --- for example, treated units paying a discounted price may bring less value to the platform. 

We will consider two ways to model the cost of the treatment intervention: (i) a \emph{fixed-cost} model, in which the platform pays a fixed cost $\kappa$ for each matched demand unit exposed to treatment; (ii) a \emph{proportional-cost} model, in which the platform pays a proportional cost $\alpha v_{i,j}$ for each demand unit of type $i$ exposed to treatment and subsequently matched to a supply unit of type $j$. We assume that $\kappa < \min_{i\in[n_d],j\in[n_s]}v_{i,j}$, and $0 \le \alpha < 1$.

Denoting by $\bm{d}=(d_1,\ldots,d_{n_d})\in\mathbb{Z}^{n_d}$ the number of demand units of each type, and by $\bm{s}\in\mathbb{Z}^{n_s}$ the number of supply units of each type, and assuming that all demand units are exposed to control, the platform solves the following value-maximizing linear program to match demand and supply units:
\begin{subequations} 
\begin{align}
    \label{lp:stochastic-primal}
    \Phi_{\text{CE}}(\bm{d}, \bm{s}) =  \max \quad&\sum_{i=1}^{n_d} \sum_{j=1}^{n_s} v_{i,j} x_{i,j} &  \\
                   \text{s.t.} \quad&   \sum_{j=1}^{n_s} x_{i,j} \leq d_i  & \forall i \in [n_d], \label{seq:demand-constraints}\\
     &        \sum_{i=1}^{n_d} x_{i,j} \leq s_j  & \forall j \in [n_s], \label{seq:supply-constraints}\\
      &       x_{i,j}\geq 0  & \forall i \in [n_d] , \; \forall j  \in [n_s]
\end{align}
\end{subequations}
where constraints~\eqref{seq:demand-constraints} and \eqref{seq:supply-constraints} respectively ensure feasibility from the demand and supply sides.

If all demand units are exposed to treatment, the matching weights change to $(1-\alpha)v_{i,j}$ under the proportional-cost model, and to $v_{i,j}-\kappa$ under the fixed-cost model. However, it is easy to show that the optimal matching does not change.

\begin{lemma} \label{lem:CE-unique}
    Under both cost models, for any allowed values of $\kappa$ or $\alpha$, the optimal matching is obtained by solving $\Phi_{\text{\emph{CE}}}(\cdot)$, ignoring all treatment costs.
\end{lemma}

The ``CE'' subscript thus refers to the fact that \textbf{C}ost considerations are \textbf{E}xcluded from the matching LP. Given the optimal matching, the platform can then subtract the total cost of the intervention, given by $\alpha \Phi_{\text{CE}}(\bm{d}, \bm{s})$ under the proportional-cost model and $\kappa\sum_{i=1}^{n_d}\sum_{j=1}^{n_s} x_{i,j}$ under the fixed-cost model. In general, we expect the increase in demand arrival rate $\bm{\beta}$ to depend on the discount $\alpha$ or $\kappa$. However, because we are measuring the impact of a \emph{particular value} of the discount $\alpha$ or $\kappa$, we do not need to specify any parametric relationship between $\bm{\beta}$ and $\alpha$ or $\kappa$.

The platform seeks to estimate the difference in total system value between the state where all demand units are exposed to control and the state where all demand units are exposed to treatment. However, it cannot simultaneously observe both of these states. Instead, it will run a randomized experiment, in which each demand intent unit of type $i$ is randomly exposed to control or treatment with some probability. In the experiment state, both treated and untreated units co-exist in the platform --- something which does not happen in either the control state or the treatment state. The platform therefore faces a choice on how to treat these coexisting units: the first option is to include the treatment cost directly in the matching function; the second is to ignore the treatment cost at match time, and estimate it ex-post using the experimental data.

In the first case, given $\bm{d}^{\control}$ control demand units of each type, $\bm{d}^{\treatment}$ treatment demand units of each type, and $\bm{s}$ supply units of each type, the platform solves the following modified matching linear program:
\begin{subequations}
\begin{align}\label{lp:stochastic-primal-coupondiscount}
    \Phi_{\text{CI}} (\bm{d}^{\control}, \bm{d}^{\treatment}, \bm{s})= \max \quad & \sum_{i=1}^{n_d}\sum_{j=1}^{n_s} v^\control_{i,j} x_{i,j}^{\control} + v^\treatment_{i,j} x_{i,j}^{\treatment} & \\
        \text{s.t.} \quad& \sum_{j=1}^{n_s} x_{i,j}^\treatment \leq d_i^{\treatment} & \forall i \in [n_d], \label{seq:demand-constraints-treatment-coupondiscount}\\
                         & \sum_{j=1}^{n_s} x_{i,j}^\control \leq d_i^{\control} & \forall i \in [n_d], \label{seq: demand-constraints-control-coupondiscount}\\
                         & \sum_{i=1}^{n_d} x_{i,j}^\control+x_{i,j}^\treatment \leq s_j & \forall j \in [n_s], \label{seq:supply-constraints-coupondiscount} \\
                         & x_{i,j}^\control \geq 0 \, , \, x_{i,j}^\treatment \geq 0  & \forall i \in [n_d] , \; \forall  j \in [n_s], 
\end{align}
\end{subequations}
where $v_{i,j}^\control=v_{i,j}$ designates the value obtained from matching a control demand unit of type $i$ to a supply unit of type $j$, while $v_{i,j}^\treatment$ denotes the value obtained from matching a treatment demand unit of type $i$ to a supply unit of type $j$.  Under the fixed-cost model, $v_{i,j}^\treatment=v_{i,j} - \kappa$, while under the proportional-cost model, $v_{i,j}^\treatment=(1-\alpha)v_{i,j}$. The ``CI'' subscript indicates that \textbf{C}osts are \textbf{I}ncluded in the matching formulation. To differentiate between the two cost models, we let $m$ designate the cost model under study, with $m=\text{``fixed"}$ referring to the fixed-cost model and $m=\text{``prop"}$ referring to the proportional-cost model, and we use $\Phi_{\text{CI}}^{m}(\cdot)$ to refer to the cost-included matching problem under cost model $m$, as well as $v_{i,j}^{\treatment,m}$ to refer to the matching values for treated demand units under cost model $m$.

In the second case, the platform instead solves $\Phi_{\text{CE}}(\bm{d}^{\control}+ \bm{d}^{\treatment}, \bm{s})$ to match supply and demand, then separately estimates the cost of the intervention using experimental data. We describe both cases in detail in the following subsection.

\subsection{Estimating Treatment Effects}

\paragraph{Global treatment effect.} The platform's overarching objective is to understand whether the intervention under study boosts total system value, and by how much. The simplest way to accomplish this is by estimating the global treatment effect, i.e., the difference in total system value between the state where all demand units are exposed to control and the state where all demand units are exposed to treatment.

\begin{definition}
The global treatment effect of an intervention under cost model $m$ is defined as
\begin{equation}\label{eq:gte}
    \Delta^{\tau,m} = \frac{1}{\tau} \mathbb{E}\left[
    \Phi^{m}_{\text{CI}}(\bm{0},\bm{D}^{\tau,\bm{\lambda}+\bm{\beta}}, \bm{S}^{\tau, \bm{\gamma}}) -
    \Phi^{m}_{\text{CI}}(\bm{D}^{\tau,\bm{\lambda}},\bm{0}, \bm{S}^{\tau, \bm{\gamma}})\right]
\end{equation}
\end{definition}

The total system value includes both the value obtained by the platform from matching supply and demand units and the possible cost of treating some or all of the demand units. In global control, there are $D_i^{\tau,\bm{\lambda}}$ control units of type $i$, zero treatment units of any type, and ${S}_j^{\tau, \bm{\gamma}}$ supply units of type $j$, so the expected total system value under cost model $m$ is $\mathbb{E}\left[\Phi^{m}_{\text{CI}}(\bm{D}^{\tau,\bm{\lambda}},\bm{0}, \bm{S}^{\tau, \bm{\gamma}})\right]$. Similarly, in global treatment, there are zero control units of any type, $D_i^{\tau,\bm{\lambda}+\bm{\beta}}$ treatment units of type $i$, and ${S}_j^{\tau, \bm{\gamma}}$ supply units of type $j$, so the expected total system value under cost model $m$ is $\mathbb{E}\left[\Phi^{m}_{\text{CI}}(\bm{0},\bm{D}^{\tau,\bm{\lambda}+\bm{\beta}}, \bm{S}^{\tau, \bm{\gamma}})\right]$. 

An alternative definition of the global treatment effect separates the value component from the cost component. Under cost model $m$, let $c^m(\bm{d}^\control, \bm{d}^\treatment, \bm{s})$ designate the cost of an intervention that leads to $d_i^\control$ demand units of type $i$ in control, $d_i^\treatment$ demand units of type $i$ in treatment, and $s_j$ supply units of type $j$. By definition,
\begin{align*}
    c^{\text{fixed}}(\bm{d}^\control, \bm{d}^\treatment, \bm{s}) &= \kappa \sum_{i=1}^{n_d}\sum_{j=1}^{n_s} x_{i,j}^{\treatment, \text{fixed} }\\
    c^{\text{prop}}(\bm{d}^\control, \bm{d}^\treatment, \bm{s}) &= \sum_{i=1}^{n_d}\sum_{j=1}^{n_s} \alpha v_{i,j}^\control x_{i,j}^{\treatment, \text{prop} }
\end{align*}
where $x_{i,j}^{\treatment, \text{prop} }$ refers to the optimal solution of the cost-included matching problem $\Phi^{\text{prop}}_{\text{CI}}(\bm{d}^\control, \bm{d}^\treatment, \bm{s})$ with the proportional cost model, and  $x_{i,j}^{\treatment, \text{fixed} }$ to the optimal solution of the cost-included matching problem $\Phi^{\text{fixed}}_{\text{CI}}(\bm{d}^\control, \bm{d}^\treatment, \bm{s})$  with the fixed cost model.

\begin{remark}
    \label{rk:gte-equivalent-definition}
    The global treatment effect is equivalently defined as
    \[
    \Delta^{\tau,m} = \frac{1}{\tau} \mathbb{E}\left[\Phi_{\text{CE}}(\bm{D}^{\tau,\bm{\lambda}+\bm{\beta}}, \bm{S}^{\tau, \bm{\gamma}}) - c^m(\bm{0}, \bm{D}^{\tau,\bm{\lambda}+\bm{\beta}}, \bm{S}^{\tau, \bm{\gamma}}) - \Phi_{\text{CE}}(\bm{D}^{\tau,\bm{\lambda}}, \bm{S}^{\tau, \bm{\gamma}})\right]
    \]
\end{remark}

Remark~\ref{rk:gte-equivalent-definition} relies on the observation that the cost of no intervention (i.e., no treated demand units) is always zero, and on the fact that the optimal matching solution does not change under either cost model as long as the marketplace hosts only control or treatment units, but not both.

\paragraph{Experimental data.} As mentioned, the platform cannot simultaneously observe the total system value under global control and global treatment. Instead, it can run a randomized experiment, observe the state of the system in an intermediate configuration, and gather data to estimate the global treatment effect. More precisely, in the experiment state, each demand intent unit of type $i$ is randomly exposed to control or treatment with probability $\rho$. As a result, the number of treatment units of each type is given by $\bm{D}^{\tau,\treatment}\sim \text{Poisson}(\tau\rho(\bm\lambda + \bm\beta))$, and the number of control units of each type is given by $\bm{D}^{\tau,\control}\sim\text{Poisson}(\tau(1-\rho)\bm\lambda)$. The vector $\bm{D}^{\tau,\ex}=\bm{D}^{\tau,\treatment}+\bm{D}^{\tau,\control}$ designates the total number of units of each type.

When designing the experiment, the platform must decide whether to match demand to supply using (i) cost-included or (ii) cost-excluded matching. In case (i), the platform solves $\Phi_{\text{CI}}^{m}(\bm{D}^{\tau,\control}, \bm{D}^{\tau,\treatment},\bm{S}^{\tau,\bm\gamma})$; it obtains \emph{primal data}, i.e., the number of control demand units of type $i$ assigned to supply units of type $j$, denoted by $X_{i,j}^{\tau,\control,\text{CI},m}$ (and analogously, $X_{i,j}^{\tau,\treatment,\text{CI},m}$). It also obtains \emph{dual data}, corresponding to the certificate of optimality provided by most LP solvers. By strong duality, the matching function can be equivalently written as
\begin{subequations}
    \begin{align}
        \Phi_{\text{CI}}^m (\bm{d}^{\control}, \bm{d}^{\treatment}, \bm{s})= \min \quad & \sum_{i=1}^{n_d} a_i^{\control, m} d_i^\control + \sum_{i=1}^{n_d}a_i^{\treatment, m} d_i^\treatment + \sum_{j=1}^{n_s} b_j^{m} s_j & \\
        \text{s.t.}\quad & a_i^{\control, m} + b_j^{m} \ge v_{i,j}^\control & \forall i\in[n_d], j\in[n_s]\\
        & a_i^{\treatment, m } + b_j^{m} \ge v_{i,j}^{\treatment, m} & \forall i\in[n_d], j\in[n_s]\\
        & a_i^{\control, m }, a_i^{\treatment, m} \ge 0 & \forall i\in[n_d]\\
        & b_j^{m } \ge 0 & \forall j\in[n_s].
    \end{align}
\end{subequations}
We denote by $\bm{A}^{\tau,\control,\text{CI},m}$, $\bm{A}^{\tau,\treatment,\text{CI},m}$ and $\bm{B}^{\tau,\exp,\text{CI},m}$ the optimal shadow prices of demand and supply obtained from the optimal cost-included matching at the experiment state.

In case (ii), the platform will solve $\Phi_{\text{CE}}(\bm{D}^{\tau,\ex}, \bm{S}^{\tau,\bm\gamma})$; it will obtain different primal data, this time corresponding to the number of units of type $i$ (both control and treatment) matched to supply units of type $j$, and denoted by $X^{\tau,\exp,\text{CE}}_{i,j}$. From these optimal decision variables, the platform can deduce (in expectation) the number of control and treatment units of each type matched to supply units of each type:
\begin{equation}
    X_{i,j}^{\tau,\control,\text{CE}}=\frac{D_i^{\tau,\control}}{D_i^{\tau,\ex}}X_{i,j}^{\tau,\ex,\text{CE}}\text{ and } X_{i,j}^{\tau,\treatment,\text{CE}}=\frac{D_i^{\tau,\treatment}}{D_i^{\tau,\ex}}X_{i,j}^{\tau,\ex,\text{CE}}.
\end{equation}
The platform also obtains the optimal dual solution of the cost-excluded matching problem:
\begin{subequations}
    \begin{align}
        \Phi_{\text{CE}} (\bm{d}, \bm{s})= \min \quad & \sum_{i=1}^{n_d} a_i d_i + \sum_{j=1}^{n_s}b_js_j & \\
        \text{s.t.}\quad & a_i + b_j \ge v_{i,j} & \forall i\in[n_d], j\in[n_s]\\
        & a_i \ge 0 & \forall i\in[n_d]\\
        & b_j \ge 0 & \forall j\in[n_s].
    \end{align}
\end{subequations}
We denote by $\bm{A}^{\tau,\text{CE}}$ and $\bm{B}^{\tau,\text{CE}}$ the optimal shadow prices of demand and supply obtained in the experiment state from the optimal cost-excluded matching. Interestingly, while Lemma~\ref{lem:CE-unique} guarantees that the \emph{primal} solution of the cost-excluded matching problem does not change when each weight $v_{i,j}$ is replaced by $v_{i,j}^m$, the same is not true for the \emph{dual} solution. We denote by $\bm{\tilde{A}}^{\tau,m}$ and $\bm{\tilde{B}}^{\tau,m}$ the optimal shadow prices of demand and supply in this altered matching problem, and refer to them as \emph{discounted shadow prices}. We show that they can be computed from the standard cost-excluded shadow prices.

\begin{lemma}\label{lem:fixed-cost-shadow-price}
    The discounted shadow prices $(\bm{\tilde{A}}^{\tau,m},\bm{\tilde{B}}^{\tau,m})$ can be computed as follows.
    \begin{itemize}
        \item Under the proportional-cost model: $(\bm{\tilde{A}}^{\tau,\text{\emph{prop}}},\bm{\tilde{B}}^{\tau,\text{\emph{prop}}})=(1-\alpha)(\bm{A}^{\tau,\text{CE}},\bm{B}^{\tau,\text{CE}})$.
        \item Under the fixed-cost model:
        \[
            (\bm{\tilde{A}}^{\tau,\text{\emph{fixed}}},\bm{\tilde{B}}^{\tau,\text{\emph{fixed}}}) = \begin{cases}(\bm{A}^{\tau,\text{\emph{CE}}} - \kappa\bm{e},\bm{B}^{\tau,\text{\emph{CE}}}), & \text{ if } \sum_{i=1}^{n_d}D_i^{\tau,\expemph} < \sum_{j=1}^{n_s}S_j^{\tau,\bm\gamma},\\
            (\bm{A}^{\tau,\text{\emph{CE}}},\bm{B}^{\tau,\text{\emph{CE}}}-\kappa\bm{e}), & \text{ if } \sum_{i=1}^{n_d}D_i^{\tau,\expemph} \ge \sum_{j=1}^{n_s}S_j^{\tau,\bm\gamma},\end{cases}
        \]
        \item[]where $\bm{e}$ designates the vector of all ones.
    \end{itemize}
\end{lemma}

\subsection{Matching Fluid Limit}

Analyzing our marketplace models in a stochastic setting, particularly with small samples, can be challenging due to the randomness of the Poisson processes generating demand and supply. This challenge is especially pronounced when the low density of demand and supply amplifies or attenuates disproportionately small variations in Poisson arrivals. As these effects diminish with increasing density, we focus on a high-density limit of the system, where $\tau$ tends to infinity.

\begin{proposition}[Fluid limit of matching problems]
    \label{prop:fluid-limit}
    Let the marketplace density $\tau$ tend to infinity. Then the total system value and the optimal matching converge to the optimal objective and solution of a deterministic linear program. More precisely:
    \begin{enumerate}
        \item In the cost-excluded setting, for any $\bm\lambda$ and $\bm\gamma$:
        \[
        \lim_{\tau\to\infty} \E\left[\frac{1}{\tau}\Phi_{\text{\emph{CE}}}(\bm{D}^{\tau,\bm\lambda},\bm{S}^{\tau,\bm\gamma})\right] 
        = \Phi_{\text{\emph{CE}}}(\bm\lambda, \bm\gamma),
        \]
        and the optimal primal solution $X_{i,j}^{\tau, \expemph, \text{\emph{CE}}}$ and dual solution $(A_i^{\tau, \text{\emph{CE}}}$, $B_j^{\tau, \text{\emph{CE}}})$ of $\Phi_{\text{\emph{CE}}}(\frac{1}{\tau}\bm{D}^{\tau,\bm\lambda + \rho \bm{\beta}},\frac{1}{\tau}\bm{S}^{\tau,\bm\gamma})$ converge almost surely to the optimal primal solution $x_{i,j}^{\expemph,\text{\emph{CE}}}$ and dual solution $(a_i^{\text{\emph{CE}}}$, $b_j^{\text{\emph{CE}}})$ of $\Phi_{\text{\emph{CE}}}(\bm\lambda+\rho\bm\beta, \bm\gamma)$.
        
        \item In the cost-included setting, for any $\bm\lambda^\controlemph$, $\bm\lambda^\treatmentemph$ and $\bm\gamma$:
        \[
        \lim_{\tau\to\infty} \E\left[\frac{1}{\tau}\Phi_{\text{\emph{CI}}}^m(\bm{D}^{\tau,\bm\lambda^\controlemph},  \bm{D}^{\tau,\bm\lambda^\treatmentemph},\bm{S}^{\tau,\bm\gamma})\right] 
        = \Phi_{\text{\emph{CI}}}^m( \bm\lambda^\controlemph, \bm\lambda^\treatmentemph, \bm\gamma),
        \]
        and the optimal primal solution $(X_{i,j}^{\tau, \control, \text{\emph{CI}},m}, X_{i,j}^{\tau, \treatmentemph, \text{\emph{CI}}, m})$ and dual solution $(A_i^{\tau, \controlemph, \text{\emph{CI}},m}$,  $A_i^{\tau, \treatmentemph, \text{\emph{CI}},m}$, $B_j^{\tau, \expemph, \text{\emph{CI}},m})$ of $\Phi_{CI}^m(\frac{1}{\tau}\bm{D}^{\tau,  \bm\lambda^\controlemph}, \frac{1}{\tau}\bm{D}^{\tau,\bm\lambda^\treatmentemph} , \frac{1}{\tau}\bm{S}^{\tau,\bm\gamma})$ converge almost surely to the optimal primal solution $(x_{i,j}^{\controlemph, \text{\emph{CI}},m}, x_{i,j}^{\treatmentemph, \text{\emph{CI}},m})$ and dual solution $(a_i^{\controlemph, \text{\emph{CI}},m}$, $a_i^{\treatmentemph, \text{\emph{CI}},m}$ , $b_j^{exp, \text{\emph{CI}},m})$ of $\Phi_{\text{\emph{CI}}}^m( \bm\lambda^\controlemph, \bm\lambda^\treatmentemph, \bm\gamma)$. 
        \end{enumerate}
\end{proposition}

The first part of Proposition~\ref{prop:fluid-limit} restates Theorem 1 of \cite{bright2024reducing}. The second part extends it to the cost-included matching linear program. Proposition~\ref{prop:fluid-limit} allows us to replace a complex stochastic system with a simple deterministic optimization problem in our analysis, simplifying key quantities such as the global treatment effect.

\begin{corollary}
In the fluid limit as $\tau\to\infty$
    \begin{align*}
            \Delta^{m} &= \lim_{\tau\to\infty} \Delta^{\tau, m} = \Phi_{\text{\emph{CI}}}^{m}(\bm{0}, \bm\lambda+\bm\beta, \bm\gamma) -  \Phi_{\text{\emph{CI}}}^{m}(\bm\lambda, \bm{0}, \bm\gamma) = \Phi_{\text{\emph{CE}}}(\bm\lambda+\bm\beta, \bm\gamma) -  \Phi_{\text{\emph{CE}}}(\bm\lambda, \bm\gamma) - c^m(\bm{0}, \bm{\lambda}+\bm{\beta}, \bm{\gamma}).
            \end{align*}
\end{corollary}

For simplicity, we assume that the optimal solutions of $\Phi_{CE}(\bm\lambda + \rho \bm\beta, \bm\gamma)$ and $\Phi_{CI}^m( (1-\rho)\bm\lambda,  \rho (\bm\lambda+\bm\beta), \bm\gamma)$ at the experiment point are unique and nondegenerate. From linear programming sensitivity analysis, we know that the functions $\Phi_{CE}(\cdot,\cdot)$  and $\Phi_{CI}^m(\cdot, \cdot,\cdot)$ are piecewise linear and concave, and thus their derivatives are Riemann-integrable. Their derivatives with respect to the arrival rate of a particular demand type are given by the corresponding demand shadow prices. We can therefore write the global treatment effect as
\begin{equation}
    \label{eq:fundamental-theorem-calculus}
    \Delta^m = \int_0^1  (\bm{a}^{\treatment, \text{CI} , m,\eta} \cdot (\bm\lambda+\bm\beta) - \bm{a}^{\control, \text{CI}, m,\eta} \cdot \bm\lambda)         d\eta
    = \int_0^1  \bm{a}^{\text{CE},\eta} \cdot \bm\beta  \, d\eta - c^m(\bm{0}, \bm{\lambda+\beta}, \bm\gamma).
\end{equation}
\section{Dependence of Interference Bias on the Experimental Design} \label{sec:rct}

In this section, we discuss the standard estimators for the global treatment effect under both the cost-excluded and cost-included design. In particular, we show that interference bias manifests differently in each design.

\subsection{Cost-Excluded Design}

We start by studying the cost-excluded design, where we define the standard ``RCT'' estimator.

\begin{definition}
In cost-excluded experiments, the RCT estimator is defined as follows.
\begin{itemize}
\item Under the proportional-cost model:
    $$
    \hat \Delta_{\text{RCT-CE}}^{\tau, \text{{prop}}}  = \frac{1}{\tau} \left(  \frac{1}{\rho}\sum_{i=1}^{n_d} \sum_{j=1}^{n_s} (1-\alpha) v_{i,j}   X_{i,j}^{\tau, \treatment, \text{CE}}  - \frac{1}{1-\rho} \sum_{i=1}^{n_d} \sum_{j=1}^{n_s} v_{i,j}  X_{i,j}^{\tau, \control, \text{CE}} 
    \right).
    $$
    \item Under the fixed-cost model:
    $$
    \hat \Delta_{\text{RCT-CE}}^{\tau, \text{{fixed}}}  = \frac{1}{\tau} \left(  \frac{1}{\rho}\sum_{i=1}^{n_d} \sum_{j=1}^{n_s}  (v_{i,j}-\kappa)   X_{i,j}^{\tau, \treatment, \text{CE}}  - \frac{1}{1-\rho} \sum_{i=1}^{n_d} \sum_{j=1}^{n_s} v_{i,j}  X_{i,j}^{\tau, \control, \text{CE}} 
    \right).
    $$
\end{itemize}
\end{definition}

In words, the standard approach to estimate the global treatment effect is simply to sum up the total value obtained from the control group, and compare it with the total value obtained from the treatment group (appropriately discounting the value of treated units). We can equivalently write each RCT estimator as follows:
\begin{subequations}
\begin{align}\label{est:RCT-CE}
\hat \Delta_{\text{RCT-CE}}^{\tau, \text{prop}}  &= \frac{1}{\tau} \left(  \frac{1}{\rho}\sum_{i=1}^{n_d} \sum_{j=1}^{n_s}  v_{i,j}   X_{i,j}^{\tau, \treatment, \text{CE}}  - \frac{1}{1-\rho} \sum_{i=1}^{n_d} \sum_{j=1}^{n_s} v_{i,j}  X_{i,j}^{\tau, \control, \text{CE}} 
-  \alpha  \frac{1}{\rho} \sum_{i=1}^{n_d}\sum_{j=1}^{n_s}v_{i,j}   X_{i,j}^{\tau,\treatment,\text{\emph{CE}}}
\right),\\
\hat \Delta_{\text{RCT-CE}}^{\tau, \text{fixed}}  &= \frac{1}{\tau} \left(  \underbrace{\frac{1}{\rho}\sum_{i=1}^{n_d} \sum_{j=1}^{n_s}  v_{i,j}   X_{i,j}^{\tau, \treatment, \text{CE}}  - \frac{1}{1-\rho} \sum_{i=1}^{n_d} \sum_{j=1}^{n_s} v_{i,j}  X_{i,j}^{\tau, \control, \text{CE}}}_{\text{cost-excluded treatment effect estimator}} 
-  \underbrace{\kappa  \frac{1}{\rho} \sum_{i=1}^{n_d}\sum_{j=1}^{n_s}   X_{i,j}^{\tau,\treatment,\text{\emph{CE}}}}_{\text{cost estimator}}
\right),
\end{align}
\end{subequations}
where the first two terms estimate the treatment effect under the assumption that the treatment has zero cost, while the last term estimates the cost of deploying the treatment to the entire platform. In the fluid limit, \cite{bright2024reducing} showed that the first part of this estimator has positive bias, i.e., always overestimates the true global treatment effect. However, their result implies that the cost estimator \emph{also} has positive bias, making the sign of the total cost-excluded RCT estimator \emph{a priori} unclear, since we are subtracting a quantity with positive bias from another quantity with positive bias. We will show in Theorem~\ref{thm:GTE:RCT-CE-overestimates-GTE} that the bias of the first term dominates the bias of the second term, meaning that the RCT estimator is always positively biased in the cost-excluded setting. First, we need to introduce the form of the RCT estimator in the fluid limit.

\begin{lemma}\label{lem:fluid-limit-rct-ce}
In the fluid limit as $\tau\to\infty$,
    \begin{align*}
    \hat{\Delta}_{\text{\emph{RCT-CE}}}^{\text{\emph{prop}} }&= \lim_{\tau\to\infty} \hat\Delta_{\text{\emph{RCT-CE}}}^{\tau, \text{\emph{prop}} }  
                    = (1-\alpha) \bar{\bm{v}}^{\expemph, \text{\emph{CE}}} \cdot (\bm\lambda+\bm\beta) - \bar{\bm{v}}^{\expemph, \text{\emph{CE}}} \cdot \bm\lambda , \\
            \hat{\Delta}_{\text{\emph{RCT-CE}}}^{\text{\emph{fixed}}}&= \lim_{\tau\to\infty} \hat\Delta_{\text{\emph{RCT-CE}}}^{\tau, \text{\emph{fixed}}}  
                    = (\bar{\bm{v}}^{\expemph, \text{\emph{CE}}}-\kappa\bar{\bm{x}}^{\expemph, \text{\emph{CE}}}) \cdot (\bm\lambda+\bm\beta) - \bar{\bm{v}}^{\expemph, \text{\emph{CE}}} \cdot \bm\lambda ,
    \end{align*}
    where $ \bar{v}_i^{\expemph, \text{\emph{CE}}}= \frac{ \sum_{j=1}^{n_s}  v_{i,j} x_{i,j}^{\expemph, \text{\emph{CE}}}  }{\lambda_i + \rho\beta_i} $ is the average value obtained by demand units of type $i$, and $ \bar{x}_i^{\expemph, \text{\emph{CE}}}= \frac{ \sum_{j=1}^{n_s}  x_{i,j}^{\expemph, \text{\emph{CE}}}  }{\lambda_i + \rho\beta_i} $ is the matching rate of demand units of type $i$, in the experiment state.
\end{lemma}

We can now state our key bias result for the RCT estimator in the cost-excluded design.

\begin{theorem} \label{thm:GTE:RCT-CE-overestimates-GTE}
In the cost-excluded experimental setting, the RCT estimator always \textit{overestimates} the global treatment effect, i.e., 
$
\hat\Delta^m_{\text{\emph{RCT-CE}}} \ge \Delta^m.
$
\end{theorem}

\begin{figure}[h!]
\FIGURE{
    \begin{subfigure}[b]{0.33\columnwidth}
    \centering
    \begin{tikzpicture}[scale=1]

  \node[draw,minimum size=0.5cm] (demand) at (0,0) {1};
  \node[draw, circle,minimum size=0.5cm] (supply1) at (2,1.5) {1};
  \node[draw, circle,minimum size=0.5cm] (supply2) at (2,0) {2};
  \node[draw, circle,minimum size=0.5cm] (supply3) at (2,-1.5){3};
  
  \coordinate[label={[label distance=0mm]left:{$\lambda_1=\lambda$}}]() at (0.5,-0.75);
  \coordinate[label={[label distance=0mm]right:{$\gamma_1=1.5$}}]() at (2.4,1.5);
  \coordinate[label={[label distance=0mm]right:{$\gamma_2=2$}}]() at (2.4,0);
  \coordinate[label={[label distance=0mm]right:{$\gamma_3=2$}}]() at (2.4,-1.5);

  \draw[] (demand) -- node[midway, above, sloped] {\footnotesize 2} (supply1);
  \draw[] (demand) -- node[midway, above, sloped] {\footnotesize 1} (supply2);
  \draw[] (demand) -- node[midway, above, sloped] {\footnotesize 0.25} (supply3);

\end{tikzpicture}
    \vspace{40pt}
    \caption{Simple matching instance}
    \end{subfigure}
    \begin{subfigure}[b]{0.66\columnwidth}
    \centering
    \begin{tikzpicture}
    \begin{axis}[
        width=11cm,
        height=8cm,
        xlabel={$\lambda$},
        ylabel={$\Phi$},
        grid=major,
        axis lines=left,
        axis line style={thick},
        xtick={0,1.5,3.5,5.5,7},
        ytick={0,3,5,5.5},
        every axis x label/.style={at={(ticklabel* cs:1)},anchor=north},
        every axis y label/.style={at={(ticklabel* cs:1)},anchor=east},
        enlargelimits,
    ]
        \addplot[
            color=RoyalBlue,
            very thick,
            domain=0:7,
            samples=100
        ]
        coordinates {
            (0, 0)
            (1.5, 3)
            (3.5, 5)
            (5.5, 5.5)
            (7, 5.5)
        };
        \node[anchor=south] at (axis cs:6,5.5) {\color{RoyalBlue}$\Phi_{\text{CE}}(\lambda,\bm\gamma)$};
        \addplot[
            color=Orange,
            very thick,
            domain=0:7,
            samples=100,
            densely dashed
        ]
        coordinates {
            (0, 0)
            (1.5, 2.55)
            (3.5, 4.25)
            (5.5, 4.675)
            (7, 4.675)
        };
        \node[anchor=north] at (axis cs:6,4.2) {\color{Orange}$(1-\alpha)\Phi_{\text{CE}}(\lambda,\bm\gamma)$};
    \end{axis}
\end{tikzpicture}
    \caption{Total value as a function of demand}
    \end{subfigure}
    }
    {Motivating example: matching with one demand type.\label{fig:setup}}
    {The left panel shows the simple matching instance with the arrival rates and the (control) matching values for each supply type. The right panel shows the total matching value as a function of the demand arrival rate, for control units (solid line) and treatment units (dashed line).}
\end{figure}

Theorem~\ref{thm:GTE:RCT-CE-overestimates-GTE} admits an intuitive graphical explanation. To illustrate it, we introduce a simple matching instance, with a single demand type with arrival rate $\lambda$, and three supply types with arrival rates $\gamma_1=1.5$ and $\gamma_2=\gamma_3=2$. We will consider a treatment with proportional cost, parametrized by $\alpha=0.15$, to provide intuition, but all the results also hold for the fixed-cost setting. Figure~\ref{fig:setup} shows a diagram of the matching instance, along with a plot of the value function (undiscounted for control and discounted for treatment) as a function of the demand arrival rate.

We assume that the arrival rate of demand units is given by $\lambda^{\control}=1$ under control, and $\lambda^{\treatment}=4$ under treatment. We assume that $\rho=0.5$: in the experiment state, the total demand arrival rate is thus given by $\rho\lambda^{\treatment} + (1-\rho)\lambda^{\control}=2.5$. From Proposition~\ref{prop:fluid-limit}, we see that the RCT estimator relies on the average value obtained from each type in the experiment state. Implicitly, the RCT estimator assumes that the average value obtained from each demand type is the same in global control and in global treatment as it is in the intermediate experiment state --- a false assumption in the presence of marketplace interference. We can show how this assumption breaks in Figure~\ref{fig:rct-ce-breaks}. Extrapolating the value function from the experiment state assuming the average value from each demand unit remains constant ignores the diminishing returns of additional demand induced by the limited availability of supply. As a result, the linear approximation for the value under global control provides an under-estimate, while the linear approximation for the value under global treatment provides an over-estimate --- so their difference overestimates the global treatment effect.

\begin{figure}
    \FIGURE{
        \begin{tikzpicture}
    \begin{axis}[
        width=11cm,
        height=8cm,
        xlabel={$\lambda$},
        ylabel={$\Phi$},
        grid=major,
        axis lines=left,
        axis line style={thick},
        xtick={0,1,2.5,4},
        xticklabels={0,$\lambda^{\text{con}}=1$, $\lambda^{\text{exp}}=2.5$, $\lambda^{\text{tre}}=4$},
        ytick={0,1.6,2,3.4,4,4.36,5.44},
        yticklabels={0,
        \raisebox{-1ex}{\color{BlueViolet}$\hat{\Phi}_{\text{RCT-CE}}(\lambda^{\text{con}})$},
        {\color{RoyalBlue}$\Phi_{\text{CE}}(\lambda^{\text{con}})$},
        {\color{Orange}$(1-\alpha)\Phi_{\text{CE}}(\lambda^{\text{exp}})$},
        \raisebox{-1ex}{\color{RoyalBlue}$\Phi_{\text{CE}}(\lambda^{\text{exp}})$},
        \raisebox{2ex}{\color{Orange}$(1-\alpha)\Phi_{\text{CE}}(\lambda^{\text{tre}})$},
        {\color{RawSienna}$(1-\alpha)\hat{\Phi}_{\text{RCT-CE}}(\lambda^{\text{tre}})$}
        },
        every tick label/.append style={font=\footnotesize},
        every axis x label/.style={at={(ticklabel* cs:1)},anchor=north},
        every axis y label/.style={at={(ticklabel* cs:1)},anchor=east},
        enlargelimits,
    ]
    \addplot[
            color=RoyalBlue,
            very thick,
            domain=0:7,
            samples=100
        ]
        coordinates {
            (0, 0)
            (1.5, 3)
            (3.5, 5)
            (5.5, 5.5)
            (7, 5.5)
        };
        \node[anchor=south] at (axis cs:6,4.8) {\footnotesize\color{RoyalBlue}$\Phi_{\text{CE}}(\lambda)$};
    \addplot[
            color=Orange,
            very thick,
            domain=0:7,
            samples=100,
            densely dashed
        ]
        coordinates {
            (0, 0)
            (1.5, 2.55)
            (3.5, 4.25)
            (5.5, 4.675)
            (7, 4.675)
        };
        \node[anchor=north] at (axis cs:6,4.6) {\footnotesize\color{Orange}$(1-\alpha)\Phi_{\text{CE}}(\lambda)$};
        \addplot[
            color=RawSienna,
            very thick,
            domain=0:7,
            samples=100,
            densely dotted
        ]
        coordinates {
            (0, 0)
            (4.1, 5.576)
        };
        \node[anchor=north] at (axis cs:5.2,6.4) {\footnotesize\color{RawSienna}$(1-\alpha)\hat{\Phi}_{\text{RCT-CE}}(\lambda)
        $};
        \addplot[
            color=BlueViolet,
            very thick,
            domain=0:7,
            samples=100,
            densely dashdotted
        ]
        coordinates {
            (0, 0)
            (4.1, 6.56)
        };
        \node[anchor=north] at (axis cs:4.5,7.4) {\footnotesize\color{BlueViolet}$\hat{\Phi}_{\text{RCT-CE}}(\lambda)=\bar{v}^{\text{exp},\text{CE}}\cdot\lambda
        $};
        \draw [<->,thick] (axis cs:0,2) -- (axis cs:0,4.36) node[midway, left] {\footnotesize$\Delta$};
        \draw [<->,thick] (axis cs:0.2,1.6) -- (axis cs:0.2,5.44) node[midway, right] {\footnotesize$\hat{\Delta}_{\text{RCT-CE}}$};
    \end{axis}
\end{tikzpicture}
    }
    {Bias of the cost-excluded RCT estimator.\label{fig:rct-ce-breaks}}
    {In the cost-excluded setting, the RCT estimator relies on a linear approximation of the value function under control (dash-dotted line) and the discounted value function under treatment (dotted line). This linear approximation induces bias. For clarity we drop the final argument of $\Phi$ since it is always $\bm\gamma$.}
\end{figure} 

Our proof of Theorem~\ref{thm:GTE:RCT-CE-overestimates-GTE} overcomes the fact that the RCT estimator is the difference between two positively-biased quantities, by observing that even though the total cost of deploying treatment is overestimated, the RCT estimator still overestimates the total discounted system value in global treatment. Thus, if we run an experiment ignoring treatment costs in the matching before factoring them in later, the standard estimator will systematically overestimate the true effect.

\subsection{Cost-Included Design}

We now consider the standard estimator in the cost-included design, defined as follows.

\begin{definition}
For cost-included experiments, the RCT estimator is given by
\begin{equation*}\label{est:RCT-CI}
\hat \Delta_{\text{RCT-CI}}^{\tau, m}  = 
\frac{1}{\tau} 
\left(  
\frac{1}{\rho}\sum_{i=1}^{n_d} \sum_{j=1}^{n_s}  v_{i,j}^{\treatment , m}   X_{i,j}^{\tau, \treatment, \text{CI},m}  - \frac{1}{1-\rho} \sum_{i=1}^{n_d} \sum_{j=1}^{n_s} v_{i,j}^{\control}  X_{i,j}^{\tau, \control, \text{CI},m}  
\right).
\end{equation*}
\end{definition}

As in the cost-excluded setting, the standard estimator simply compares the total value obtained from the treatment group with the total value obtained from the control group. Unlike the cost-excluded setting, the treatment values are already discounted in the matching formulation. The same pattern holds in the fluid limit, as we see in the following lemma.

\begin{lemma}
    \label{lem:fluid-limit-rct-ci}
    In the fluid limit as $\tau\to\infty$,
    \[
    \hat{\Delta}_{\text{\emph{RCT-CI}}}^{m} = \lim_{\tau\to\infty} \hat\Delta_{\text{\emph{RCT-CI}}}^{\tau, m} 
                    = \bar{\bm{v}}^{\treatmentemph, \text{\emph{CI}},m} \cdot (\bm\lambda+\bm\beta) - \bar{\bm{v}}^{\controlemph, \text{\emph{CI}},m} \cdot \bm\lambda ,
    \]
    where $ \bar{v}_i^{\treatmentemph, \text{\emph{CI}},m}= \frac{ \sum_{j=1}^{n_s} v_{i,j}^{\treatmentemph, m} x_{i,j}^{\treatmentemph, \text{\emph{CI}},m}  }{\rho(\lambda_i + \beta_i)} $ is the average value obtained by treated demand units of type $i$ and $ \bar{v}_i^{\controlemph, \text{\emph{CI}},m}= \frac{ \sum_{j=1}^{n_s}v^\controlemph_{i,j} x_{i,j}^{\controlemph, \text{\emph{CI}},m}  }{(1-\rho)\lambda_i } $ is the average value obtained by untreated demand units of type $i$ in the cost-inclusive experiment.
\end{lemma}

Superficially, the cost-included RCT estimator looks a lot like the cost-excluded RCT estimator. It turns out, however, that its behavior is much more complex under interference. In particular, it depends on the overall ratio between supply and demand in the marketplace. We define the total marketplace supply as the sum of all supply arrival rates, i.e., $\Gamma=\sum_{j=1}^{n_s}\gamma_j$. For a particular set of arrival rates $(\bm\lambda,\bm\beta,\bm\gamma)$, we can use $\Gamma$ as a scaling factor --- when we do this, we always assume that the supply arrival rates for each type are all scaled by the same amount.
The following theorem then describes how the cost-included RCT estimator performs under various supply configurations.

\begin{theorem}\label{thm:GTE:Bias-RCT-CI-can-be-high}
In the cost-included setting, for any arrival rates $(\bm\lambda,\bm\beta,\bm\gamma)$ and any cost model $m$, the RCT estimator has the following properties:
\begin{enumerate}
\item It is negatively biased in the low-supply limit, i.e., there exists $\Gamma_{\min}$ such that for any $\Gamma < \Gamma_{\min}$,
\[
\hat{\Delta}^m_{\text{\emph{RCT-CI}}} - \Delta^m < 0.
\]
\item It can also be positively biased, i.e., there always exists a scaling factor $\Gamma_0$ such that for $\Gamma = \Gamma_0$,
\[
\hat{\Delta}^m_{\text{\emph{RCT-CI}}} - \Delta^m > 0.
\]
\item Its relative bias in the low-supply limit tends to infinity as the treatment cost tends to zero, i.e., if $\Gamma < \Gamma_{\min}$,
\begin{align*}
\lim_{\alpha \to 0} \frac{\hat{\Delta}^{\text{\emph{prop}}}_{\text{\emph{RCT-CI}}}-\Delta^{\text{\emph{prop}}}}{\abs{\Delta^{\text{\emph{prop}}}}} 
= 
\lim_{\kappa \to 0} \frac{\hat{\Delta}^{\text{\emph{fixed}}}_{\text{\emph{RCT-CI}}}-\Delta^{\text{\emph{fixed}}}}{\abs{\Delta^{\text{\emph{fixed}}}}} 
= 
-\infty.
\end{align*}
\end{enumerate}
\end{theorem}

\begin{figure}[h!]
    \caption{Bias of the RCT-CI estimator on the simple example from Figure~\ref{fig:setup}.}
    \begin{tikzpicture}
    \begin{axis}[
        clip=false,
        width=7cm,
        height=7cm,
        xlabel={$\eta$},
        ylabel={$\Phi$},
        grid=major,
        axis lines=left,
        legend style={cells={align=left}},
        legend cell align={left},
        legend pos=outer north east,
        axis line style={thick},
        xtick={0,0.5,1},
        xticklabels={0,$\rho=0.5$,1},
        ytick={0,2,4.36,5.1},
        yticklabels={
            0,
            {\color{RoyalBlue}$\Phi_{\text{CI}}(\lambda^{\text{con}}{,}\,0)$},
            {\color{Orange}$\Phi_{\text{CI}}(0{,}\,\lambda^{\text{tre}})$},
            {\color{RawSienna}$\hat{\Phi}_{\text{RCT-CI}}(0{,}\,\lambda^{\text{tre}})$}
        },
        every tick label/.append style={font=\footnotesize},
        every axis x label/.style={at={(ticklabel* cs:1)},anchor=north},
        every axis y label/.style={at={(ticklabel* cs:1)},anchor=east},
        enlargelimits,
    ]
        \addplot[
            color=RoyalBlue,
            very thick,
            domain=0:1,
            samples=100
        ]
        coordinates {
            (0, 2)
            (1, 0)
        };
        \addplot[
            color=Orange,
            very thick,
            domain=0:1,
            samples=100,
            densely dashed
        ]
        coordinates {
            (0, 0)
            (0.17, 1.13)
            (0.83, 3.97)
            (1, 4.36)
        };
        \addplot[
            color=ForestGreen,
            very thick,
            domain=0:1,
            samples=100,
            dashdotted
        ]
        coordinates {
            (0, 2)
            (0.17, 2.79)
            (0.83, 4.31)
            (1, 4.36)
        };
        \addplot[
            color=RawSienna,
            very thick,
            domain=0:1,
            samples=100,
            densely dotted
        ]
        coordinates {
            (0, 0)
            (1, 5.1)
        };
        \draw [<->,thick] (axis cs:-0.02,2) -- (axis cs:-0.02,4.36) node[midway, right,xshift=-0.15cm] {\footnotesize$\Delta$};
        \draw [<->,thick] (axis cs:1.02,2) -- (axis cs:1.02,5.1) node[midway, left,yshift=-0.5cm] {\footnotesize$\hat{\Delta}_{\text{RCT-CI}}$};
        \node[anchor=east] at (-0.11,1.6) {\color{BlueViolet}\scriptsize$=\hat{\Phi}_{\text{RCT-CI}}(\lambda^{\text{con}}{,}\,0)$};
        \node[anchor=center] at (0.5,-2) {(a) $\lambda^{\text{con}} = 1,~\lambda^{\text{tre}}=4$};
    \end{axis}
\end{tikzpicture}%
\begin{tikzpicture}
    \begin{axis}[
        clip=false,
        width=7cm,
        height=7cm,
        xlabel={$\eta$},
        ylabel={$\Phi$},
        grid=major,
        axis lines=left,
        legend style={cells={align=left}},
        legend cell align={left},
        legend pos=outer north east,
        axis line style={thick},
        xtick={0,0.5,1},
        xticklabels={0,$\rho=0.5$,1},
        ytick={0,3.6125,4.5,4.5675,6},
        yticklabels={
            0,
            {\color{RawSienna}$\hat{\Phi}_{\text{RCT-CI}}(0{,}\,\lambda^{\text{tre}})$},
            \raisebox{-3ex}{\color{RoyalBlue}$\Phi_{\text{CI}}(\lambda^{\text{con}}{,}\,0)$},
            \raisebox{3ex}{\color{Orange}$\Phi_{\text{CI}}(0{,}\,\lambda^{\text{tre}})$},
            {\color{BlueViolet}$\hat{\Phi}_{\text{RCT-CI}}(\lambda^{\text{con}}{,}\,0)$}
        },
        every tick label/.append style={font=\footnotesize},
        every axis x label/.style={at={(ticklabel* cs:1)},anchor=north},
        every axis y label/.style={at={(ticklabel* cs:1)},anchor=east},
        enlargelimits,
    ]
        \addplot[
            color=RoyalBlue,
            very thick,
            domain=0:1,
            samples=100
        ]
        coordinates {
            (0, 4.5)
            (0.5, 3)
            (1, 0)
        }; \label{plot1}
        \addplot[
            color=Orange,
            very thick,
            domain=0:1,
            samples=100,
            densely dashed
        ]
        coordinates {
            (0, 0)
            (0.25, 1.0625)
            (0.5, 1.80625)
            (1, 4.56875)
        };\label{plot2}
        \addplot[
            color=ForestGreen,
            very thick,
            domain=0:1,
            samples=100,
            dashdotted
        ]
        coordinates {
            (0, 4.5)
            (0.25, 4.8125)
            (0.5, 4.80625)
            (1, 4.5675)
        };\label{plot3}
        \addplot[
            color=RawSienna,
            very thick,
            domain=0:1,
            samples=100,
            densely dotted
        ]
        coordinates {
            (0, 0)
            (1, 3.6125)
        };\label{plot4}
        \addplot[
            color=BlueViolet,
            very thick,
            domain=0:1,
            samples=100,
            densely dashdotted
        ]
        coordinates {
            (0, 6)
            (1, 0)
        };\label{plot5}
        \draw [|-|,semithick] (axis cs:-0.02,4.5) -- (axis cs:-0.02,4.5675) node[midway, left,xshift=0.05cm] {\footnotesize$\Delta$};
        \draw [<->,thick] (axis cs:1.02,6) -- (axis cs:1.02,3.6125) node[midway, left,yshift=0.5cm] {\footnotesize$\hat{\Delta}_{\text{RCT-CI}}$};
        \node[anchor=center] at (0.5,-2.35) {(b) $\lambda^{\text{con}} = 3,~\lambda^{\text{tre}}=5$};
    \end{axis}
\end{tikzpicture}
\vspace{-10pt}
\begin{center}
\begin{tikzpicture}
\matrix(first)[
  matrix of nodes,
  execute at empty cell={\node[draw=none]{};},
  anchor=north,
  inner sep=0.2em,
  name=table,
  font=\scriptsize,
  nodes={anchor=west,align=left}
]at(0,0)
{  
    \ref{plot1} & $\sum_{j=1}^{3}v^{\text{con}}_{1,j}x^{\text{con}}_{1,j}(\eta)$ & 
    \ref{plot2} & $\sum_{j=1}^{3}v^{\text{tre}}_{1,j}x^{\text{tre}}_{1,j}(\eta)$ &
    \ref{plot3} & $\Phi_{\text{CI}}((1-\eta)\lambda^{\text{con}},\eta\lambda^{\text{tre}})$
    \\
    & (value of control units) & & (value of treatment units) & & (total value)\\
};
\matrix(first)[
  matrix of nodes,
  execute at empty cell={\node[draw=none]{};},
  anchor=north,
  inner sep=0.2em,
  name=table,
  font=\scriptsize,
  nodes={anchor=west,align=left}
]at(0,-1.1)
{  
    \ref{plot5} & $\bar{v}_1^{\text{con}}(\rho)\cdot(1-\eta)\lambda^{\text{con}}$ &
    \ref{plot4} & $\bar{v}_1^{\text{tre}}(\rho)\cdot\eta\lambda^{\text{tre}}$ 
    \\
    & (estimated value of control units) & & (estimated value of treatment units)\\
};
\end{tikzpicture}
\end{center}
    \label{fig:rct-ci}
    {\footnotesize\textit{Note.} The left panel shows an instance of positive bias. Linearly extrapolating the average value obtained from the control and treatment groups at the experiment state correctly estimates the value at global control but overestimates the value at global treatment. The right panel shows an instance of negative bias: linearly extrapolating the average value obtained from the control and treatment groups at the experiment state both overestimates the value at global control and underestimates the value at global treatment. Note that the horizontal axis now represents the continuous spectrum of possible experiments, parametrized by $\eta$, where $\eta=0$ designates global control, $\eta=1$ global treatment, and $\eta=\rho$ the actual experiment under study. For clarity we drop the final argument of $\Phi$ since it is always $\bm\gamma$.}
\end{figure}
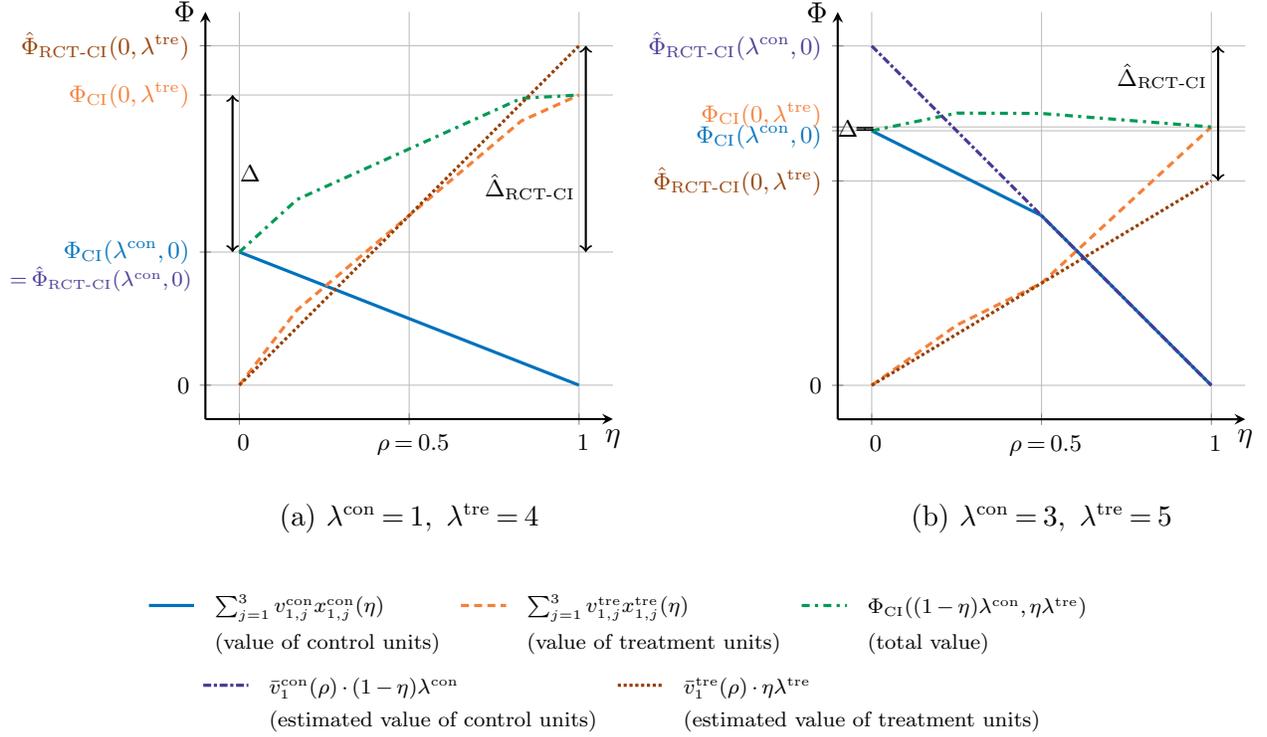

We can illustrate the first two properties of Theorem~\ref{thm:GTE:Bias-RCT-CI-can-be-high} using the example from Figure~\ref{fig:setup}, with two different levels of demand. If we keep $\lambda^{\control}=1$ and $\lambda^{\treatment}=4$ as in Figure~\ref{fig:rct-ce-breaks}, the true treatment effect is $\Delta=2.36$, but the cost-included RCT estimator gives a value of $\hat{\Delta}_{\text{RCT-CI}}=3.1$. However, if $\lambda^{\control}=3$ and $\lambda^{\treatment}=5$ (a higher aggregate level of demand is equivalent to a lower aggregate level of supply), then the true treatment effect is $\Delta=0.06875$, but the estimator gives a value of $\hat{\Delta}_{\text{RCT-CI}}=-2.3875$. We illustrate the workings of the cost-included estimator on these examples in Figure~\ref{fig:rct-ci}.
Even in our simple example, the matching value function depends on the number of control and treatment demand units, which means that even without considering supply, its input space is two-dimensional. For ease of visualization, we consider the one-dimensional line between the point $(\lambda^{\text{con}},0)$ and  $(0,\lambda^{\text{tre}})$, parametrized by $\eta\in[0,1]$. The global control state $(\eta=0)$, global treatment state $(\eta=1)$ and experiment state $(\eta=\rho)$ all lie along this line. The RCT estimator uses information about the experiment state (in particular, the average value obtained from each type) to estimate the value at the endpoints of this line. Because of interference, the matching value function is piecewise linear and concave. The value obtained from each group is also piecewise linear, but not necessarily concave. As a result, unlike the cost-excluded setting, in the cost-included setting the bias can be either positive or negative.

In the cost-excluded setting, marketplace interference happens when demand units contend for scarce supply: not taking into account this contention means overestimating the marginal value of additional demand units, which leads to positive bias in the RCT estimator. In the cost-included setting, however, a second effect materializes: because treated units are costly, they are deprioritized in the matching algorithm. This leads to overemphasizing the value of control units, and therefore, negative bias. The third part of Theorem~\ref{thm:GTE:Bias-RCT-CI-can-be-high} highlights this second effect. Surprisingly, the effect appears for any nonzero cost, leading to a discontinuity. Even an infinitesimal cost leads to deprioritization of treatment units, which makes them look significantly worse than control units to the matching algorithm. We would expect that a very small cost also implies a small difference between the cost-excluded and cost-included setting. In fact, the opposite is true: the difference between the two settings is largest when the cost of treatment is small.

One interesting takeaway from this section is that, if we use the standard estimator, it is not clear whether the cost-excluded or cost-included design is preferable. The cost-excluded design always has positive bias, while the cost-included design could have negative or positive (or zero) bias. Neither is an especially appealing option. We investigate alternatives in the following section.
\section{Bias Reduction Via Shadow Prices} \label{sec:shadow-prices}

The examples in Figures~\ref{fig:rct-ce-breaks} and \ref{fig:rct-ci} suggest a key reason for the breakdown of the standard estimator in both experimental designs: the assumption that the average value obtained by a demand unit is the same under global control and under global treatment (after accounting for price discounts). This assumption does not hold under marketplace interference. We may therefore prefer another way to estimate the total system value under global control and under global treatment --- by using the optimal dual solution in the experiment state.

\subsection{Cost-Excluded Design}

Recall that in the cost-excluded experiment state, we obtain \emph{two} sets of optimal shadow prices, depending on whether we use the original (control) or discounted (treatment) matching weights. Though this choice does not affect the primal solution (this is a \emph{cost-excluded} matching, after all), it does affect the dual solution. We use both the original shadow prices and the discounted shadow prices in defining the following cost-excluded shadow price estimator.

\begin{definition}
For cost-excluded experiments, the shadow price (SP) estimator is defined as
\[
\hat{\Delta}^{\tau,m}_{\text{SP-CE}} = \frac{1}{\tau}\left[\left(\frac{1}{\rho}\bm{\tilde{A}}^{\tau,m}\cdot\bm{D}^{\tau,\treatment} + \bm{\tilde{B}}^{\tau,m}\cdot \bm{S}^{\tau,\bm\gamma}\right) - \left(\frac{1}{1-\rho} \bm{A}^{\tau,\text{CE}}\cdot \bm{D}^{\tau,\control} + \bm{B}^{\tau,\text{CE}}\cdot\bm{S}^{\tau,\bm\gamma} \right)\right].
\]
\end{definition}

The key idea behind the shadow price estimator is that it assumes that while the optimal matching will obviously be different under global control and/or global treatment (because there will be different numbers of units), the \emph{dual} of the optimal matching may not change much, especially since the treatment effect in any experiment is likely to be small relative to the total system value. We can use Lemma~\ref{lem:fixed-cost-shadow-price} to obtain more specific expressions for each cost model. In particular, under the proportional-cost model, we obtain
\begin{equation}
    \hat{\Delta}^{\tau,\text{prop}}_{\text{SP-CE}} = \frac{1}{\tau}\left[\bm{A}^{\tau,\text{CE}}\cdot\left(\frac{1}{\rho}\bm{D}^{\tau,\treatment} - \frac{1}{1-\rho}\bm{D}^{\tau,\control}\right) - \alpha\left(\frac{1}{\rho}\bm{A}^{\tau,\text{CE}}\cdot \bm{D}^{\tau,\treatment} + \bm{B}^{\tau,\text{CE}}\cdot\bm{S}^{\tau,\bm\gamma}\right)\right],
\end{equation}
where the first term is exactly the original shadow price estimator from \cite{bright2024reducing}, while the second term is an estimate of the total cost of the intervention. Similarly, under the fixed-cost model, we obtain
\begin{equation}
    \hat{\Delta}^{\tau,\text{fixed}}_{\text{SP-CE}} = \begin{cases}
        \frac{1}{\tau}\left[\bm{A}^{\tau,\text{CE}}\cdot\left(\frac{1}{\rho}\bm{D}^{\tau,\treatment} - \frac{1}{1-\rho}\bm{D}^{\tau,\control}\right)- \kappa\frac{1}{\rho}\sum_{i=1}^{n_d}D_i^{\tau,\treatment}\right]& \text{if } \sum_{i=1}^{n_d}D_i^{\tau,\ex} < \sum_{j=1}^{n_s}S_j^{\tau,\bm\gamma},\\
        \frac{1}{\tau}\left[\bm{A}^{\tau,\text{CE}}\cdot\left(\frac{1}{\rho}\bm{D}^{\tau,\treatment} - \frac{1}{1-\rho}\bm{D}^{\tau,\control}\right)- \kappa\sum_{j=1}^{n_s}S_j^{\tau,\bm\gamma}\right]& \text{if } \sum_{i=1}^{n_d}D_i^{\tau,\ex} \ge \sum_{j=1}^{n_s}S_j^{\tau,\bm\gamma}.
    \end{cases}
\end{equation}
Once again, the second term estimates the total cost of the intervention, which is driven by total demand if the experiment is demand-constrained (first case) or total supply if the experiment is supply-constrained (second case).
In order to analyze whether this new estimator is biased, we again turn to the fluid limit of the matching system.

\begin{lemma} \label{lem:sp-ce-fluid-limit}
    In the cost-excluded design, the fluid-limit form of the shadow price estimator is
    \[
        \hat{\Delta}^{m}_{\text{\emph{SP-CE}}} = \lim_{\tau\to\infty}\hat{\Delta}^{\tau,m}_{\text{\emph{SP-CE}}}=\left(\bm{\tilde{a}}^{m}\cdot(\bm\lambda + \bm\beta) + \bm{\tilde{b}}^{m}\cdot \bm\gamma\right) - \left(\bm{a}^{\text{\emph{CE}}}\cdot \bm\lambda + \bm{b}^{\text{\emph{CE}}}\cdot\bm\gamma \right),
    \]
    where we define the discounted shadow prices $(\bm{\tilde{a}}^{m}, \bm{\tilde{b}}^{m})$ of the fluid-limit matching $\Phi_{\text{\emph{CE}}}(\bm{\lambda}+\rho\bm\beta, \bm\gamma)$.
\end{lemma}

In the fluid limit, the discounted shadow prices $(\bm{\tilde{a}}^{m}, \bm{\tilde{b}}^{m})$ can be computed from the original shadow prices following the same reasoning from Lemma~\ref{lem:fixed-cost-shadow-price}. We can further specify the form of the fluid-limit shadow price estimator as follows.

\begin{corollary}
    The cost-excluded shadow price estimator can equivalently be written as
    \begin{align*}
        \hat{\Delta}_{\text{\emph{SP-CE}}}^{\text{\emph{prop}} }&=\bm{a}^{\text{\emph{CE}}} \cdot \bm\beta - \alpha \, (   \bm{a}^{\text{\emph{CE}}}\cdot ( \bm\lambda+\bm\beta) +\bm{b}^{\text{\emph{CE}}}\cdot \bm\gamma     ), \\
        \hat{\Delta}_{\text{\emph{SP-CE}}}^{\text{\emph{fixed}} }&= \begin{cases}
            \bm{a}^{\text{\emph{CE}}} \cdot \bm\beta - \kappa\sum_{i=1}^{n_d}\lambda_i + \beta_i, & \text{if } \sum_{i=1}^{n_d}\lambda_i + \rho\beta_i < \sum_{j=1}^{n_s}\gamma_j,\\
            \bm{a}^{\text{\emph{CE}}} \cdot \bm\beta - \kappa\sum_{j=1}^{n_s}\gamma_j, & \text{if } \sum_{i=1}^{n_d}\lambda_i + \rho\beta_i \ge \sum_{j=1}^{n_s}\gamma_j.\\
        \end{cases}
    \end{align*}
\end{corollary}

Our next result shows that the shadow price estimator in the cost-excluded design reduces bias compared to the standard estimator.

\begin{theorem}\label{thm:RCT-CE-vs-SP-CE}
In the cost-excluded setting, under cost model $m$, the SP estimator reduces bias from the RCT estimator as long as $\rho \le \frac{1-\zeta^m}{2-\zeta^m}$, where
\begin{align*}
\zeta^{\text{\emph{prop}}} &=
    \alpha,\\
\zeta^{\text{\emph{fixed}}} &=
    \frac{\kappa}{\min_{i,j}v_{i,j}}.
\end{align*}
Furthermore, for any $\alpha$ or $\kappa$, there exists an instance for which the bound is tight.
\end{theorem}

Theorem~\ref{thm:RCT-CE-vs-SP-CE} provides a clear guarantee for bias reduction, which only depends on known quantities (the discounts $\alpha$ or $\kappa$) and a controllable quantity (the treatment fraction $\rho$). The condition is not very restrictive: as long as the discount remains small, the upper bound on the treatment fraction is close to $0.5$. In addition, requiring a small $\rho$ is not very limiting in practice since platforms often prefer experiments with smaller treatment groups (so as to limit experimentation costs).

\begin{figure}
    \FIGURE{
            \begin{tikzpicture}
        \begin{axis}[
            width=11cm,
            height=8cm,
            xlabel={$\lambda$},
            ylabel={$\Phi$},
            grid=major,
            axis lines=left,
            axis line style={thick},
            xtick={0,1,2.5,4},
            xticklabels={0,$\lambda^{\text{con}}=1$, $\lambda^{\text{exp}}=2.5$, $\lambda^{\text{tre}}=4$},
            ytick={0,2.5,2,3.4,4,4.36,4.675},
            yticklabels={0,
            \raisebox{-1ex}{\color{BlueViolet}$\hat{\Phi}_{\text{SP-CE}}(\lambda^{\text{con}})$},
            {\color{RoyalBlue}$\Phi_{\text{CE}}(\lambda^{\text{con}})$},
            {\color{Orange}$(1-\alpha)\Phi_{\text{CE}}(\lambda^{\text{exp}})$},
            \raisebox{-1ex}{\color{RoyalBlue}$\Phi_{\text{CE}}(\lambda^{\text{exp}})$},
            \raisebox{2ex}{\color{Orange}$(1-\alpha)\Phi_{\text{CE}}(\lambda^{\text{tre}})$},
            \raisebox{3.5ex}{\phantom{a}\color{RawSienna}$(1-\alpha)\hat{\Phi}_{\text{SP-CE}}(\lambda^{\text{tre}})$}
            },
            every tick label/.append style={font=\footnotesize},
            every axis x label/.style={at={(ticklabel* cs:1)},anchor=north},
            every axis y label/.style={at={(ticklabel* cs:1)},anchor=east},
            enlargelimits,
        ]
        \addplot[
                color=RoyalBlue,
                very thick,
                domain=0:7,
                samples=100
            ]
            coordinates {
                (0, 0)
                (1.5, 3)
                (3.5, 5)
                (5.5, 5.5)
                (7, 5.5)
            };
            \node[anchor=south] at (axis cs:6,5.4) {\footnotesize\color{RoyalBlue}$\Phi_{\text{CE}}(\lambda)$};
        \addplot[
                color=Orange,
                very thick,
                domain=0:7,
                samples=100,
                densely dashed
            ]
            coordinates {
                (0, 0)
                (1.5, 2.55)
                (3.5, 4.25)
                (5.5, 4.675)
                (7, 4.675)
            };
            \node[anchor=north] at (axis cs:6,4.6) {\footnotesize\color{Orange}$(1-\alpha)\Phi_{\text{CE}}(\lambda)$};
            \addplot[
                color=RawSienna,
                very thick,
                domain=0:7,
                samples=100,
                densely dotted
            ]
            coordinates {
                (0.8, 1.955)
                (4.1, 4.76)
            };
            \node[anchor=north] at (axis cs:5.2,5.4) {\footnotesize\color{RawSienna}$(1-\alpha)\hat{\Phi}_{\text{SP-CE}}(\lambda)
            $};
            \addplot[
                color=BlueViolet,
                very thick,
                domain=0:7,
                samples=100,
                densely dashdotted
            ]
            coordinates {
                (0.5, 2)
                (4.1, 5.6)
            };
            \node[anchor=north] at (axis cs:3.5,6.3) {\footnotesize\color{BlueViolet}$\hat{\Phi}_{\text{SP-CE}}(\lambda)=a^{\text{exp},\text{CE}}\cdot\lambda + C
            $};
            \draw [<->,thick] (axis cs:0,2) -- (axis cs:0,4.36) node[midway, left] {\footnotesize$\Delta$};
            \draw [<->,thick] (axis cs:0.2,2.5) -- (axis cs:0.2,4.675) node[midway, right] {\footnotesize$\hat{\Delta}_{\text{SP-CE}}$};
            \addplot[
                color=White,
                thick,
                domain=0:7,
                samples=100,
                densely dotted
            ]
            coordinates {
                (4.05, 6.56)
                (4.1, 6.56)
            };
        \end{axis}
    \end{tikzpicture}
    }
    {Bias of the cost-excluded shadow price estimator.\label{fig:sp-ce}}
    {In the cost-excluded setting, we still build linear approximations of the discounted and undiscounted value functions, but now the linear approximation is correct in a Taylor sense.}
\end{figure} 

What Theorem~\ref{thm:RCT-CE-vs-SP-CE} does not provide is an indication of \emph{how much} the SP estimator reduces bias. Quantifying bias reduction in general requires some additional assumptions about the state of the marketplace --- while these assumptions may not be directly verifiable from the single experiment state, the platform may be able to assess their validity using general knowledge about the marketplace. We lay out the result in the following theorem.

\begin{theorem}\label{thm:fraction-bias-SP-CE-vs-RCT-CE}
    In the cost-excluded design, the ratio of the shadow price estimator bias to the standard estimator bias is bounded as follows.
    \begin{enumerate}
    \item In the proportional cost setting, assuming that $(\bar{\bm{v}}^{\exp,\text{\emph{CE}}} - \bm{a}^{\text{\emph{CE}},0})\cdot\bm\beta>0$,
    \[
    \frac{\abs{\hat{\Delta}^{\text{\emph{prop}}}_{\text{\emph{SP-CE}}} - \Delta^{\text{\emph{prop}}}}}{\abs{\hat{\Delta}^{\text{\emph{prop}}}_{\text{\emph{RCT-CE}}} - \Delta^{\text{\emph{prop}}}} }\le \frac{(\bm{a}^{\text{\emph{CE}},0} - \bm{a}^{\text{\emph{CE}},1})\cdot\bm\beta}{(\bar{\bm{v}}^{\exp,\text{\emph{CE}}} - \bm{a}^{\text{\emph{CE}},0})\cdot\bm\beta}.
    \]
    \item In the fixed cost setting, assuming that $(\bar{\bm{v}}^{\exp,\text{\emph{CE}}} - (1-\rho)\kappa\bm{e} - \bm{a}^{\text{\emph{CE}},0})\cdot\bm\beta>0$ and either (i) $\sum_{i=1}^{n_d}\lambda_i + \beta_i \le \sum_{j=1}^{n_s}\gamma_j$ or (ii) $\sum_{j=1}^{n_s}\gamma_j \le \sum_{i=1}^{n_d}\lambda_i$,
    \[
    \frac{\abs{\hat{\Delta}^{\text{\emph{fixed}}}_{\text{\emph{SP-CE}}} - \Delta^{\text{\emph{fixed}}}}}{\abs{\hat{\Delta}^{\text{\emph{fixed}}}_{\text{\emph{RCT-CE}}} - \Delta^{\text{\emph{fixed}}}} }\le \frac{(\bm{a}^{\text{\emph{CE}},0} - \bm{a}^{\text{\emph{CE}},1})\cdot\bm\beta}{(\bar{\bm{v}}^{\exp,\text{\emph{CE}}} - (1-\rho)\kappa\bm{e} - \bm{a}^{\text{\emph{CE}},0})\cdot\bm\beta}.
    \]
    \end{enumerate}
    In both cases, $a_i^{\text{\emph{CE}},0}$ and $a_i^{\text{\emph{CE}},1}$ denote the demand shadow prices of each type in the global control and global treatment states, respectively.
\end{theorem}

Theorem~\ref{thm:fraction-bias-SP-CE-vs-RCT-CE} bounds the SP estimator bias as a fraction of the RCT estimator bias with the ratio of two quantities. The numerator depends on the difference between demand marginal values in global control and global treatment, while the denominator depends on the difference between the average and marginal value of demand in the experiment and control states. The denominator is large when there is a lot of interference between units, meaning marginal values are much smaller than average values. The numerator is small when the treatment itself does not significantly change the amount of interference in the marketplace. In other words, the shadow price estimator is most effective when the treatment effect is marginal (common in practice) and/or when the marketplace suffers from a lot of interference. We will verify these conclusions numerically in Section~\ref{sec:numerical-experiments}.

\subsection{Cost-Included Design}

We can analogously define a shadow price estimator in the cost-included design. In this case, there is no need for discounted shadow prices, as the dual matching formulation already has distinct shadow prices for treated and untreated demand units of each type.

\begin{definition}
    For cost-included experiments, the shadow price (SP) estimator is defined as
\begin{equation*}\label{est:SP-CI}
     \hat\Delta^{\tau,m}_{\text{SP-CI}} = \frac{1}{\tau}
    \left(
    \frac{1}{\rho}\sum_{i=1}^{n_d}A_i^{\tau, \treatment, \text{CI}, m} D_i^{\tau,\treatment}
    -
    \frac{1}{1-\rho} \sum_{i=1}^n A_i^{\tau, \control, \text{CI},m}  D_i^{\tau, \control}
    \right).
\end{equation*}
\end{definition}

Similarly, obtaining the fluid limit form of the shadow price estimator is a straightforward application of Proposition~\ref{prop:fluid-limit}.

\begin{lemma}\label{lem:fluid-limit-sp-ci}
    The cost-included shadow price estimator in the fluid limit becomes
    \begin{align*}
            \hat{\Delta}_{\text{\emph{SP-CI}}}^{m} &= \lim_{\tau\to\infty} \hat\Delta_{\text{\emph{SP-CI}}}^{\tau, m} 
                    = \bm{a}^{\treatmentemph, \text{\emph{CI}},m} \cdot (\bm\lambda+\bm\beta) - \bm{a}^{\controlemph, \text{\emph{CI}},m} \cdot \bm\lambda .
                    \end{align*}
\end{lemma}

As with the cost-excluded design, it is of interest to verify that the cost-included shadow price estimator has bias-reducing properties. Unfortunately, such a property is more difficult to establish in this setting, since, as discussed in Section~\ref{sec:rct}, the cost-included RCT estimators can sometimes be unbiased through luck, when the positive bias resulting from overestimating the marginal value of an additional unit cancels with the negative bias resulting from deprioritizing lower-value treated units in the matching. Nevertheless, we can obtain the following guarantees.

\begin{theorem}\label{thm:sp-ci-reduces-bias}
    In the cost-included setting, for any cost model $m$, the SP estimator verifies:
    \begin{enumerate}
        \item If $\Gamma<\Gamma_{\min}$ (RCT estimator always negatively biased), then
        \[
        \abs{\hat{\Delta}^m_{\text{\emph{SP-CI}}}-\Delta^m}\le \abs{\hat{\Delta}^m_{\text{\emph{RCT-CI}}}-\Delta^m},
        \]
        and furthermore
        \[
        \lim_{\alpha \to 0} \frac{\abs{\hat{\Delta}^{\text{\emph{prop}}}_{\text{\emph{SP-CI}}}-\Delta^{\text{\emph{prop}}} }}{\abs{\hat{\Delta}^{\text{\emph{prop}}}_{\text{\emph{RCT-CI}}}-\Delta^{\text{\emph{prop}}} }} 
        = 
        \lim_{\kappa \to 0} \frac{\abs{\hat{\Delta}^{\text{\emph{fixed}}}_{\text{\emph{SP-CI}}}-\Delta^{\text{\emph{fixed}}} }}{\abs{\hat{\Delta}^{\text{\emph{fixed}}}_{\text{\emph{RCT-CI}}}-\Delta^{\text{\emph{fixed}}} }} = 0,
        \]
        as long as $\alpha \le \frac{1}{2}$ in the proportional-cost model and $\kappa\le \frac{1}{2}\min_j\max_i v_{i,j}$ in the fixed-cost model.
        \item If $\Gamma=\Gamma_0$ (RCT estimator always positively biased), then
        \[
            \hat{\Delta}^m_{\text{\emph{SP-CI}}} = \hat{\Delta}^m_{\text{\emph{RCT-CI}}}.
        \]
    \end{enumerate}
\end{theorem}

Theorem~\ref{thm:sp-ci-reduces-bias} revisits the conditions from Theorem~\ref{thm:GTE:Bias-RCT-CI-can-be-high}, in which we established nonzero bias for the standard estimator. Under these same conditions, the shadow price estimator is no worse (and can be substantially better) than the standard estimator. In the first case in particular, while the relative bias of the standard estimator counter-intuitively grew to infinity as the treatment cost tended to zero, the bias of the shadow price estimator relative to that of the standard estimator tends to zero, a much more desirable behavior. The only condition is that the cost not be too high, a very reasonable condition since practical experiments rarely change prices by more than a few percentage points.

\begin{figure}[h!]
    \caption{Bias of the SP-CI estimator on the simple example from Figure~\ref{fig:setup}.}
    \begin{minipage}{0.49\columnwidth}
\begin{tikzpicture}
    \begin{axis}[
        clip=false,
        width=7cm,
        height=7cm,
        xlabel={$\eta$},
        ylabel={$\Phi$},
        grid=major,
        axis lines=left,
        legend style={cells={align=left}},
        legend cell align={left},
        legend pos=outer north east,
        axis line style={thick},
        xtick={0,0.5,1},
        xticklabels={0,$\rho=0.5$,1},
        ytick={0,2,2.4,4.36,4.7},
        yticklabels={
            0,
            {\color{ForestGreen}$\Phi_{\text{CI}}(\lambda^{\text{con}}{,}\,0)$},
            {\color{Aquamarine}$\hat{\Phi}_{\text{SP-CI}}(\lambda^{\text{con}}{,}\,0)$},
            {\color{ForestGreen}$\Phi_{\text{CI}}(0{,}\,\lambda^{\text{tre}})$},
            {\color{Aquamarine}$\hat{\Phi}_{\text{SP-CI}}(0{,}\,\lambda^{\text{tre}})$}
        },
        every tick label/.append style={font=\footnotesize},
        every axis x label/.style={at={(ticklabel* cs:1)},anchor=north},
        every axis y label/.style={at={(ticklabel* cs:1)},anchor=east},
        enlargelimits,
    ]
        \addplot[
            color=White,
            domain=0:1,
            samples=100
        ]
        coordinates {
            (0.9, 0.1)
            (1, 0)
        };
        \addplot[
            color=Aquamarine,
            very thick,
            domain=0:1,
            samples=100,
            solid
        ]
        coordinates {
            (0, 2.4)
            (1, 4.7)
        };
        \addplot[
            color=ForestGreen,
            very thick,
            domain=0:1,
            samples=100,
            dashdotted
        ]
        coordinates {
            (0, 2)
            (0.17, 2.79)
            (0.83, 4.31)
            (1, 4.36)
        };
        \draw [<->,thick] (axis cs:-0.02,2) -- (axis cs:-0.02,4.36) node[midway, right,xshift=-0.15cm] {\footnotesize$\Delta$};
        \draw [<->,thick] (axis cs:1.02,2.4) -- (axis cs:1.02,4.7) node[midway, left,yshift=-0.25cm] {\footnotesize$\hat{\Delta}_{\text{SP-CI}}$};
        \node[anchor=center] at (0.5,-1.8) {(a) $\lambda^{\text{con}} = 1,~\lambda^{\text{tre}}=4$};
    \end{axis}
\end{tikzpicture}%
\end{minipage}\hfill
\begin{minipage}{0.49\columnwidth}
 \begin{tikzpicture}
        \begin{axis}[
            clip=false,
            width=7cm,
            height=7cm,
            xlabel={$\eta$},
            ylabel={$\Phi$},
            grid=major,
            axis lines=left,
            legend style={cells={align=left}},
            legend cell align={left},
            legend pos=outer north east,
            axis line style={thick},
            xtick={0,0.5,1},
            xticklabels={0,$\rho=0.5$,1},
            ytick={4.5,4.5675,4.79375,4.81875},
            axis y discontinuity=crunch,
            ymin=3.5,
            ymax=5,
            yticklabels={
                \raisebox{-5ex}{\color{ForestGreen}$\Phi_{\text{CI}}(\lambda^{\text{con}}{,}\,0)$},
                \raisebox{1ex}{\color{ForestGreen}${\Phi}_{\text{CI}}(0{,}\,\lambda^{\text{tre}})$},
                \raisebox{-4ex}{\color{Aquamarine}$\hat{\Phi}_{\text{SP-CI}}(0{,}\,\lambda^{\text{tre}})$},
                \raisebox{4ex}{\color{Aquamarine}$\hat{\Phi}_{\text{SP-CI}}(\lambda^{\text{con}}{,}\,0)$},
            },
            every tick label/.append style={font=\footnotesize},
            every axis x label/.style={at={(ticklabel* cs:1)},anchor=north},
            every axis y label/.style={at={(ticklabel* cs:1)},anchor=east},
            enlargelimits,
        ]
            \addplot[
                color=Aquamarine,
                very thick,
                domain=0:1,
                samples=100,
                solid
            ]
            coordinates {
                (0, 4.81875)
                (0.25, 4.8125)
                (0.5, 4.80625)
                (1, 4.79375)
            };\label{plot7}
            \addplot[
                color=ForestGreen,
                very thick,
                domain=0:1,
                samples=100,
                dashdotted
            ]
            coordinates {
                (0, 4.5)
                (0.25, 4.8125)
                (0.5, 4.80625)
                (1, 4.5675)
            };\label{plot6}
            \draw [|-|,semithick] (axis cs:-0.02,4.5) -- (axis cs:-0.02,4.5675) node[midway, left,xshift=0.05cm,yshift=0.05cm] {\footnotesize$\Delta$};
            \draw [|-|,semithick] (axis cs:1.02,4.79375) -- (axis cs:1.02,4.81875) node[midway, left,yshift=0.25cm,xshift=0.25cm] {\footnotesize$\hat{\Delta}_{\text{SP-CI}}$};
            \node[anchor=center] at (0.5,2.92) {(b) $\lambda^{\text{con}} = 3,~\lambda^{\text{tre}}=5$};
        \end{axis}
    \end{tikzpicture}%
    \end{minipage}
    \vspace{2pt}
\begin{center}
\begin{tikzpicture}
\matrix(first)[
  matrix of nodes,
  execute at empty cell={\node[draw=none]{};},
  anchor=north,
  inner sep=0.2em,
  name=table,
  font=\scriptsize,
  nodes={anchor=west,align=left}
]at(0,-1.1)
{  
    \ref{plot6} & $\Phi_{\text{CI}}((1-\eta)\lambda^{\text{con}},\eta\lambda^{\text{tre}})$ &
    \ref{plot7} & $\Phi_{\text{CI}}((1-\rho)\lambda^{\text{con}},\rho\lambda^{\text{tre}}) + (\eta - \rho)\frac{d}{d\eta}\Phi_{\text{CI}}((1-\eta)\lambda^{\text{con}},\eta\lambda^{\text{tre}})$ 
    \\
    & (total value) & & (first-order approximation of value function)\\
};
\end{tikzpicture}
\end{center}
    \label{fig:sp-ci}
    {\footnotesize\textit{Note.} The shadow price estimator builds a single linear approximation of the value function by using the shadow prices at the experiment point.}
\end{figure}
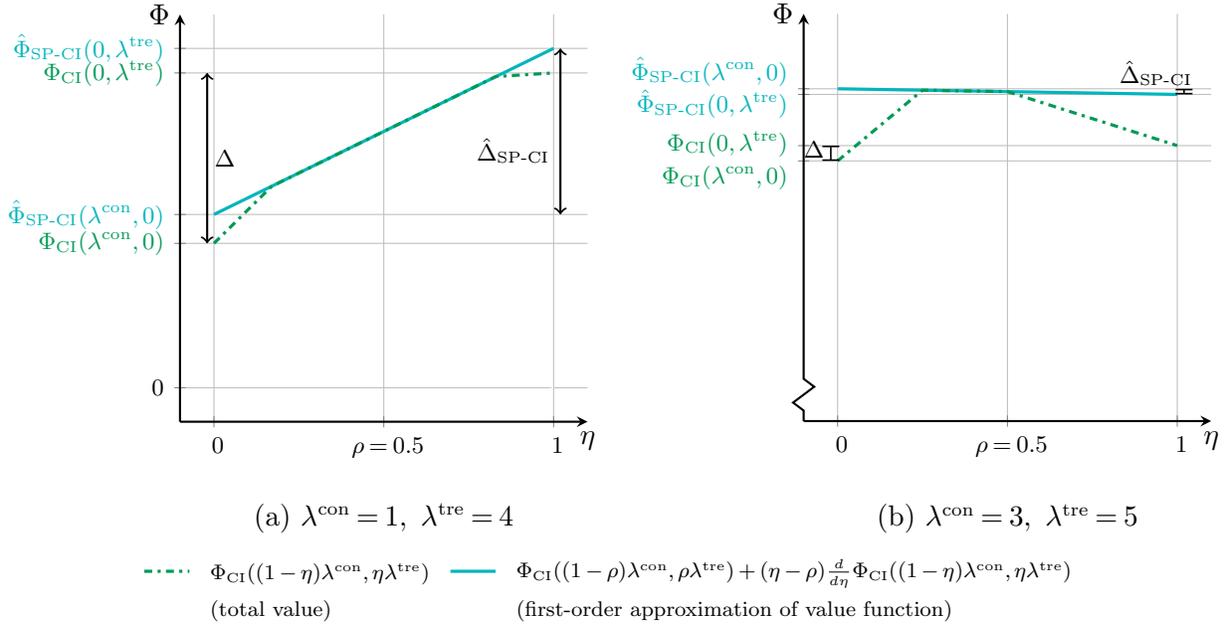

The main takeaway from this section is that, regardless of which design the platform chooses, using the shadow price estimator is likely to reduce the bias caused by interference. In the next section, we compare the performance of the two designs in terms of bias reduction.
\section{Design Recommendations for Platforms}

\label{sec:design-recommendations}

Having established the bias-reducing properties of the shadow price estimator in both designs, we now consider the problem of which design may be preferable.

\subsection{Include cost or exclude it?}

Platforms may have flexibility in choosing whether to implement a cost-included or cost-excluded design when evaluating pricing interventions. While the shadow price estimator provably reduces bias in both situations (under some conditions), it is of interest to understand whether one design is inherently less biased than another. One way to do so is to identify necessary and sufficient conditions for each design to lead to a zero-bias estimator, regardless of the treatment fraction $\rho$.

\begin{theorem}
    \label{thm:ce-vs-ci}
    The shadow price estimator is unbiased for every $0\le\rho\le1$ if and only if one of the following conditions holds:
    \begin{enumerate}
        \item The platform uses the cost-excluded design and 
            \begin{enumerate}
                \item $\bm{a}^{\text{\emph{CE}},0}=\bm{a}^{\text{\emph{CE}},1} $  under the proportional-cost model, 
                \item $\bm{a}^{\text{\emph{CE}},0}=\bm{a}^{\text{\emph{CE}},1}$, and additionally, either (i) 
                $\sum_{i=1}^{n_d}\lambda_i + \beta_i \le \sum_{j=1}^{n_s}\gamma_j $
                or  (ii)  
                $\sum_{j=1}^{n_s}\gamma_j \le \sum_{i=1}^{n_d}\lambda_i $
                under the fixed-cost model.
            \end{enumerate}
        \item The platform uses the cost-included design and $\bm{a}^{\treatmentemph,\text{\emph{CI}},m,0}=\bm{a}^{\treatmentemph,\text{\emph{CI}},m,\rho}$ and $\bm{a}^{\controlemph,\text{\emph{CI}},m,0}=\bm{a}^{\controlemph,\text{\emph{CI}},m,\rho}$ for every $0\le\rho\le 1$.
    \end{enumerate}
\end{theorem}

The first part of Theorem~\ref{thm:ce-vs-ci} is not very surprising. It follows almost directly from Theorem~\ref{thm:fraction-bias-SP-CE-vs-RCT-CE} that the first two conditions are sufficient to ensure the shadow price estimator is unbiased in the cost-excluded design. It turns out that if we want the estimator to be unbiased \emph{for any choice} of $\rho$, these conditions are also necessary. Intuitively, the conditions simply mean that the interference bias should be the same under global control and global treatment --- in other words, the treatment itself is not leading to more or less interference. This is reasonable in a world of marginal treatments.

The second part of Theorem~\ref{thm:ce-vs-ci} looks superficially similar to the first part. However, it tells a very different story. In particular, it requires that the marginal value of an untreated demand unit is exactly the same in global control (when there are lots of untreated demand units) and in global treatment (when there are no untreated demand units at all). Under mild assumptions, we can show that this is a much stronger requirement, as described in the following corollary.

\begin{corollary}
\label{cor:ci-unbiased}
Assume that the matching problem in the global control and global treatment states has a unique and nondegenerate optimum. Then the SP-CI estimator is unbiased for any $\rho$ if and only if the optimal matching solution in global treatment verifies
\[
x_{i,j}^{\treatmentemph,\text{\emph{CI}},m,1} = \begin{cases} \lambda_i + \beta_i, & \text{if } v_{i,j} = \max_{1\le k\le n_s} v_{i,k},\\
0,&\text{otherwise.}\end{cases}
\]
\end{corollary}

Corollary~\ref{cor:ci-unbiased} establishes that most of the time, the SP-CI estimator is only unbiased for any $\rho$ \emph{if there is no interference bias to begin with}. In other words, even if the treatment itself has no effect on the structure of interference in the marketplace, the shadow price estimator in the cost-included design may still be biased for some $\rho$. Corollary~\ref{cor:ci-unbiased} is a much more stringent condition --- indeed, if it holds, the standard estimator would also be unbiased for every $\rho$. We can therefore provide an initial answer to platforms wondering which design is preferable for pricing interventions: the cost-excluded design is much more likely to be unbiased for treatments with small effects (as most practical treatments are), whereas the cost-included design is only unbiased when the standard estimator is also unbiased. We will also consider this decision numerically in Section~\ref{sec:numerical-experiments}.

\subsection{Simulation-based Estimation}

An alternative approach to the standard estimator is to simply estimate the arrival rates of each type under control and treatment, and evaluate the counterfactual matching values under these estimated arrival rates. This approach is called the ``Two-LP estimator'' by \citet{bright2024reducing}, since it relies on solving two matching problems. We refer to it as ``simulation-based'' (SB) since it involves reconstructing the matching problem with two different demand realizations and comparing the total obtained values.

\begin{definition}\label{def:simulation-based}
The simulation-based estimator is given by
\begin{equation*}
    \hat\Delta_{\text{SB}}^{\tau, m}  =  
    \Phi^m_{\text{CI}}\left(\bm{0},  \frac{1}{\tau \rho} \bm{D}^{\tau, \treatment}, \frac{1}{\tau}\bm{S}^{\tau, \bm{\gamma}} \right) - 
    \Phi^m_{\text{CI}}\left( \frac{1}{\tau (1-\rho) } \bm{D}^{\tau, \control}, \bm{0}, \frac{1}{\tau}\bm{S}^{\tau, \bm{\gamma}} \right).
\end{equation*}
\end{definition}

Similarly to Remark~\ref{rk:gte-equivalent-definition}, there is an equivalent definition of the SB estimator using the cost-excluded matching formulation. It is easy to show that the simulation-based estimator is unbiased in the fluid limit.

\begin{lemma}
    In the fluid limit,
    \begin{align}
            \hat{\Delta}^m_{\text{\emph{SB}}} &= \lim_{\tau\to\infty} \hat \Delta_{\text{\emph{SB}}}^{\tau, m}  
            = \Phi_{\text{\emph{CI}}}^m (\bm{0}, \bm\lambda+\bm\beta, \bm\gamma) - \Phi_{\text{\emph{CI}}}^m (\bm{\lambda}, \bm{0}, \bm\gamma) = \Delta^m.
        \end{align}
\end{lemma}

The simulation-based approach is attractive, and it allows the platform to sidestep the cost-included vs. cost-excluded design question. However, it also has drawbacks. On the practical side, it requires running the matching algorithm for counterfactual configurations --- but such a counterfactual matching means that much more data must be collected from the experimental state in order to replicate the matching function. For example, \cite{azagirre2024better} describe how the matching values in the dispatch system at Lyft are updated dynamically throughout the day using reinforcement learning --- which makes reproducing a matching configuration offline a complex engineering task.

On the theoretical side, \cite{bright2024reducing} showed that it may lead to significant bias in the low-data setting (i.e., when the fluid limit is a poor approximation of marketplace dynamics), especially when the experiment is asymmetric, i.e., $\rho$ is close to 0 or 1. We will verify numerically that this behavior continues to occur in pricing interventions in Section~\ref{sec:numerical-experiments}.

\section{Numerical Experiments} \label{sec:numerical-experiments}

The majority of our analysis is conducted in the fluid limit. We now validate our key findings numerically in a finite-sample setting. We focus on the proportional-cost model, and present similar results for the fixed-cost model in Appendix~\ref{Appendix-sec:fixed-cost-numerical-experiments}.

\subsection{Setup}

We consider a synthetic data setting, motivated by geographical matching (e.g., in ride-hailing). We generate 50 matching instances with $n_s=10$ demand types and $n_d=10$ supply types, each time by sampling $n_d+n_s$ points uniformly in the unit square. We associate each type with one point (or \emph{location}), and define the value of matching a demand unit of type $i$ (located at $\bm{z}_d^i$) to a supply unit of type $j$ (located at $\bm{}z_s^j$) as $v_{i,j}= e^{-\norm{\bm{z}_d^i-\bm{z}_s^j}_2}$. We let $\tau=1$, $\lambda_i=13$ and $\beta_i=3$ for each demand type $i$. We also let $\gamma_j=\Gamma \lambda_i$, where $\Gamma$ takes 30 possible values ranging from 0.3 (undersupply) to 3 (oversupply). Due to similarity in results, we only consider the proportional-cost model in the main text and present analogous results for the fixed-cost setting in the appendix. We consider three possible values of the proportional discount parameter $\alpha$ (5\%, 10\%, 20\%). For each matching setting, we compute the global treatment effect and estimator values from 50 Poisson samples.

\subsection{Bias Reduction}

\begin{figure}[t]
\caption{Comparing Estimator Biases.} 
\begin{center}
\includegraphics[height=4.7in]{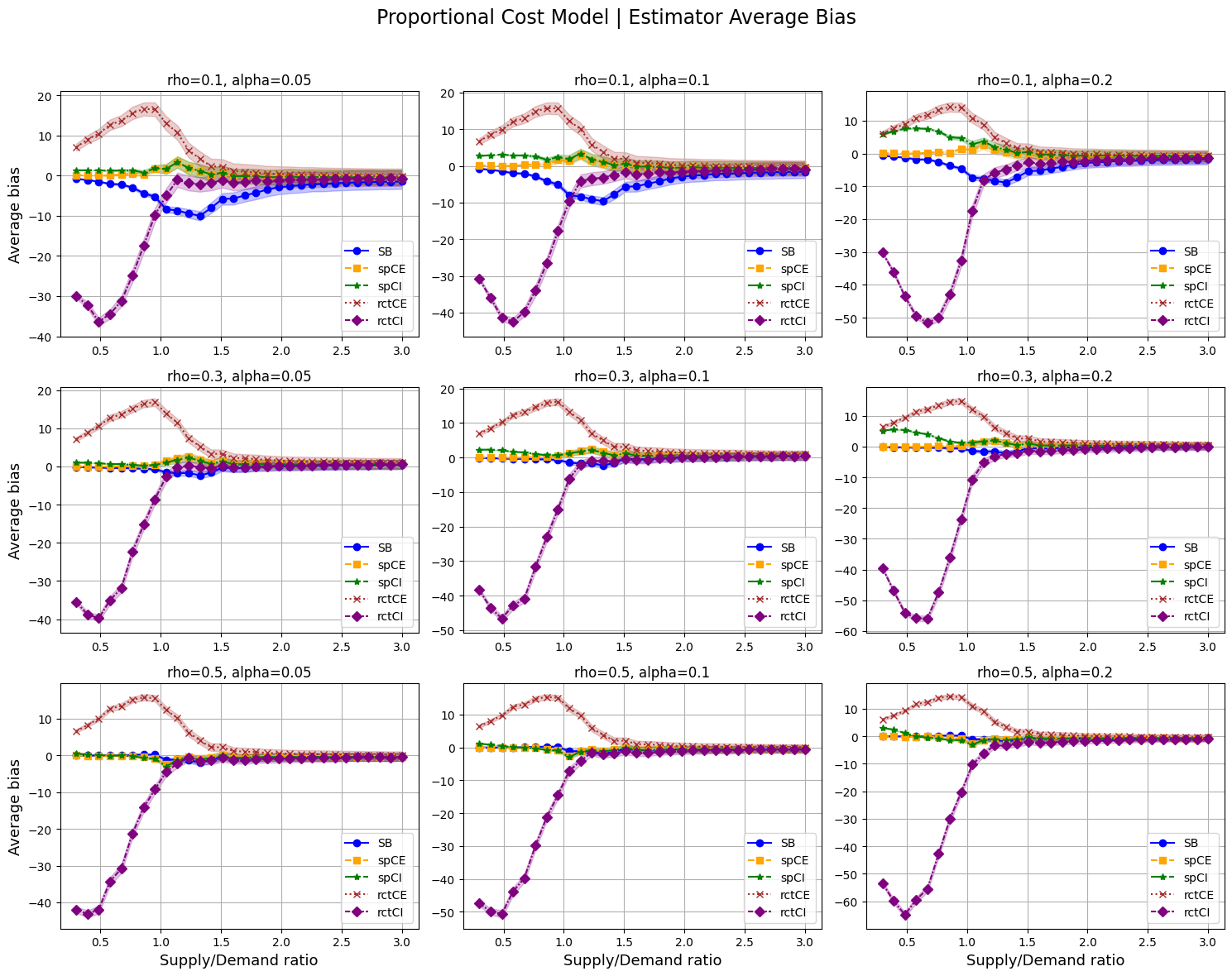}
\label{fig:simulations-singed-biases}
\end{center}
{\footnotesize\textit{Note.} The shadow price estimator in the cost-excluded setting reduces bias across all scenarios. The RCT-CE estimator is always positively biased, while the RCT-CI estimator shows significant negative bias in the undersupply regime. The simulation-based estimator (SB) exhibits significant bias in asymmetric experiments (small $\rho$).}
\end{figure}

Figure~\ref{fig:simulations-singed-biases} displays the (signed) average biases of all estimators along with their standard errors for simulated experiments with treatment fractions ($\rho$) of $0.1, 0.3$, and $0.5$, and for three discount levels ($\alpha$) of $5\%$, $10\%$, and $20\%$. As predicted by theory in the fluid limit, the RCT-CE estimator overestimates the global treatment effect, inducing positive bias. Meanwhile, the RCT-CI estimator exhibits negative bias in the undersupply regime and transitions to positive bias at certain supply-to-demand ratios. We observe that both shadow price estimators consistently reduce bias compared to their standard counterparts. In addition, the simulation-based estimator (SB) shows significant bias for small $\rho$ values (asymmetric experiments).

\begin{figure}[t]
\caption{Bias Reduction from the Shadow Price Estimator Across Designs.} 
\begin{center}
\includegraphics[height=4.7in]{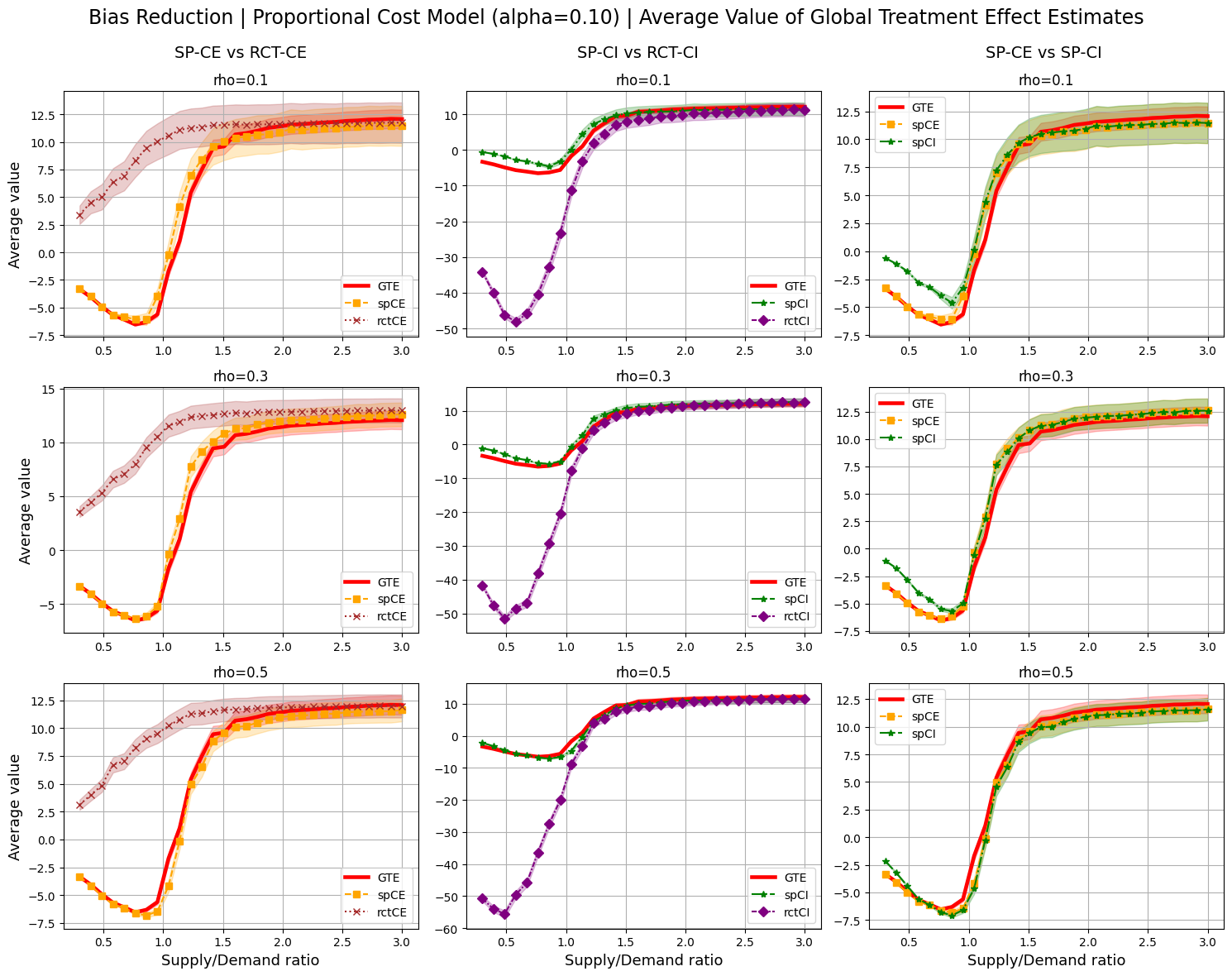}
\label{fig:simulations-bias-reduction}
\end{center}
{\footnotesize\textit{Note.} Pairwise comparisons of shadow price and standard estimators across designs. The shadow price estimator in the cost-excluded design outperforms all others in bias reduction across all experiment configurations, particularly in the low-supply regime.}
\end{figure}

Figure~\ref{fig:simulations-bias-reduction} provides pairwise comparisons of shadow price and standard estimators across designs, as well as the shadow price estimators across designs, all in our finite-sample setting. The left column compares the SP and RCT estimators in the cost-excluded setting, with the former much more closely tracking the global treatment effect. Similarly, in the middle column we compare the SP and RCT estimators in the cost-included setting, with a similar observation. Finally, the third column illustrates the message of Section~\ref{sec:design-recommendations}, namely that in a head-to-head comparison between shadow price estimators in the cost-included and cost-excluded settings, the latter tends to have lower bias. Overall, the SP estimator in the cost-excluded design outperforms others in bias reduction across most experimental configurations.

\subsection{Bias Reduction Guarantees in Practice}

\begin{figure}[t]
\caption{Effect of Treatment Fraction on Bias Reduction Guarantees.} 
\begin{center}
\includegraphics[height=2in]{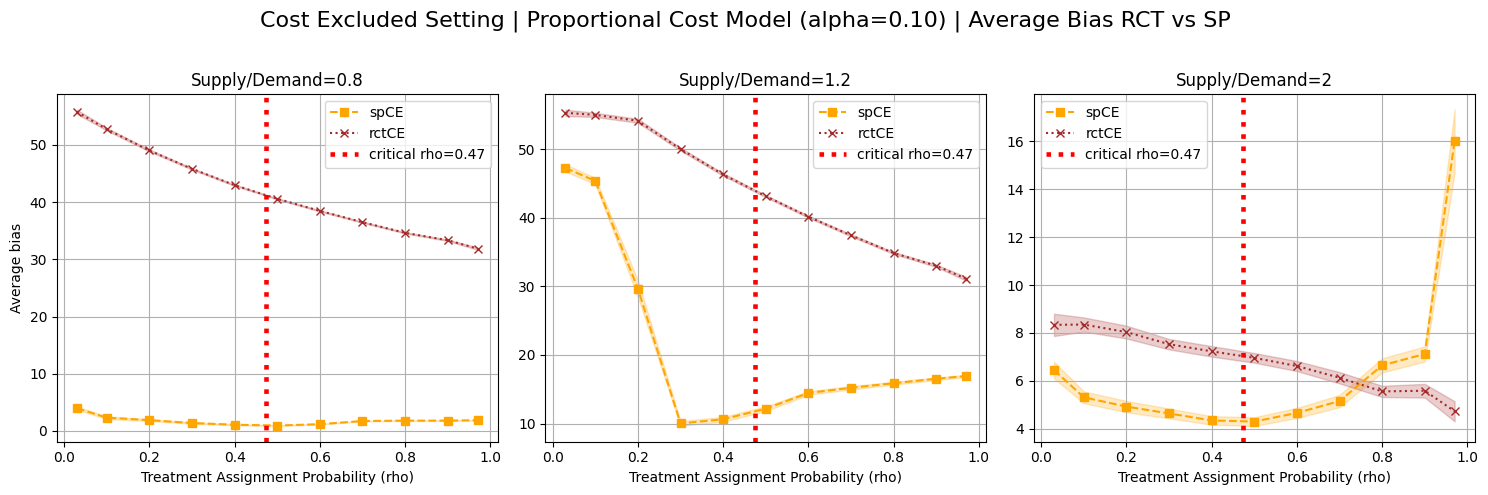}
\label{fig:simulations-thm3}
\end{center}
\end{figure}

In Theorems~\ref{thm:RCT-CE-vs-SP-CE} and \ref{thm:fraction-bias-SP-CE-vs-RCT-CE}, we established conditions for bias reduction from the shadow price estimator and bounds on the amount of bias reduction, in the cost-excluded setting. Theorem~\ref{thm:RCT-CE-vs-SP-CE} provides a sufficient condition on the treatment assignment probability $\rho$ for the shadow price estimator to reduce bias compared to the standard estimator. Figure~\ref{fig:simulations-thm3} shows that the theoretical guarantee can be somewhat conservative, particularly in undersupply or balanced supply/demand configurations, where the shadow price estimator reduces bias far beyond the maximum allowed $\rho$. In oversupply, the bound becomes tighter, though there is still room to increase $\rho$ without compromising bias reduction. This suggests that additional information on the structure of the value function (to which our theory is entirely agnostic) may yield tighter guarantees on bias reduction. The choice of $\rho$ is an important consideration since the SP estimator bias can indeed increase bias in some highly asymmetric cases.

\begin{figure}[t]
\caption{Effect of Supply-to-demand Ratio on Bias Reduction Guarantees.} 
\begin{center}
\includegraphics[height=2.9in]{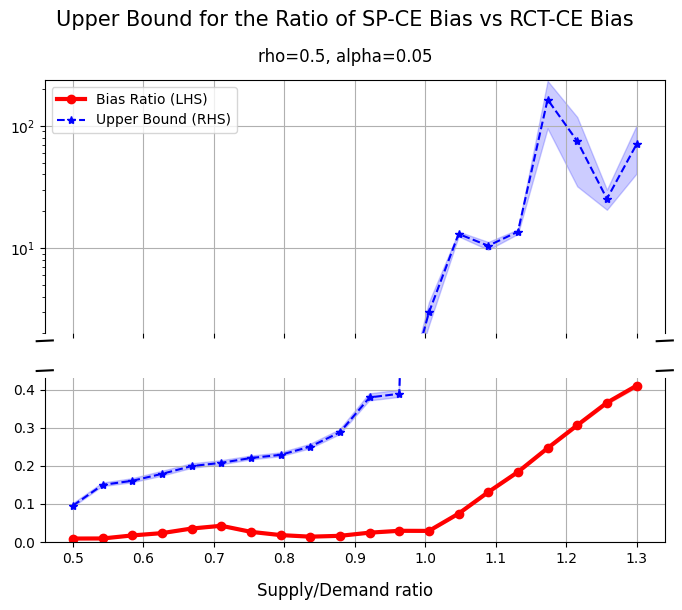}
\label{fig:simulations-thm4}
\end{center}
{\footnotesize\textit{Note.} Tightness of the upper bound of the bias ratio between SP and RCT in the cost-excluding setting. The bound is most meaningful in undersupply regimes providing a reasonable order of magnitude, while it becomes large in high-supply settings.}
\end{figure}

Meanwhile, Theorem~\ref{thm:fraction-bias-SP-CE-vs-RCT-CE} established an upper bound on the ratio between the biases of the shadow price and standard estimators in the cost-excluded setting. Figure~\ref{fig:simulations-thm4} explores the tightness of this upper bound. As expected, the bound is conservative, though reasonably meaningful in the undersupply setting: for example, when the supply/demand ratio is $0.7$, the theory guarantees that the SP bias removes at least 80\% of the RCT bias, while in practice it removes around 95\%. In high-supply regimes, the bound's performance deteriorates sharply, but the actual estimator performance decays more gracefully as both estimators converge to the unbiased ground truth.

\section{Conclusion} \label{sec:conclusion}

This paper addresses the issue of interference bias in pricing experiments, in marketplaces with a matching component. We first identify a crucial design decision --- whether to discriminate, at match time, between units paying full price and units paying a discounted price, which only coexist in the marketplace as an artifact of the experimentation process. We first show that a simple duality-based correction on the standard estimator (called the shadow price estimator) reduces bias significantly no matter what the platform decides. Our analysis further recommends treating the treated and untreated units as indistinguishable to reduce bias, even though this is not revenue-maximizing in the short term. 

\bibliographystyle{informs2014} 
\SingleSpacedXI

\DoubleSpacedXI

\newpage

\begin{APPENDICES}
\section{Matching Properties}

We first introduce the following simple lemma.

\begin{lemma} \label{lem:positive-weights}
    Consider the optimal cost-excluded or cost-included matching, with positive edge weights. At least one of the following statements is true:
    \begin{enumerate}
        \item Every demand unit is matched to a supply unit.
        \item Every supply unit is matched to a demand unit.
    \end{enumerate}
\end{lemma}

\proof{Proof of Lemma~\ref{lem:positive-weights}.}
Assume neither statement is true. Then there is at least one unmatched demand unit and one unmatched supply unit. Matching them together yields positive payoff which increases the objective value, contradicting optimality of the matching.
\\
\Halmos
\endproof

\proof{Proof of Lemma~\ref{lem:CE-unique}.}
Under the proportional model, the proof is trivial, since multiplying all weights by $(1-\alpha)$ is akin to a change in units.

Under the fixed-cost model, let $x^*_{i,j}$ designate the optimal solution of $\Phi_{\text{CE}}(\bm{d}, \bm{s})$, and let $\hat{x}_{i,j}$ designate the optimal solution of $\bar{\Phi}_{\text{CE}}(\bm{d}, \bm{s})$ (the latter designating the cost-excluded matching problem with all weights replaced by $v_{i,j}-\kappa$). Because each optimal solution is feasible for the other problem, we can write
\begin{subequations}
\begin{align}
    \sum_{i=1}^{n_d}\sum_{j=1}^{n_s}(v_{i,j}-\kappa)x_{i,j}^* &\le \sum_{i=1}^{n_d}\sum_{j=1}^{n_s}(v_{i,j}-\kappa)\hat{x}_{i,j},\label{eq:lemma1-eq1}\\
    \sum_{i=1}^{n_d}\sum_{j=1}^{n_s}v_{i,j}\hat{x}_{i,j} &\le \sum_{i=1}^{n_d}\sum_{j=1}^{n_s}v_{i,j}{x}_{i,j}^*.\label{eq:lemma1-eq2}
\end{align}
\end{subequations}

By Lemma~\ref{lem:positive-weights}, we know that
\[
\sum_{i=1}^{n_d}\sum_{j=1}^{n_d}\hat{x}_{i,j}=\sum_{i=1}^{n_d}\sum_{j=1}^{n_d}{x}_{i,j}^*=\min\left(\sum_{i=1}^{n_d}d_i,\sum_{j=1}^{n_s}s_j\right),
\]
which means we can re-write~\eqref{eq:lemma1-eq1} as
\[
\sum_{i=1}^{n_d}\sum_{j=1}^{n_s}v_{i,j}x_{i,j}^* \le \sum_{i=1}^{n_d}\sum_{j=1}^{n_s}v_{i,j}\hat{x}_{i,j},
\]
completing the proof together with \eqref{eq:lemma1-eq2}.
\\
\Halmos
\endproof

\proof{Proof of Lemma~\ref{lem:fixed-cost-shadow-price}.}
Under the proportional-cost model, multiplying the dual variables by $(1-\alpha)$ preserves feasibility and correctly scales the objective value by $(1-\alpha)$.

Under the fixed-cost model, removing $\kappa$ from either all the demand shadow prices or all the supply shadow prices also preserves that $\tilde{A}_i^{\tau,\text{fixed}} + \tilde{B}_j^{\tau,\text{fixed}}\ge v_{i,j}-\kappa$. Additionally, the objective value is reduced by $\kappa\min\left(\sum_{i=1}^{n_d}D_i^{\tau,\ex} , \sum_{j=1}^{n_s}S_j^{\tau,\bm\gamma}\right)$, which we know from Lemma~\ref{lem:positive-weights} is the correct amount. All that remains to be shown is that the modified shadow prices remain nonnegative.

Consider the first case without loss of generality, i.e., where $\sum_{i=1}^{n_d}D_i^{\tau,\ex} < \sum_{j=1}^{n_s}S_j^{\tau,\bm\gamma}$. Pick any demand type $k$, and let $0 < \varepsilon < \sum_{j=1}^{n_s}S_j^{\tau,\bm\gamma} - \sum_{i=1}^{n_d}D_i^{\tau,\ex}$. From linear programming sensitivity analysis we can write:
\[
\sum_{i=1}^{n_d}A_i^{\tau,\text{CE}}(D_i^{\tau,\ex} + \varepsilon\delta_{ik}) + \sum_{j=1}^{n_s}B_j^{\tau,\text{CE}}S_j^{\tau,\bm\gamma} \ge \Phi_{\text{CE}}(\bm{D}^{\tau,\ex}+\varepsilon \bm{e}_k,\bm{S}^{\tau,\bm\gamma}),
\]
which we can use strong duality to re-write as
\[
A_k^{\tau,\text{CE}}\varepsilon \ge \Phi_{\text{CE}}(\bm{D}^{\tau,\ex}+\varepsilon \bm{e}_k,\bm{S}^{\tau,\bm\gamma}) - \Phi_{\text{CE}}(\bm{D}^{\tau,\ex},\bm{S}^{\tau,\bm\gamma}) \ge \varepsilon\min_{i,j}v_{i,j} \ge \varepsilon\kappa,
\]
where the penultimate inequality follows from the construction of $\varepsilon$.
\\
\Halmos
\endproof

\proof{Proof of Proposition \ref{prop:fluid-limit}.} 

The first part of Proposition~\ref{prop:fluid-limit} exactly restates Theorem 1 of \cite{bright2024reducing}, so we refer readers to this proof and only consider the second part.

We first prove the following claim:
\begin{equation}
\label{eq:scaling}
\frac{1}{\tau}\Phi_{\text{\emph{CI}}}^m\left(\bm{D}^{\tau,\bm\lambda^\control},  \bm{D}^{\tau,\bm\lambda^\treatment},\bm{S}^{\tau,\bm\gamma}\right) = \Phi_{\text{\emph{CI}}}^m\left(\frac{1}{\tau}\bm{D}^{\tau,\bm\lambda^\control},  \frac{1}{\tau}\bm{D}^{\tau,\bm\lambda^\treatment},\frac{1}{\tau}\bm{S}^{\tau,\bm\gamma}\right).
\end{equation}
This claim is easily proven from the dual:
\begin{align*}
    \frac{1}{\tau}\min \quad & \sum_{i=1}^{n_d} A_i^{\tau,\control, m} D_i^{\tau,\lambda^\control} + A_i^{\tau,\treatment, m} D_i^{\tau,\lambda^\treatment} + \sum_{j=1}^{n_s} B_j^{\tau,m} S_j^{\tau,\bm\gamma} \\
    = \min \quad &\sum_{i=1}^{n_d} A_i^{\tau,\control, m} \frac{1}{\tau}D_i^{\tau,\lambda^\control} + A_i^{\tau,\treatment, m} \frac{1}{\tau}D_i^{\tau,\lambda^\treatment} + \sum_{j=1}^{n_s} B_j^{\tau,m} \frac{1}{\tau}S_j^{\tau,\bm\gamma}\\
        \text{s.t.}\quad & A_i^{\tau,\control, m} + B_j^{\tau,m} \ge v_{i,j}^\control & \forall i\in[n_d], j\in[n_s]\\
        & A_i^{\tau,\treatment, m } + B_j^{\tau,m} \ge v_{i,j}^\treatment & \forall i\in[n_d], j\in[n_s]\\
        & A_i^{\tau,\control, m }, A_i^{\tau,\treatment, m} \ge 0 & \forall i\in[n_d]\\
        & B_j^{\tau,m } \ge 0 & \forall j\in[n_s].
\end{align*}

By the strong law of large numbers, $\left(\frac{1}{\tau}\bm{D}^{\tau,\bm\lambda^\control},\frac{1}{\tau}\bm{D}^{\tau,\bm\lambda^\treatment},\frac{1}{\tau}\bm{S}^{\tau,\bm\gamma}\right)$ converges to $(\bm\lambda^\control,\bm\lambda^\treatment,\bm\gamma)$ almost surely, and furthermore, considering all integer values of $\tau$, we have
\[
    \Phi_{\text{\emph{CI}}}^m\left(\frac{1}{\tau}\bm{D}^{\tau,\bm\lambda^\control},  \frac{1}{\tau}\bm{D}^{\tau,\bm\lambda^\treatment},\frac{1}{\tau}\bm{S}^{\tau,\bm\gamma}\right) \le \max_{i,j}v_{i,j} \max_{\tau=1}^\infty \frac{1}{\tau}\sum_{j=1}^{n_s}{S_j^{\tau,\bm\gamma}}.
\]
Lemma 2 of \cite{bright2024reducing} shows that the right-hand side has bounded mean, so we can apply claim~\eqref{eq:scaling} and the dominated convergence theorem to conclude that
\begin{equation}
\label{eq:convergence-fluid-limit-expectation}
\lim_{\tau\to\infty} \E\left[\frac{1}{\tau}\Phi_{\text{\emph{CI}}}^m(\bm{D}^{\tau,\bm\lambda^\control},  \bm{D}^{\tau,\bm\lambda^\treatment},\bm{S}^{\tau,\bm\gamma})\right] 
        = \Phi_{\text{\emph{CI}}}^m( \bm\lambda^\control, \bm\lambda^\treatment, \bm\gamma).
\end{equation}
To prove convergence of the optimal primal and dual variables, we first assume uniqueness of the optimal primal and dual solutions of $\Phi_{\text{\emph{CI}}}^m( \bm\lambda^\controlemph, \bm\lambda^\treatmentemph, \bm\gamma)$ without loss of generality. It is easier to start with convergence of the optimal dual variables: assume we can find an infinite subsequence $\tau_1,\ldots$ such that $(\bm{A}^{\tau_i,\control,\text{CI},m},\bm{A}^{\tau_i,\treatment,\text{CI},m},\bm{B}^{\tau_i,\ex,\text{CI},m})\neq (\bm{a}^{\control,\text{CI},m},\bm{a}^{\control,\text{CI},m},\bm{b}^{\ex,\text{CI},m})$. The existence of this subsequence violates either the optimality of $(\bm{a}^{\control,\text{CI},m},\bm{a}^{\control,\text{CI},m},\bm{b}^{\ex,\text{CI},m})$ in $\Phi_{\text{\emph{CI}}}^m( \bm\lambda^\control, \bm\lambda^\treatment, \bm\gamma)$ or the convergence of the objective value~\eqref{eq:convergence-fluid-limit-expectation}. This proves almost sure convergence of the dual optimal solution. Almost sure convergence of the primal solution follows from the complementary slackness conditions: the optimal primal solution is a continuous function (solving a linear system) of the optimal dual solution, therefore almost sure convergence of one implies almost sure convergence of the other.
\\
\Halmos
\endproof

\section{RCT Estimator Analysis}

\proof{Proof of Lemma~\ref{lem:fluid-limit-rct-ce}}
We focus on the proportional-cost case, since both cases follow the same reasoning.
\begin{align*}
    \hat \Delta_{\text{RCT-CE}}^{\tau, \text{{prop}}}  &= \frac{1}{\tau} \left(  \frac{1}{\rho}\sum_{i=1}^{n_d} \sum_{j=1}^{n_s} (1-\alpha) v_{i,j}   X_{i,j}^{\tau, \treatment, \text{CE}}  - \frac{1}{1-\rho} \sum_{i=1}^{n_d} \sum_{j=1}^{n_s} v_{i,j}  X_{i,j}^{\tau, \control, \text{CE}} 
    \right)\\
    &=\frac{1}{\tau} \left(  \frac{1}{\rho}\sum_{i=1}^{n_d} \sum_{j=1}^{n_s} (1-\alpha) v_{i,j}   \frac{D_i^{\tau,\treatment}}{D_i^{\tau,\ex}}X_{i,j}^{\tau,\ex,\text{CE}}   - \frac{1}{1-\rho} \sum_{i=1}^{n_d} \sum_{j=1}^{n_s} v_{i,j}  \frac{D_i^{\tau,\control}}{D_i^{\tau,\ex}}X_{i,j}^{\tau,\ex,\text{CE}} 
    \right)\\
    &=   \frac{1}{\rho}\sum_{i=1}^{n_d} (1-\alpha) \bar{V}_i^{\tau,\ex}  \frac{D_i^{\tau,\treatment}}{\tau}  - \frac{1}{1-\rho} \sum_{i=1}^{n_d} \bar{V}_i^{\tau,\ex}  \frac{D_i^{\tau,\control}}{\tau},
\end{align*}
where $\bar{V}_i^{\tau,\ex}=(\sum_{j=1}^{n_s}v_{i,j}X_{i,j}^{\tau,\ex,\text{CE}})/D_i^{\tau,\ex}$ is the average value obtained from demand type $i$ in the experiment. From Proposition~\ref{prop:fluid-limit} and the strong law of large numbers, we obtain the limit as
\begin{align*}
    \hat \Delta_{\text{RCT-CE}}^{\text{{prop}}} &= (1-\alpha) \frac{1}{\rho}\bar{\bm{v}}^{\ex, \text{{CE}}} \cdot \rho(\bm\lambda+\bm\beta) - \frac{1}{1-\rho}  \bar{\bm{v}}^{\ex, \text{{CE}}} \cdot (1-\rho)\bm\lambda\\
    &= (1-\alpha) \bar{\bm{v}}^{\ex, \text{{CE}}} \cdot (\bm\lambda+\bm\beta) - \bar{\bm{v}}^{\ex, \text{{CE}}} \cdot \bm\lambda.
\end{align*}
\Halmos
\endproof

We are now ready to prove our first main result. We first introduce some additional notation.

Let $\bar{\Phi}^{m}_{\text{CE}}(\cdot,\cdot)$ designate the cost-excluded matching problem, where the matching weights $v_{i,j}$ are replaced with  $v_{i,j}^m \;$, where $v_{i,j}^{\text{prop}}=(1-\alpha)v_{i,j}$ and $v_{i,j}^{\text{fixed}} = v_{i,j}-\kappa \;$, i.e., 
\begin{subequations} 
\begin{align}
    \label{lp:CE-with-cost}
    \bar{\Phi}^{m}_{\text{CE}}(\bm{d}, \bm{s})=  \max \quad&\sum_{i=1}^{n_d} \sum_{j=1}^{n_s} v_{i,j}^m x_{i,j} &  \\
                   \text{s.t.} \quad&   \sum_{j=1}^{n_s} x_{i,j} \leq d_i  & \forall i \in [n_d], \label{seq:demand-constraints-CE-with-cost}\\
     &        \sum_{i=1}^{n_d} x_{i,j} \leq s_j  & \forall j \in [n_s], \label{seq:supply-constraints-CE-with-cost}\\
      &       x_{i,j}\geq 0  & \forall i \in [n_d] , \; \forall j  \in [n_s], 
\end{align}
\end{subequations}
while the dual takes the form
\begin{subequations} 
\begin{align}
    \label{lp:dual-CE-with-cost}
    \bar{\Phi}^{m}_{\text{CE}}(\bm{d}, \bm{s})=  \min \quad&\sum_{i=1}^{n_d} a_i^{m} d_i + \sum_{j=1}^{n_s} b_j^{m} s_j &  \\
                   \text{s.t.} \quad&  a_i^{m}  +  b_j^{m} \geq v_{i,j}^m  & \forall i \in [n_d], \, j\in[n_s] \label{seq:dual-constraints-CE-with-cost}\\
     &        a_i^{m}\geq 0 & \forall i \in [n_d], \\
      &       b_j^{m}\geq 0  & \forall j \in [n_s]. 
\end{align}
\end{subequations}

\proof{Proof of Theorem \ref{thm:GTE:RCT-CE-overestimates-GTE}.} 

Consider the optimal solution $x_{i,j}^{\ex,\text{CE}}$ of LP~\eqref{lp:stochastic-primal} (which is also the optimal solution of~\eqref{lp:CE-with-cost}) in the experiment state where the total demand arrival rate in the fluid limit is $\bm\lambda + \rho\bm\beta$. Let us define $y_{i,j}$ as
\[
y_{i,j} = \frac{\lambda_i}{\lambda_i+\rho\beta_i}x_{i,j}^{\ex,\text{CE}}.
\]
Clearly, $y_{i,j}$ is feasible in both ${\Phi}^m_{\text{CE}}(\bm\lambda,\bm\gamma)$ and $\bar{\Phi}^m_{\text{CE}}(\bm\lambda,\bm\gamma)$ since
$$
\sum_{j=1}^{n_s} \frac{\lambda_i}{\lambda_i+\rho\beta_i}x_{i,j}^{\ex,\text{CE}}
=\frac{\lambda_i}{\lambda_i+\rho\beta_i} \sum_{j=1}^{n_s}x_{i,j}^{\ex,\text{CE}} 
\leq
\frac{\lambda_i}{\lambda_i+\rho\beta_i} (\lambda_i+\rho\beta_i)=\lambda_i,
$$
and 
$$
\sum_{i=1}^{n_d} \underbrace{\frac{\lambda_i}{\lambda_i+\rho\beta_i}}_{\le 1}  x_{i,j}^{\ex,\text{CE}} \leq \sum_{i=1}^{n_d} x_{i,j}^{\ex,\text{CE}} \le \gamma_j,
$$
so we can write:
\[
\bar{\Phi}^m_{\text{CE}}(\bm\lambda,\bm\gamma) \ge \sum_{i=1}^{n_d}\sum_{j=1}^{n_s}v_{i,j}^my_{i,j} 
= \sum_{i=1}^{n_d}\lambda_i \underbrace{\frac{\sum_{j=1}^{n_s}v_{i,j}^m x_{i,j}^{\ex,\text{CE}}}{\lambda_i + \rho\beta_i}}_{\bar{v}^m_i} = \sum_{i=1}^{n_d}\bar{v}^m_i\lambda_i,
\]
where $\bar{v}^m_i$ designates the average per-unit value obtained from demand type $i$ under matching $\bar{\Phi}^m_{\text{CE}}(\bm\lambda + \rho\bm\beta,\bm\gamma)$. By Lemma~\ref{lem:CE-unique}, we know that $\bar{v}_i^{\text{prop}}=(1-\alpha)\bar{v}_i^{\ex,\text{CE}}$, while $\bar{v}_i^{\text{fixed}}=\bar{v}_i^{\ex,\text{CE}}-\kappa\bar{x}_i^{\ex,\text{CE}}$. Similarly, we can write
\[
{\Phi}^m_{\text{CE}}(\bm\lambda,\bm\gamma) \ge  \sum_{i=1}^{n_d}\bar{v}_i^{\ex,\text{CE}}\lambda_i,
\]

We also know that by definition, $\bar{\Phi}_{\text{CE}}^m(\bm\lambda+\rho\bm\beta,\bm\gamma) = \sum_{i=1}^{n_d}\bar{v}^m_i(\lambda_i+\rho\beta_i)$. By concavity of $\bar{\Phi}^m(\cdot,\cdot)$, we can conclude
\[
\bar{\Phi}^m_{\text{CE}}(\bm\lambda+\bm\beta,\bm\gamma) \le \sum_{i=1}^{n_d}\bar{v}^m_i(\lambda_i + \beta_i).
\]

Finally, we can write:
\[
\Delta^m = \bar{\Phi}_{\text{CE}}^m(\bm\lambda+\bm\beta,\bm\gamma) - \Phi_{\text{CE}}(\bm\lambda,\bm\gamma) \le \sum_{i=1}^{n_d}\bar{v}^m_i(\lambda_i + \beta_i)-\sum_{i=1}^{n_d}\bar{v}_i^{\ex,\text{CE}}\lambda_i = \hat{\Delta}^m_{\text{RCT-CE}}.
\]

\Halmos
\endproof

\proof{Proof of Lemma~\ref{lem:fluid-limit-rct-ci}}
From Definition \ref{est:RCT-CI}  we know:
$$
\hat \Delta_{\text{{RCT-CI}}}^{\tau, m} 
=
\frac{1}{\tau} \left(  \frac{1}{\rho}\sum_{i=1}^{n_d} \sum_{j=1}^{n_s} v_{i,j}^{\treatment, m} X_{i,j}^{\tau, \treatment, \text{CI}, m}  - \frac{1}{1-\rho} \sum_{i=1}^{n_d} \sum_{j=1}^{n_s} v_{i,j}^{\control}  X_{i,j}^{\tau, \control, \text{CI}, m}  \right).
$$
Define $\bar{V}_i^{\tau,\treatment,\text{CI},m}$ as
\[
\bar{V}_i^{\tau,\treatment,\text{CI},m}=\frac{\sum_{j=1}^{n_s}   v_{i,j}^{\treatment, m} X_{i,j}^{\tau, \treatment, \text{CI}, m}  }{ D_i^{\tau, \treatment} }, 
\]
and analogously define $\bar{V}_i^{\tau,\control,\text{CI},m}$.

Thus,
$$
\hat{\Delta}^{\tau, m}_{\text{{RCT-CI}}} 
= 
\frac{1}{\tau} 
\left(   
\frac{1}{\rho}  \bar{\bm{V}}^{\tau, \treatment, \text{CI}, m} \cdot \bm{D}^{\tau, \treatment} 
- 
\frac{1}{1-\rho} \bar{\bm{V}}^{\tau, \control, \text{CI}, m} \cdot \bm{D}^{\tau, \control}
\right).
$$
By Proposition \ref{prop:fluid-limit} we can conclude
$$
 \lim_{\tau \rightarrow \infty} \hat \Delta_{\text{{RCT-CI}}}^{\tau, m} 
= 
 \hat \Delta_{\text{{RCT-CI}}}^m =
\frac{1}{\rho}  \rho ( \bm{\lambda} + \bm{\beta} ) \cdot \bar{\bm{v}}^{\treatment, \text{CI}, m}  
- \frac{1}{1-\rho} (1-\rho)   \bm{\lambda} \cdot \bar{\bm{v}}^{\control, \text{CI}, m} 
$$
$$
\hat\Delta_{\text{{RCT-CI}}}^m=\bar{\bm{v}}^{\treatment, \text{CI}, m}\cdot (\bm\lambda+\bm\beta) - 
\bar{\bm{v}}^{\control, \text{CI}, m}\cdot \bm\lambda.
$$
\Halmos
\endproof

\proof{Proof of Theorem \ref{thm:GTE:Bias-RCT-CI-can-be-high}.}   

\textbf{Part 1: first and last result.} Recall that in the experiment state, the cost-included (CI) matching LP~\eqref{lp:stochastic-primal-coupondiscount} gives the platform a total value $\Phi_{CI}^m ((1-\rho)\bm\lambda, \rho(\bm\lambda+\bm\beta) ,\bm\gamma)$. In global control, the total value equals $\Phi_{CI}^m ( \bm\lambda,\bm{0} , \bm\gamma)$, while in global treatment it is $\Phi_{CI}^m ( \bm{0},\bm\lambda+\bm\beta , \bm\gamma)$. 

As stated in the main text, let $\Gamma$ designate the sum of all incoming supply, i.e., $  \Gamma = \sum_{j=1}^{n_s} \gamma_j$. Fix $\bm\lambda$ and $\rho$, and let $\Gamma < \min_{i\in[n_d]} (1-\rho)\lambda_i$. Let $x_{i,j}^{\treatment, \text{CI}, m}(\eta)$ (resp. $x_{i,j}^{\control, \text{CI}, m}(\eta)$) designate the optimal treatment (resp. control) matching variables for treatment fraction $\eta$.

In both the experiment and global control states, every supply unit is matched with its highest-value demand unit, which must be a control unit. So we know $x_{i,j}^{\treatment, \text{CI}, m}(\rho)=0$ for all $i$ and $j$. Moreover, for each supply type $j \in[n_s]$ we have that 
$$
x_{i,j}^{\control, \text{CI}, m}(\rho)=x_{i,j}^{\control, \text{CI}, m}(0) = 
\begin{cases}
    \gamma_j, & \text{ if } v_{i,j} > v_{k,j} \; ,\, \forall k \in [n_d] \setminus \{ i \}, \\
    0,        & \text{ otherwise. } 
\end{cases}
$$

Finally, in global treatment, every supply unit is still matched with its highest-value demand unit, but this time it must be the treated version of that unit, since there are no control units. In other words,
$$
x_{i,j}^{\treatment, \text{CI}, m}(1) = 
\begin{cases}
    \gamma_j, & \text{ if } v_{i,j} > v_{k,j} \; ,\, \forall k \in [n_d] \setminus \{ i \}, \\
    0,        & \text{ otherwise. } 
\end{cases}
$$

For each supply type $j$, let $i(j)$ designate the demand type that verifies $v_{i(j),j}>v_{k,j}$ for every $k\neq i(j)$. Then we can write the global treatment effect as
\[
\Delta^{m} = \Phi_{CI}^m ( \bm{0}, \bm\lambda+\bm\beta ,\bm\gamma) - \Phi_{CI}^m (\bm\lambda,\bm{0}, \bm\gamma)=\sum_{j=1}^{n_s}v^m_{i(j),j}\gamma_j-\sum_{j=1}^{n_s}v_{i(j),j}\gamma_j,
\]
and more specifically, $\Delta^{\text{prop}}=-\alpha \sum_{j=1}^{n_s}v_{i(j),j}\gamma_j$, and $\Delta^{\text{fixed}}=-\kappa\Gamma$.

At the experiment point $\eta=\rho$ we get that 
$$
\bar{v}_i^{\treatment, \text{CI}, m}(\rho) = 0 
\quad \text{ and }\quad 
\bar{v}_i^{\control, \text{CI}, m}(\rho)=\frac{\sum\limits_{j=1; i=i(j)}^{n_s} v_{i(j),j}  \gamma_j }{(1-\rho)\lambda_i}.
$$
Therefore, 
$$
\hat\Delta_{\text{RCT-CI}}^m = - \sum_{i=1}^{n_d} \bar{v}_i^{\control, \text{CI}, m}(\rho) \lambda_i = - \frac{1}{1-\rho} \sum_{j=1}^{n_s}v_{i(j),j}\gamma_j.
$$
The signed bias for the proportional cost becomes 
$$
\hat\Delta_{RCT-CI}^{\text{prop}}-\Delta^{\text{prop}} = - \frac{1}{1-\rho} \sum_{j=1}^{n_s}v_{i(j),j}\gamma_j + \alpha \sum_{j=1}^{n_s}v_{i(j),j}\gamma_j
= - \left(\frac{1}{1-\rho}-\alpha\right) \sum_{j=1}^{n_s}v_{i(j),j}\gamma_j,
$$
which is clearly negative since $\alpha \le 1$. The relative bias equals 
$$
\frac{\hat\Delta_{RCT-CI}^{\text{prop}}-\Delta^{\text{prop}}}{\left| \Delta^{\text{prop}}  \right|} =
 \frac{  - (\frac{1}{1-\rho}-\alpha) }{\left| - \alpha  \right|}  
= 1 - \frac{1}{\alpha(1-\rho)},
$$
which tends to $-\infty$ as $\alpha\to 0^+$.
Similarly, for the fixed-cost model we get
$$
\hat\Delta_{\text{RCT-CI}}^{\text{fixed}}-\Delta^{\text{fixed}} = - \frac{1}{1-\rho} \sum_{j=1}^{n_s}v_{i(j),j}\gamma_j 
+ \kappa\Gamma,
$$
which is negative since $\kappa < \min_{i,j}v_{i,j}$. The relative bias is given by
$$
\frac{\hat\Delta_{\text{RCT-CI}}^{\text{fixed}}(\rho)-\Delta^{\text{fixed}}}{\left| \Delta^{\text{fixed}}  \right|} =
 \frac{ - \frac{1}{1-\rho} \sum_{j=1}^{n_s}v_{i(j),j}\gamma_j   +\kappa\Gamma }{\left| - \kappa\Gamma \right|}  
=   1- \frac{1}{1-\rho} \underbrace{\frac{\sum_{j=1}^{n_s}v_{i(j),j}\gamma_j}{\kappa\Gamma} }_{>1},
$$
which again tends to $-\infty$ as $\kappa \rightarrow 0^+$.

\textbf{Part 2: second result.} Fix $\bm\lambda$, $\bm\beta$, and fix $\bm\gamma$ up to the scaling factor $\Gamma$. Let $\Gamma_m$ be the smallest value of $\Gamma$ such that in the global treatment state, every demand unit is matched to its highest-value supply unit, i.e.,
\[
x_{i,j}^{\treatment,\text{CI},m}(1) = \begin{cases}
    \lambda_i + \beta_i, & \text{ if } v_{i,j} > v_{i,l}, \forall l \in [n_s] \setminus \{ j \},  \\
    0,        & \text{ otherwise. }
\end{cases}
\]
For each demand type $i$, let $j(i)$ designate the highest-value matching supply type, i.e., verifying $v_{i,j(i)}>v_{i,l}$ for all $l\neq j(i)$.

We define $\Gamma_0=(1-\varepsilon)\Gamma_m$, with very small $\varepsilon>0$. By definition of $\Gamma_m$, it is no longer possible for every demand type to be matched with its top supply type, i.e., there exists some $i$ such that $x_{i,j(i)}^{\treatment,\text{CI},m}(1) < \lambda_i + \beta_i$. We choose $\varepsilon$ small enough that there is only one such type, which we denote by $i^*$, and also so that it does not affect the experiment or global control states. In the global control state, we can therefore write:
\[
x_{i,j}^{\control,\text{CI},m}(0) = \begin{cases}
    \lambda_i, & \text{ if } v_{i,j} > v_{i,l}, \forall l \in [n_s] \setminus \{ j \},  \\
    0,        & \text{ otherwise, }
\end{cases}
\]
and similarly in the experiment state:
\begin{align*}
x_{i,j}^{\control,\text{CI},m}(\rho) &= \begin{cases}
    (1-\rho)\lambda_i, & \text{ if } v_{i,j} > v_{i,l}, \forall l \in [n_s] \setminus \{ j \},  \\
    0,        & \text{ otherwise. }
    \end{cases}\\
    x_{i,j}^{\treatment,\text{CI},m}(\rho) &= \begin{cases}
    \rho(\lambda_i + \beta_i), & \text{ if } v_{i,j} > v_{i,l}, \forall l \in [n_s] \setminus \{ j \},  \\
    0,        & \text{ otherwise. }
\end{cases}
\end{align*}

We can then compute the following quantities:
\begin{itemize}
    \item Total value in global control:
    \[
    \Phi^m_{\text{CI}}(\bm{\lambda},\bm{0},\bm\gamma) = \sum_{i=1}^{n_d}v_{i,j(i)}\lambda_i.
    \]
    \item Cost-included RCT estimator:
    \[
    \hat{\Delta}_{\text{RCT-CI}}^m = \sum_{i=1}^{n_d}v^m_{i,j(i)}(\lambda_i + \beta_i) - v_{i,j(i)}\lambda_i.
    \]
    \item Total value in global treatment:
    \[
    \Phi^m_{\text{CI}}(\bm{0},\bm{\lambda}+\bm\beta,\bm\gamma) = \sum_{i=1;i\neq i^*}^{n_d}v^m_{i,j(i)}(\lambda_i + \beta_i) + v^m_{i^*,j(i^*)}x_{i^*,j(i^*)}^{\treatment,\text{CI},m}(1) + \sum_{j=1;j \neq j(i^*)}^{n_s}v^m_{i^*,j}x_{i^*,j}^{\treatment,\text{CI},m}(1).
    \]
\end{itemize}

Therefore, we can write the signed bias of the RCT estimator as:
\begin{align*}
    \hat{\Delta}_{\text{RCT-CI}}^m - \Delta^m &= \sum_{i=1}^{n_d}v^m_{i,j(i)}(\lambda_i + \beta_i) - v_{i,j(i)}\lambda_i - \sum_{i=1;i\neq i^*}^{n_d}v^m_{i,j(i)}(\lambda_i + \beta_i) \\&- v^m_{i^*,j(i^*)}x_{i^*,j(i^*)}^{\treatment,\text{CI},m}(1) - \sum_{j=1;j \neq j(i^*)}^{n_s}v^m_{i^*,j}x_{i^*,j}^{\treatment,\text{CI},m}(1) + \sum_{i=1}^{n_d}v_{i,j(i)}\lambda_i\\
    &= v^m_{i^*,j(i^*)}(\lambda_i + \beta_i - x_{i^*,j(i^*)}^{\treatment,\text{CI},m}(1)) - \sum_{j=1;j \neq j(i^*)}^{n_s}v^m_{i^*,j}x_{i^*,j}^{\treatment,\text{CI},m}(1)\\
    &> v^m_{i^*,j(i^*)}(\lambda_i + \beta_i - x_{i^*,j(i^*)}^{\treatment,\text{CI},m}(1)) - v^m_{i^*,j(i^*)}\sum_{j=1;j\neq j(i^*)}^{n_s}x_{i^*,j}^{\treatment,\text{CI},m}(1)\\
    &\ge v^m_{i^*,j(i^*)}\left(\lambda_i + \beta_i - \sum_{j=1}^{n_s}x_{i^*,j}^{\treatment,\text{CI},m}(1)\right)\ge 0.
\end{align*}

Note that the strict inequality holds because by construction $x_{i^*,j(i^*)}^{\treatment,\text{CI},m}(1) < \lambda_i + \beta_i$.
\\
\Halmos
\endproof

\section{SP Estimator Analysis}

\proof{Proof of Lemma~\ref{lem:sp-ce-fluid-limit}.}
The lemma follows directly from Proposition~\ref{prop:fluid-limit} and the strong law of large numbers. In particular, the convergence of $(\bm{\tilde{A}}^{\tau,m}, \bm{\tilde{B}}^{\tau,m})$ to $(\bm{\tilde{a}}^{m}, \bm{\tilde{b}}^m)$ follows directly from Proposition~\ref{prop:fluid-limit} and Lemma~\ref{lem:fixed-cost-shadow-price}.
\Halmos
\endproof

\proof{Proof of Theorem \ref{thm:RCT-CE-vs-SP-CE}.}
We first introduce the following useful form of the global treatment effect:
\begin{align*}
    \Delta^m &= \bar{\Phi}^m_{\text{CE}}(\bm\lambda+\bm\beta,\bm\gamma) - \Phi_{\text{CE}}(\bm\lambda,\bm\gamma)\\
    &=\bar{\Phi}^m_{\text{CE}}(\bm\lambda+\bm\beta,\bm\gamma) - \bar{\Phi}^m_{\text{CE}}(\bm\lambda+\rho\bm\beta,\bm\gamma) + \bar{\Phi}^m_{\text{CE}}(\bm\lambda+\rho\bm\beta,\bm\gamma)  - \Phi_{\text{CE}}(\bm\lambda+\rho\bm\beta,\bm\gamma)+ \Phi_{\text{CE}}(\bm\lambda+\rho\bm\beta,\bm\gamma) - \Phi_{\text{CE}}(\bm\lambda,\bm\gamma)\\
    &= \int_\rho^1 \bm{\tilde{a}}^{m,\eta}\cdot\bm\beta d\eta + \bar{\Phi}^m_{\text{CE}}(\bm\lambda+\rho\bm\beta,\bm\gamma)  - \Phi_{\text{CE}}(\bm\lambda+\rho\bm\beta,\bm\gamma) + \int_0^\rho \bm{a}^{\eta}\cdot\bm\beta d\eta .
\end{align*}
Similarly, recalling that $\bar{\bm{v}}^{\ex, m}$ is defined as 
\begin{align*}
    \bar{\bm{v}}^{\ex, \text{prop}}&= (1-\alpha)\bar{\bm{v}}^{\ex, \text{{CE}}}\\
    \bar{\bm{v}}^{\ex, \text{fixed}}&= \bar{\bm{v}}^{\ex, \text{{CE}}} - \kappa\bar{\bm{x}}^{\ex, \text{{CE}}},
\end{align*}
we can write
\begin{align*}
    \hat{\Delta}^m_{\text{RCT-CE}} &= \bar{\bm{v}}^{\ex, m} \cdot (\bm\lambda+\bm\beta) - \bar{\bm{v}}^{\ex, \text{{CE}}} \cdot \bm\lambda\\
    &=\bar{\bm{v}}^{\ex, m} \cdot (\bm\lambda+\bm\beta) - \bar{\bm{v}}^{\ex, m} \cdot (\bm\lambda+\rho\bm\beta) + \bar{\bm{v}}^{\ex, m} \cdot (\bm\lambda+\rho\bm\beta) - \bar{\bm{v}}^{\ex, \text{CE}} \cdot (\bm\lambda+\rho\bm\beta) + \bar{\bm{v}}^{\ex, \text{CE}} \cdot (\bm\lambda+\rho\bm\beta) - \bar{\bm{v}}^{\ex, \text{{CE}}} \cdot \bm\lambda\\
    &= (1-\rho)\bar{\bm{v}}^{\ex, m} \cdot \bm\beta + \bar{\Phi}^m_{\text{CE}}(\bm\lambda+\rho\bm\beta,\bm\gamma)  - \Phi_{\text{CE}}(\bm\lambda+\rho\bm\beta,\bm\gamma) + \rho\bar{\bm{v}}^{\ex, \text{{CE}}} \cdot \bm\beta.
\end{align*}

And analogously,
\begin{align*}
    \hat{\Delta}^m_{\text{SP-CE}}&=\left(\bm{\tilde{a}}^{m}\cdot(\bm\lambda + \bm\beta) + \bm{\tilde{b}}^{m}\cdot \bm\gamma\right) - \left(\bm{a}^{\text{\emph{CE}}}\cdot \bm\lambda + \bm{b}^{\text{\emph{CE}}}\cdot\bm\gamma \right)\\
    &=\left(\bm{\tilde{a}}^{m}\cdot(\bm\lambda + \bm\beta) + \bm{\tilde{b}}^{m}\cdot \bm\gamma\right) - \left(\bm{\tilde{a}}^{m}\cdot(\bm\lambda + \rho\bm\beta) + \bm{\tilde{b}}^{m}\cdot \bm\gamma\right) + \left(\bm{\tilde{a}}^{m}\cdot(\bm\lambda + \rho\bm\beta) + \bm{\tilde{b}}^{m}\cdot \bm\gamma\right) - \left(\bm{a}^{\text{\emph{CE}}}\cdot \bm\lambda + \bm{b}^{\text{\emph{CE}}}\cdot\bm\gamma \right)\\
    &= (1 - \rho)\bm{\tilde{a}}^{m}\cdot\bm\beta + \bar{\Phi}^m_{\text{CE}}(\bm\lambda+\rho\bm\beta,\bm\gamma) - \left(\bm{a}^{\text{\emph{CE}}}\cdot (\bm\lambda+\rho\bm\beta) + \bm{b}^{\text{\emph{CE}}}\cdot\bm\gamma \right) + \left(\bm{a}^{\text{\emph{CE}}}\cdot (\bm\lambda+\rho\bm\beta) + \bm{b}^{\text{\emph{CE}}}\cdot\bm\gamma \right) - \left(\bm{a}^{\text{\emph{CE}}}\cdot \bm\lambda + \bm{b}^{\text{\emph{CE}}}\cdot\bm\gamma \right)\\
    &= (1 - \rho)\bm{\tilde{a}}^{m}\cdot\bm\beta + \bar{\Phi}^m_{\text{CE}}(\bm\lambda+\rho\bm\beta,\bm\gamma) - \Phi_{\text{CE}}(\bm\lambda+\rho\bm\beta,\bm\gamma) + \rho \bm{a}^{\text{\emph{CE}}}\cdot\bm\beta.
\end{align*}

Recalling that $\bm{\bar{v}}^{\ex,m}\ge \bm{\tilde{a}}^m$, we can bound the bias of the RCT estimator as
\begin{align}
\hat{\Delta}^m_{\text{RCT-CE}} - \Delta^m &= (1-\rho)\bar{\bm{v}}^{\ex, m} \cdot \bm\beta + \rho\bar{\bm{v}}^{\ex, \text{{CE}}} \cdot \bm\beta - \int_\rho^1 \bm{\tilde{a}}^{m,\eta}\cdot\bm\beta d\eta - \int_0^\rho \bm{a}^{\eta}\cdot\bm\beta d\eta\nonumber\\
&=\int_\rho^1(\bar{\bm{v}}^{\ex, m}-\bm{\tilde{a}}^{m,\eta})\cdot\bm\beta d\eta + \underbrace{\rho\bar{\bm{v}}^{\ex, \text{{CE}}} \cdot \bm\beta - \int_0^\rho \bm{a}^{\eta}\cdot\bm\beta d\eta}_{\ge 0\text{ (Theorem~\ref{thm:GTE:RCT-CE-overestimates-GTE})}}\nonumber\\
&\ge \int_\rho^1(\bar{\bm{v}}^{\ex, m}-\bm{\tilde{a}}^m+\bm{\tilde{a}}^m-\bm{\tilde{a}}^{m,\eta})\cdot\bm\beta d\eta\nonumber\\
&=(1-\rho)(\bar{\bm{v}}^{\ex, m}-\bm{\tilde{a}}^m)\cdot\bm\beta + \int_\rho^1(\bm{\tilde{a}}^m-\bm{\tilde{a}}^{m,\eta})\cdot\bm\beta d\eta.\label{eq:intermediate-step-sp-ce-reduce-bias}
\end{align}

Separately, we know from concavity of the matching value function and from Theorem~\ref{thm:GTE:RCT-CE-overestimates-GTE} that
\[
    \bar{\bm{v}}^{\ex, \text{{CE}}}\cdot\bm\lambda \le \bar{\Phi}_{\text{CE}}(\bm\lambda,\bm\gamma) \le \bar{\Phi}_{\text{CE}}(\bm\lambda+\rho\bm\beta,\bm\gamma) - \rho \bm{a}^{\text{CE}}\cdot\bm\beta,
\]
therefore
\begin{align}
    \bar{\Phi}_{\text{CE}}(\bm\lambda+\rho\bm\beta,\bm\gamma) - \rho \bm{a}^{\text{CE}}\cdot\bm\beta - \bar{\Phi}_{\text{CE}}(\bm\lambda,\bm\gamma) &\le \bar{\bm{v}}^{\ex, \text{{CE}}}\cdot(\bm\lambda+\rho\bm\beta)- \rho \bm{a}^{\text{CE}}\cdot\bm\beta - \bar{\bm{v}}^{\ex, \text{{CE}}}\cdot\bm\lambda\nonumber\\
    \int_0^\rho (\bm{a}^\eta -\bm{a}^{\text{CE}})\cdot\bm\beta d\eta &\le \rho (\bar{\bm{v}}^{\ex, \text{{CE}}} - \bm{a}^{\text{CE}})\cdot\bm\beta\label{eq:intermediate-step-sp-ce-reduce-bias-2}
\end{align}

Now, we are ready to upper-bound the bias of the SP estimator, as follows:
\begin{align}
    \hat{\Delta}^m_{\text{SP-CE}} - \Delta^m &= (1 - \rho)\bm{\tilde{a}}^{m}\cdot\bm\beta + \rho \bm{a}^{\text{\emph{CE}}}\cdot\bm\beta - \int_\rho^1 \bm{\tilde{a}}^{m,\eta}\cdot\bm\beta d\eta - \int_0^\rho \bm{a}^{\eta}\cdot\bm\beta d\eta \nonumber\\
    &= \int_0^\rho (\bm{a}^{\text{\emph{CE}}} - \bm{a}^{\eta}) \cdot\bm\beta d\eta + \int_\rho^1 (\bm{\tilde{a}}^{m} - \bm{\tilde{a}}^{m,\eta})\cdot\bm\beta d\eta\nonumber.
    \end{align}
Using the triangle inequality, we can write
    \begin{align}
    \abs{\hat{\Delta}^m_{\text{SP-CE}} - \Delta^m}&\le \abs{\int_0^\rho (\bm{a}^{\text{\emph{CE}}} - \bm{a}^{\eta}) \cdot\bm\beta d\eta} + \abs{\int_\rho^1 (\bm{\tilde{a}}^{m} - \bm{\tilde{a}}^{m,\eta})\cdot\bm\beta d\eta}\nonumber\\
    &=\int_0^\rho (\bm{a}^{\eta}-\bm{a}^{\text{\emph{CE}}}) \cdot\bm\beta d\eta + \int_\rho^1 (\bm{\tilde{a}}^{m} - \bm{\tilde{a}}^{m,\eta})\cdot\bm\beta d\eta\nonumber\\
    &\le \rho (\bar{\bm{v}}^{\ex, \text{{CE}}} - \bm{a}^{\text{CE}})\cdot\bm\beta + \int_\rho^1 (\bm{\tilde{a}}^{m} - \bm{\tilde{a}}^{m,\eta})\cdot\bm\beta d\eta,\label{eq:intermediate-step-sp-ce-reduce-bias-3}
\end{align}
where the last step follows from~\eqref{eq:intermediate-step-sp-ce-reduce-bias-2}. Putting together~\eqref{eq:intermediate-step-sp-ce-reduce-bias} and \eqref{eq:intermediate-step-sp-ce-reduce-bias-3}, we obtain a sufficient condition for $\hat{\Delta}^m_{\text{RCT-CE}} - \Delta^m \ge\abs{\hat{\Delta}^m_{\text{SP-CE}} - \Delta^m} $ as
\[
\rho (\bar{\bm{v}}^{\ex, \text{{CE}}} - \bm{a}^{\text{CE}})\cdot\bm\beta \le (1-\rho)(\bar{\bm{v}}^{\ex, m}-\bm{\tilde{a}}^m)\cdot\bm\beta.
\]

Under the proportional-cost model, this reduces to
\begin{align*}
    \rho & \le (1-\rho)\frac{(1-\alpha)(\bar{\bm{v}}^{\ex, \text{{CE}}} - \bm{a}^{\text{CE}})\cdot\bm\beta}{(\bar{\bm{v}}^{\ex, \text{{CE}}} - \bm{a}^{\text{CE}})\cdot\bm\beta}\\
    \rho &\le \frac{1-\alpha}{2-\alpha}.
\end{align*}

Under the fixed-cost model, there are two cases. If $\sum_{i=1}^{n_d}\lambda_i + \rho\beta_i < \sum_{j=1}^{n_s}\gamma_j$, then $\bar{\bm{x}}^{\ex,\text{CE}}=\bm{e}$ and the condition becomes
\begin{align*}
    \rho &\le (1-\rho)\frac{(\bar{\bm{v}}^{\ex, \text{{CE}}} - \kappa\bm{e} - \bm{a}^{\text{CE}}+\kappa\bm{e})\cdot\bm\beta}{(\bar{\bm{v}}^{\ex, \text{{CE}}} - \bm{a}^{\text{CE}})\cdot\bm\beta}\\
    \rho &\le \frac{1}{2}.
\end{align*}
Otherwise, if $\sum_{i=1}^{n_d}\lambda_i + \rho\beta_i \ge \sum_{j=1}^{n_s}\gamma_j$, the condition becomes
\begin{align*}
    \rho &\le (1-\rho)\frac{(\bar{\bm{v}}^{\ex, \text{{CE}}} - \kappa\bar{\bm{x}}^{\ex,\text{CE}} - \bm{a}^{\text{CE}})\cdot\bm\beta}{(\bar{\bm{v}}^{\ex, \text{{CE}}} - \bm{a}^{\text{CE}})\cdot\bm\beta}\\
    \rho &\le (1-\rho)\left(1-\frac{\kappa\bar{\bm{x}}^{\ex,\text{CE}}\cdot\bm\beta}{(\bar{\bm{v}}^{\ex, \text{{CE}}} - \bm{a}^{\text{CE}})\cdot\bm\beta}\right).
\end{align*}
This last condition is useful, but not computable before running the experiment (though it can be verified after the fact). For ex-ante bias reduction guarantees, we observe that any demand unit of type $i$ can only be matched if another demand unit of another type is unmatched (because all supply units are currently matched). Therefore, we can write $a_i^{\text{CE}}\le \bar{v}_i^{\ex,\text{CE}}-\min_{i,j}v_{i,j}$ (the maximum value we will get from the additional unit is the average value from all units of this type, minus the smallest possible value we might lose). Therefore, we can write
\[
\frac{\kappa\bar{\bm{x}}^{\ex,\text{CE}}\cdot\bm\beta}{(\bar{\bm{v}}^{\ex, \text{{CE}}} - \bm{a}^{\text{CE}})\cdot\bm\beta} \le \frac{\kappa\bm{e}\cdot\bm\beta}{\min_{i,j}v_{i,j}\bm{e}\cdot\bm\beta}=\frac{\kappa}{\min_{i,j}v_{i,j}},
\]
leading to the final condition
\[
\rho \le \frac{1-\frac{\kappa}{\min_{i,j}v_{i,j}}}{2 - \frac{\kappa}{\min_{i,j}v_{i,j}}} \le \frac{1}{2}.
\]

We can prove that this bound is tight by constructing the following instance for any $\alpha$ under the proportional-cost model. (It turns out that the exact same construction applies to the fixed-cost model, with $\kappa=\alpha$).

Consider an instance with a single demand and supply type, $\lambda=0$, $\beta=1$, $v_{1,1}=1$, and $\gamma = (1-\alpha)/(2-\alpha)$. The true effect is given by
\[
\Delta^{\text{prop}} = (1-\alpha)v_{1,1}\gamma - 0 = \frac{(1-\alpha)^2}{2-\alpha}.
\]
For any $\rho>(1-\alpha)/(2-\alpha)$, the shadow price estimator will always be $0$, meaning the absolute bias of the shadow price estimator will be
\[
\abs{\hat{\Delta}^{\text{prop}}_{\text{SP-CE}} - \Delta^{\text{prop}}} = \frac{(1-\alpha)^2}{2-\alpha}.
\]
Meanwhile the standard estimator will be upper-bounded by its value at $\rho=(1-\alpha)/(2-\alpha)$, i.e.,
\[
\hat{\Delta}^{\text{prop}}_{\text{RCT-CE}} \le (1-\alpha),
\]
meaning the absolute bias is upper-bounded by
\[
\abs{\hat{\Delta}^{\text{prop}}_{\text{RCT-CE}} - \Delta^{\text{prop}}} \le 1-\alpha-\frac{(1-\alpha)^2}{2-\alpha}=\frac{(1-\alpha)^2}{2-\alpha} = \abs{\hat{\Delta}^{\text{prop}}_{\text{SP-CE}} - \Delta^{\text{prop}}}.
\]
\Halmos
\endproof

\proof{Proof of Theorem~\ref{thm:fraction-bias-SP-CE-vs-RCT-CE}.}
In the proof of Theorem~\ref{thm:RCT-CE-vs-SP-CE}, we established the following results:
\begin{align*}
    \hat{\Delta}^m_{\text{RCT-CE}} - \Delta^m &=\int_0^\rho (\bar{\bm{v}}^{\ex, \text{{CE}}}-\bm{a}^{\eta})\cdot\bm\beta d\eta+\int_\rho^1(\bar{\bm{v}}^{\ex, m}-\bm{\tilde{a}}^{m,\eta})\cdot\bm\beta d\eta\\
    \hat{\Delta}^m_{\text{SP-CE}} - \Delta^m &= \int_0^\rho (\bm{a}^{\text{\emph{CE}}} - \bm{a}^{\eta}) \cdot\bm\beta d\eta + \int_\rho^1 (\bm{\tilde{a}}^{m} - \bm{\tilde{a}}^{m,\eta})\cdot\bm\beta d\eta.
\end{align*}

We will also use the concavity of the matching function, which implies $\bm{a}^{\text{CE},0}\cdot\bm\beta \ge \bm{a}^{\text{CE},\eta}\cdot\bm\beta \ge \bm{a}^{\text{CE},1}\cdot\bm\beta$ for any $\eta$, and in particular for $\eta=\rho$, $\bm{a}^{\text{CE},0}\cdot\bm\beta \ge \bm{a}^{\text{CE}}\cdot\bm\beta \ge \bm{a}^{\text{CE},1}\cdot\bm\beta$.

\textbf{Proportional cost model:} We can write (dropping the CE subscripts when clear from context):
\begin{align*}
    \hat{\Delta}^{\text{prop}}_{\text{RCT-CE}} - \Delta^{\text{prop}} &= \int_0^\rho (\bar{\bm{v}}^{\ex, \text{{CE}}}-\bm{a}^{\eta})\cdot\bm\beta d\eta+(1-\alpha)\int_\rho^1(\bar{\bm{v}}^{\ex, \text{CE}}-\bm{a}^{\eta})\cdot\bm\beta d\eta\\
    &\ge \left(\rho + (1-\rho)(1-\alpha)\right)(\bar{\bm{v}}^{\ex, \text{CE}}-\bm{a}^{0})\cdot\bm\beta.
\end{align*}
and similarly,
\begin{align*}
    \hat{\Delta}^{\text{prop}}_{\text{SP-CE}} - \Delta^{\text{prop}} &= \int_0^\rho (\bm{a}^{\text{\emph{CE}}} - \bm{a}^{\eta}) \cdot\bm\beta d\eta + (1-\alpha)\int_\rho^1 (\bm{{a}}^{\text{CE}} - \bm{{a}}^{\eta})\cdot\bm\beta d\eta.\\
    \abs{\hat{\Delta}^{\text{prop}}_{\text{SP-CE}}-\Delta^{\text{prop}}} & \le \left(\rho + (1-\rho)(1-\alpha)\right)(\bm{a}^{0}-\bm{a}^{1})\cdot\bm\beta.
\end{align*}

Putting the two together, as long as $(\bar{\bm{v}}^{\ex, \text{CE}}-\bm{a}^{0})\cdot\bm\beta > 0$, we obtain:
\[
    \frac{\abs{\hat{\Delta}^{\text{prop}}_{\text{SP-CE}} - \Delta^{\text{prop}}}}{\abs{\hat{\Delta}^{\text{prop}}_{\text{RCT-CE}} - \Delta^{\text{prop}}}} \le \frac{\left(\rho + (1-\rho)(1-\alpha)\right)(\bm{a}^{0}-\bm{a}^{1})\cdot\bm\beta}{\left(\rho + (1-\rho)(1-\alpha)\right)(\bar{\bm{v}}^{\ex, \text{CE}}-\bm{a}^{0})\cdot\bm\beta}=\frac{(\bm{a}^{0}-\bm{a}^{1})\cdot\bm\beta}{(\bar{\bm{v}}^{\ex, \text{CE}}-\bm{a}^{0})\cdot\bm\beta}.
\]

\textbf{Fixed-cost case:} First consider the case where $\sum_{i=1}^{n_d}\lambda_i + \beta_i < \sum_{j=1}^{n_s}\gamma_j$. In this case, we can write:
\begin{align*}
    \hat{\Delta}^{\text{fixed}}_{\text{RCT-CE}} - \Delta^{\text{fixed}} &= \int_0^\rho (\bar{\bm{v}}^{\ex, \text{{CE}}}-\bm{a}^{\eta})\cdot\bm\beta d\eta+\int_\rho^1(\bar{\bm{v}}^{\ex, \text{CE}}-\kappa\bm{e}-\bm{a}^{\eta} + \kappa\bm{e})\cdot\bm\beta d\eta\\
    &\ge (\bar{\bm{v}}^{\ex, \text{CE}}-\bm{a}^{0})\cdot\bm\beta.
\end{align*}
And similarly for the shadow price estimator:
\begin{align*}
    \hat{\Delta}^{\text{fixed}}_{\text{SP-CE}} - \Delta^{\text{fixed}} &= \int_0^\rho (\bm{a}^{\text{\emph{CE}}} - \bm{a}^{\eta}) \cdot\bm\beta d\eta + \int_\rho^1 (\bm{{a}}^{\text{CE}} - \kappa\bm{e} - \bm{{a}}^{\eta} + \kappa\bm{e})\cdot\bm\beta d\eta.\\
    \abs{\hat{\Delta}^{\text{fixed}}_{\text{SP-CE}}-\Delta^{\text{fixed}}} & \le (\bm{a}^{0}-\bm{a}^{1})\cdot\bm\beta.
\end{align*}
Assuming a positive denominator, we can conclude once again:
\[
    \frac{\abs{\hat{\Delta}^{\text{fixed}}_{\text{SP-CE}} - \Delta^{\text{fixed}}}}{\abs{\hat{\Delta}^{\text{fixed}}_{\text{RCT-CE}} - \Delta^{\text{fixed}}}} \le\frac{(\bm{a}^{0}-\bm{a}^{1})\cdot\bm\beta}{(\bar{\bm{v}}^{\ex, \text{CE}}-\bm{a}^{0})\cdot\bm\beta}.
\]
Finally, consider the case where $\sum_{i=1}^{n_d}\lambda_i \ge \sum_{j=1}^{n_s}\gamma_j$. The shadow price estimator bias can be bounded as:
\begin{align*}
    \hat{\Delta}^{\text{fixed}}_{\text{SP-CE}} - \Delta^{\text{fixed}} &= \int_0^\rho (\bm{a}^{\text{\emph{CE}}} - \bm{a}^{\eta}) \cdot\bm\beta d\eta + \int_\rho^1 (\bm{{a}}^{\text{CE}} - \bm{{a}}^{\eta})\cdot\bm\beta d\eta.\\
    \abs{\hat{\Delta}^{\text{fixed}}_{\text{SP-CE}}-\Delta^{\text{fixed}}} & \le (\bm{a}^{0}-\bm{a}^{1})\cdot\bm\beta.
\end{align*}
Meanwhile, we can write the RCT estimator bias as
\begin{align*}
    \hat{\Delta}^{\text{fixed}}_{\text{RCT-CE}} - \Delta^{\text{fixed}} &= \int_0^\rho (\bar{\bm{v}}^{\ex, \text{{CE}}}-\bm{a}^{\eta})\cdot\bm\beta d\eta+\int_\rho^1(\bar{\bm{v}}^{\ex, \text{CE}}-\kappa\bm{\bar{x}}^{\ex,\text{CE}}-\bm{a}^{\eta})\cdot\bm\beta d\eta\\
    &\ge (\bar{\bm{v}}^{\ex, \text{CE}}-\bm{a}^{0})\cdot\bm\beta - (1-\rho)\kappa\bm{e}\cdot\bm{\beta},
\end{align*}
so again assuming positivity of the denominator, we obtain:
\[
    \frac{\abs{\hat{\Delta}^{\text{fixed}}_{\text{SP-CE}} - \Delta^{\text{fixed}}}}{\abs{\hat{\Delta}^{\text{fixed}}_{\text{RCT-CE}} - \Delta^{\text{fixed}}}} \le\frac{(\bm{a}^{0}-\bm{a}^{1})\cdot\bm\beta}{(\bar{\bm{v}}^{\ex, \text{CE}}-\bm{a}^{0})\cdot\bm\beta - (1-\rho)\kappa\bm{e}\cdot\bm{\beta}},
\]
and we notice this bound implies the bound in the first case of the fixed-cost model.\\
\Halmos
\endproof

\proof{Proof of Lemma \ref{lem:fluid-limit-sp-ci}.}
The lemma follows directly from the strong law of large numbers and Proposition~\ref{prop:fluid-limit}.
\Halmos
\endproof

\proof{Proof of Theorem~\ref{thm:sp-ci-reduces-bias}.}
We revisit the construction from Theorem~\ref{thm:GTE:Bias-RCT-CI-can-be-high}. First, for $\Gamma < \Gamma_{\min}$, we know that in both the experiment and global control states, every supply unit is matched with its highest-value demand unit, which must be a control unit, and furthermore there is at least one unmatched control unit of each type. So $a_{i}^{\treatment, \text{CI}, m}(\rho)=a_{i}^{\treatment, \text{CI}, m}(\rho)=0$ for all $i$ and $j$. Recall that we can write the global treatment effect as
\[
\Delta^{m} = \Phi_{CI}^m ( \bm{0}, \bm\lambda+\bm\beta ,\bm\gamma) - \Phi_{CI}^m (\bm\lambda,\bm{0}, \bm\gamma)=\sum_{j=1}^{n_s}v^m_{i(j),j}\gamma_j-\sum_{j=1}^{n_s}v_{i(j),j}\gamma_j,
\]
and more specifically, $\Delta^{\text{prop}}=-\alpha \sum_{j=1}^{n_s}v_{i(j),j}\gamma_j$, and $\Delta^{\text{fixed}}=-\kappa\Gamma$. Because all the shadow prices are 0, the shadow price estimator gives $\hat\Delta_{\text{SP-CI}}^m = 0$.

The signed bias for the proportional cost becomes 
$$
\hat\Delta_{\text{SP-CI}}^{\text{prop}}-\Delta^{\text{prop}} = \alpha \sum_{j=1}^{n_s}v_{i(j),j}\gamma_j,
$$
so the relative bias of the shadow price estimator to the standard estimator is
$$
\frac{\abs{\hat\Delta_{\text{SP-CI}}^{\text{prop}}-\Delta^{\text{prop}}}}{\abs{\hat\Delta_{\text{RCT-CI}}^{\text{prop}}-\Delta^{\text{prop}}}} =
 \frac{\alpha}{{\frac{1}{1-\rho}-\alpha}},
$$
which is less than 1 as long as
\[
\rho \ge 1 - \frac{1}{2\alpha},
\]
and tends to zero as $\alpha$ tends to zero.

The signed bias for the fixed cost becomes $\hat\Delta_{\text{SP-CI}}^{\text{fixed}}-\Delta^{\text{fixed}} =\kappa\Gamma$, so the relative bias of the shadow price estimator to the standard estimator is
\[
\frac{\abs{\hat\Delta_{\text{SP-CI}}^{\text{fixed}}-\Delta^{\text{fixed}}}}{\abs{\hat\Delta_{\text{RCT-CI}}^{\text{fixed}}-\Delta^{\text{fixed}}}} =
 \frac{\kappa\Gamma}{\frac{1}{1-\rho} \sum_{j=1}^{n_s}v_{i(j),j}\gamma_j   -\kappa\Gamma},
\]
which again tends to $0$ as $\kappa \rightarrow 0$, and is less than one as long as
\[
\rho \ge 1 - \frac{1}{2\frac{\kappa}{\min_j v_{i(j),j}}}.
\]

The second part of the theorem is easier to show. Given that in the experiment state, every demand type is matched with its highest-value supply type, the average value obtained from each type is equal to its marginal value, and the RCT and SP estimators are equal.
\\
\Halmos
\endproof

\section{Design selection analysis}

To prove Theorem~\ref{thm:ce-vs-ci}, we first need to prove the following lemma using linear programming sensitivity analysis.

\begin{lemma}
\label{lem:sensitivity-analysis}
    For $\bm{b}\in\mathbb{R}^m$, $\bm{A}\in\mathbb{R}^{m\times n}$ and $\bm{c}\in\mathbb{R}^n$, define $P(\bm{b})=\{\bm{x}\in\mathbb{R}^n: \bm{A}\bm{x}=\bm{b},\bm{x}\ge 0\}$, and $F(\bm{b})=\min_{\bm{x}\in P(\bm{b})}\bm{c}\cdot\bm{x}$. If $F(\bm{b})$ is linear on the line segment from $\bm{b}_1$ to $\bm{b}_2$, i.e., $F(\lambda\bm{b}_1+(1-\lambda)\bm{b}_2)=\lambda F(\bm{b}_1) + (1-\lambda) F(\bm{b}_2)$ for $\lambda\in[0,1]$, then $F(\bm{b}_1)$ and $F(\bm{b}_2)$ have the same optimal dual solution.
\end{lemma}

\proof{Proof of Lemma~\ref{lem:sensitivity-analysis}.}
We can write the dual of the problem of interest as
\begin{align*}
    F(\bm{b})=\max\quad & \bm{p}\cdot\bm{b}\\
    \text{{s.t.}}\quad & \bm{A}^{\top}\bm{p}\le\bm{c}.
\end{align*}

Assume without loss of generality that the optimal primal solutions at $F(\bm{b}_1)$ and $F(\bm{b}_2)$ are both nondegenerate. Then we know that (1) the optimal dual solutions, denoted by $\bm{p}_1$ and $\bm{p}_2$, are unique; and (2), the function is locally linear around $\bm{b}_1$ and $\bm{b}_2$. Assume the two optimal dual solutions are different. The second fact means that for small enough $1>\lambda >0$, we can write
\[
    F(\lambda\bm{b}_1+(1-\lambda)\bm{b}_2) = \bm{p}_1\cdot\left(\lambda\bm{b}_1+(1-\lambda)\bm{b}_2\right).
\]
We also know that
\[
F(\lambda\bm{b}_1+(1-\lambda)\bm{b}_2) = \lambda F(\bm{b}_1) + (1-\lambda) F(\bm{b}_2),
\]
so by strong duality, and putting the above two equations together, we obtain
\begin{align*}
    \bm{p}_1\cdot\left(\lambda\bm{b}_1+(1-\lambda)\bm{b}_2\right) &= \lambda \bm{p}_1\cdot \bm{b}_1 + (1-\lambda)\bm{p}_2\cdot \bm{b}_2\\
    \bm{p}_1\bm{b}_2&=\bm{p}_2\bm{b}_2.
\end{align*}
Therefore $\bm{p}_1$ is also an optimal dual solution for $F(\bm{b}_2)$, which is a contradiction unless $\bm{p}_1=\bm{p}_2$.
\\
\Halmos
\endproof

\proof{Proof of Theorem~\ref{thm:ce-vs-ci}.}
We consider each design separately.

\textbf{Cost-excluded design:} We know from the proof of Theorem~\ref{thm:fraction-bias-SP-CE-vs-RCT-CE} that
\[
    \abs{\hat{\Delta}^{\text{prop}}_{\text{SP-CE}}-\Delta^{\text{prop}}} \le \left(\rho + (1-\rho)(1-\alpha)\right)(\bm{a}^{0}-\bm{a}^{1})\cdot\bm\beta,
\]
where the CE label is omitted for succinctness. Therefore the first condition in Theorem~\ref{thm:ce-vs-ci} is sufficient.

Similarly, in the fixed-cost model under either condition (i) or (ii), the proof of Theorem~\ref{thm:fraction-bias-SP-CE-vs-RCT-CE} tells us that
\[
\abs{\hat{\Delta}^{\text{fixed}}_{\text{SP-CE}}-\Delta^{\text{fixed}}} \le (\bm{a}^{0}-\bm{a}^{1})\cdot\bm\beta,
\]
so the second condition in Theorem~\ref{thm:ce-vs-ci} is also sufficient.

To see why both are also necessary, observe that we can write the bias of the shadow price estimator as
\[
    \hat{\Delta}^m_{\text{SP-CE}} - \Delta^m = \int_0^\rho (\bm{a}^{\text{\emph{CE}}} - \bm{a}^{\eta}) \cdot\bm\beta d\eta + \int_\rho^1 (\bm{\tilde{a}}^{m} - \bm{\tilde{a}}^{m,\eta})\cdot\bm\beta d\eta\nonumber.
\]
Because both terms are non-increasing in $\rho$, the bias can only be 0 for every $\rho$ if both integrands are zero, i.e., $(\bm{a}^{\text{\emph{CE}}} - \bm{a}^{\eta}) \cdot\bm\beta=0$ for all $\eta$ and $\rho$ and $(\bm{\tilde{a}}^{m} - \bm{\tilde{a}}^{m,\eta})\cdot\bm\beta=0$ for all $\eta$ and $\rho$. Under both the proportional-cost model and the fixed-cost model with condition (i) or (ii), this is equivalent to writing $\bm{a}^{\text{CE},0}\cdot\bm\beta = \bm{a}^{\text{CE},\eta} \cdot \bm\beta=\bm{a}^{\text{CE},1}\cdot\bm\beta$, or simply saying that the matching function is linear on the line segment from $\bm\lambda$ to $\bm\lambda+\bm\beta$. By Lemma~\ref{lem:sensitivity-analysis} this implies $\bm{a}^{\text{CE},0}=\bm{a}^{\text{CE},1}$. Note that if condition (i) or (ii) does not hold in the fixed-cost model, then we cannot satisfy that the second term of the bias is zero for every $\rho$.

\textbf{Cost-included design:} Define $\Psi^m(\eta)=\Phi_{\text{CI}}^m((1-\eta)\bm{\lambda}, \eta(\bm\lambda+\bm\beta),\bm\gamma)$ for any $\eta\in[0,1]$, such that $\Delta^m=\Psi(1)-\Psi(0)$.

We can write the derivative $\Psi'(\eta)=\bm{a}^{\treatment,\text{CI},m,\eta}\cdot(\bm\lambda + \bm\beta) - \bm{a}^{\control,\text{CI},m,\eta}\cdot\bm\lambda$, and observe that
\begin{align*}
    \hat{\Delta}^m_{\text{SP-CI}} - \Delta^m &= \int_0^1 \left(\bm{a}^{\treatment,\text{CI},m,\rho}\cdot(\bm\lambda + \bm\beta) - \bm{a}^{\control,\text{CI},m,\rho}\cdot\bm\lambda\right)-\left(\bm{a}^{\treatment,\text{CI},m,\eta}\cdot(\bm\lambda + \bm\beta) - \bm{a}^{\control,\text{CI},m,\eta}\cdot\bm\lambda\right)d\eta\\
    &=\int_0^\rho \left( \Psi'(\rho) - \Psi'(\eta) \right)d\eta + \int_\rho^1 \left( \Psi'(\rho) - \Psi'(\eta) \right)d\eta.
\end{align*}
Again, both terms are non-increasing in $\rho$, so the only way to ensure zero bias for all $\rho$ is to make sure $\Psi'(0)=\Psi'(\eta)=\Psi'(1)$ for all $\eta\in[0,1]$. From Lemma~\ref{lem:sensitivity-analysis} this implies that $\bm{a}^{\treatment,\text{CI},m,0}=\bm{a}^{\treatment,\text{CI},m,1}$ and $\bm{a}^{\control,\text{CI},m,0}=\bm{a}^{\control,\text{CI},m,1}$. This condition is therefore necessary, and it is also trivially sufficient.
\\
\Halmos
\endproof

\proof{Proof of Corollary~\ref{cor:ci-unbiased}.}
First, we observe that the ``if'' direction is trivial: if every demand type is matched to its highest-value supply type, there is no interference and all estimators are unbiased. We now focus on the ``only if'' condition.

For any demand type $i$, define $j(i)\coloneq \arg\max_{1\le k\le n_s}v_{i,k}$.
Assume that in the global treatment state, there exists a demand type $i_1$ and supply type $j_1\neq j(i)$ such that $x_{i_1,j_1}^{\treatment,CI,m,1} > 0$.

\textbf{Step 1:} We can begin by showing that if the above condition holds in global treatment, an analogous condition holds in global control. Assume for a contradiction that $
x_{i,j}^{\control,\text{{CI}},m,0} =\lambda_i \Leftrightarrow v_{i,j} = \max_{1\le k\le n_s} v_{i,k}$. Under the non-degeneracy assumption, we must also have that for all $1\le j\le n_s$,
\begin{equation}
\label{eq:supply-slack-1}
\sum_{i=1}^{n_d}x_{i,j}^{\control,\text{{CI}},m,0} < \gamma_j.
\end{equation}
Let $\varepsilon > 0$, and consider the following solution for $\rho=\varepsilon$:
\begin{align*}
    x_{i,j}^{\control,\text{CI},m,\varepsilon}&=
        (1-\varepsilon) x_{i,j}^{\control,\text{CI},m,0}\\
    x_{i,j}^{\treatment,\text{CI},m,\varepsilon} &= \varepsilon \frac{\lambda_i + \beta_i}{\lambda_i} x_{i,j}^{\control,\text{CI},m,0}.
\end{align*}
Clearly, this solution is demand-feasible:
\begin{align*}
    \sum_{j=1}^{n_s}x_{i,j}^{\control,\text{CI},m,\varepsilon} &= (1-\varepsilon)x_{i,j(i)}^{\control,\text{CI},m,0} = (1-\varepsilon)\lambda_i,\\
    \sum_{j=1}^{n_s}x_{i,j}^{\treatment,\text{CI},m,\varepsilon} &= \varepsilon \frac{\lambda_i + \beta_i}{\lambda_i} x_{i,j(i)}^{\control,\text{CI},m,0} = \varepsilon(\lambda_i + \beta_i).
\end{align*}
As for the supply constraints,
\[
\sum_{i=1}^{n_d} x_{i,j}^{\control,\text{CI},m,\varepsilon} + x_{i,j}^{\treatment,\text{CI},m,\varepsilon} = \sum_{i=1}^{n_d} \left(1+\frac{\varepsilon\beta_i}{\lambda_i}\right) x_{i,j}^{\control,\text{CI},m,0} < \gamma_j,
\]
where the final inequality follows from~\eqref{eq:supply-slack-1} as long as $\varepsilon$ is small enough. Furthermore, this constructed solution is optimal since each demand type $i$ (control or treatment) is assigned to its individually value maximizing supply type $j(i)$. Thus, for treatment fraction $\varepsilon$, $x_{i_1,j_1}^{\treatment,\text{CI},m,\varepsilon} = 0$. However, we know from Theorem~\ref{thm:ce-vs-ci} that the dual basis is the same for any treatment fraction, and this also implies that the primal basis is the same for any treatment fraction. Thus, we have a contradiction, and there must exist a demand type $i_0$ and supply type $j_0\neq j(i)$ such that in the global control state, $x_{i_0,j_0}^{\control,\text{CI},m,0} > 0$.

\textbf{Step 2: } We have shown that if there is a demand type matched to a suboptimal supply type in the global treatment state, there must also be a demand type matched to a suboptimal supply type in the global control state. We now show that this also leads to a contradiction, distinguishing between two cases.

\textit{Case 1 (no split demand type): } Assume that in the global control state, each demand type is assigned to at most one supply type, in other words, for every $i\in[n_d]$, there is at most one supply type $j$ such that $x_{i,j}^{\control,\text{CI},m,0} > 0$. The demand type $i_0$ then verifies that $x_{i_0,j(i_0)}^{\control,\text{CI},m,0} = 0$, which means that the supply constraint for type $j(i_0)$ must be tight. Therefore, there are $n_d(n_s-1) + n_d + 1$ active constraints and the optimum is degenerate at global control which is a contradiction.

\textit{Case 2 (split demand type): } There exists a demand type (call it $i_0$ again without loss of generality) such that there exist two supply types $j_0, j'_0$ with $x_{i_0,j_0}^{\control,\text{CI},m,0} > 0$ and $x_{i_0,j'_0}^{\control,\text{CI},m,0} > 0$. Choose a small $\delta > 0$ (without loss of generality let $v_{i_0,j_0} > v_{i_0, j'_0}$). By Theorem~\ref{thm:ce-vs-ci}, the optimal primal basis of the cost-included matching problem must be the same for any treatment fraction, i.e., we must have that $x_{i_0,j_0}^{\control,\text{CI},m,1-\delta} > 0$ and $x_{i_0,j'_0}^{\control,\text{CI},m,1-\delta} > 0$ for treatment fraction $1-\delta$.

We can choose $\delta$ small enough that the optimal basis does not change (recall that we assumed the optimal matching under global treatment is nondegenerate and unique). We can then construct the optimal solution $x_{i,j}^{\cdot,\text{CI},m,1-\delta}$ via a sequence of augmenting paths, none of which alter the optimal basis. Because $x_{i_0,j_0}^{\control,\text{CI},m,1-\delta} > 0$ and $x_{i_0,j'_0}^{\control,\text{CI},m,1-\delta} > 0$, there must be two possible augmenting paths to increase the control demand of type $i_0$. One augmenting path must increase the total value more than the other (by uniqueness of the matching optimum). Therefore, it is suboptimal to choose a convex combination of the two rather than the most value-increasing. This leads to a contradiction.

We see that violating the condition in Corollary~\ref{cor:ci-unbiased} means that the SP-CI duals cannot remain constant between global control and global treatment --- completing the proof.
\\
\Halmos
\endproof
\newpage
\section{Numerical Experiments for the Fixed Cost Setting} \label{Appendix-sec:fixed-cost-numerical-experiments}

In this section, we present our key findings in a finite-sample setting for the fixed cost model.

\subsection{Setup}

The setup is exactly the same as described in Section~\ref{sec:numerical-experiments}. The only difference is the focus we place now on the fixed cost model. 

\subsection{Bias Reduction}

Figure~\ref{fig:simulations-singed-biases-fixed-cost} displays the (signed) average biases of all estimators along with their standard errors for treatment assignment probabilities $\rho$ of $0.1, 0.3$, and $0.5$, and for three fixed treatment cost levels $\kappa$ of  $10\%$, $30\%$, and $50\%$ of the minimum matching value. Across all cases, the RCT-CE estimator overestimates the global treatment effect, inducing positive bias. The RCT-CI estimator exhibits negative bias in the undersupply regime and transitions to positive bias at certain supply-to-demand ratios. Notably, SP estimators consistently reduce bias compared to standard estimators. In addition, the simulation-based estimator shows significant bias for small $\rho$ values (asymmetric experiments).

\begin{figure}[t]
\caption{Signed Biases of all estimators in the fixed cost setting.} 
\begin{center}
\includegraphics[height=4.7in]{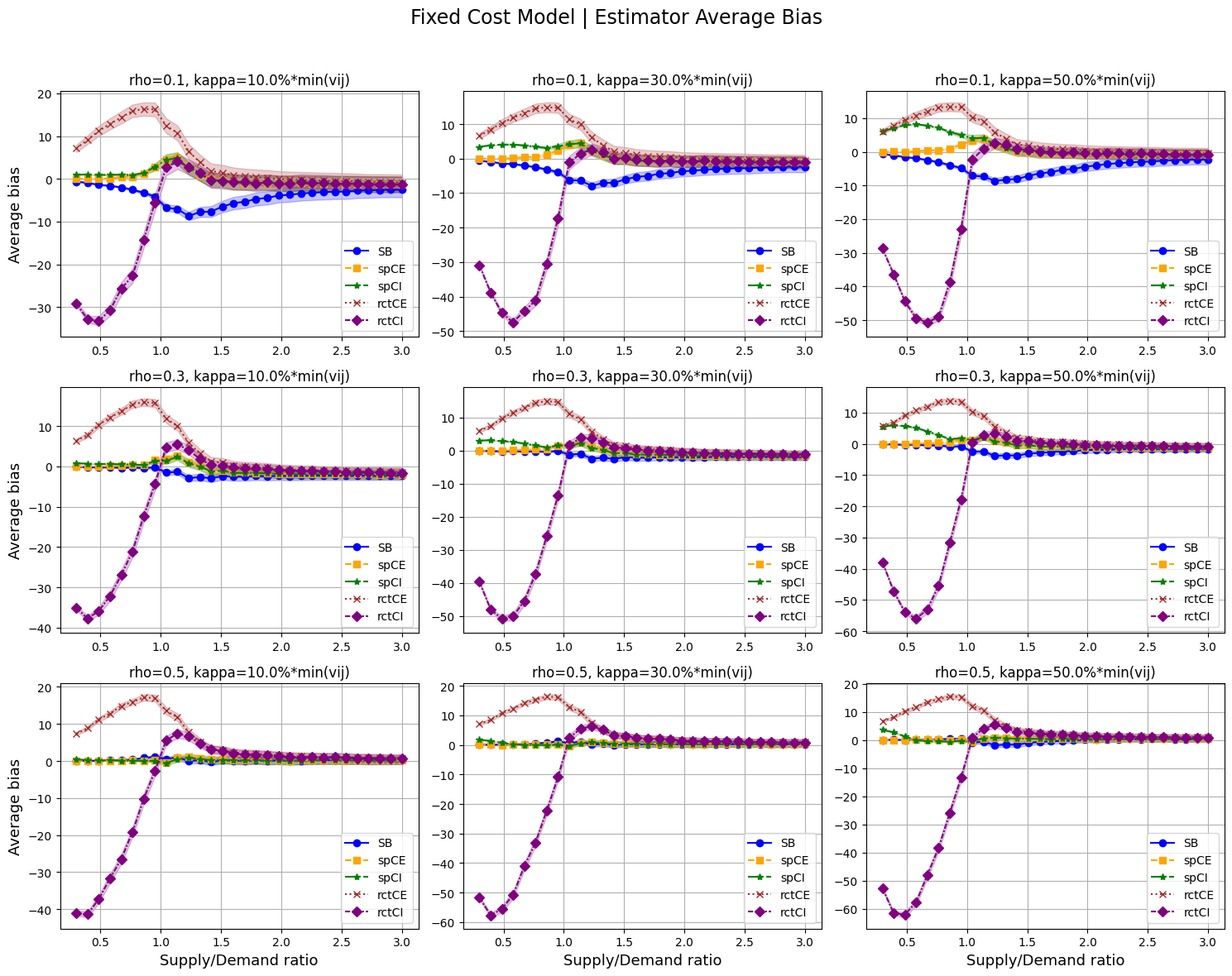}
\label{fig:simulations-singed-biases-fixed-cost}
\end{center}
{\footnotesize\textit{Note.} The shadow price estimator in the cost-excluded setting consistently reduces bias across all scenarios. The RCT-CE estimator is always positively biased, while the RCT-CI estimator shows significant negative bias in the undersupply regime. The simulation-based estimator (SB) exhibits significant bias in asymmetric experiments (small $\rho$).}
\end{figure}

Figure~\ref{fig:simulations-bias-reduction-fixed-cost} provides pairwise comparisons of shadow price and standard estimators across designs, as well as the shadow price estimators across designs in finite-sample settings. The left column highlights bias reduction by SP relative to RCT in the cost-excluding setting, while the middle column does so for the cost-included setting. Shadow price estimators significantly reduce bias in low-supply regimes. The third column illustrates improvements in bias reduction by SP-CE over SP-CI, particularly in low-supply settings. Overall, the SP estimator in the cost-excluded design outperforms others in bias reduction across all experimental configurations.

\begin{figure}[t]
\caption{Bias Reduction.} 
\begin{center}
\includegraphics[height=4.7in]{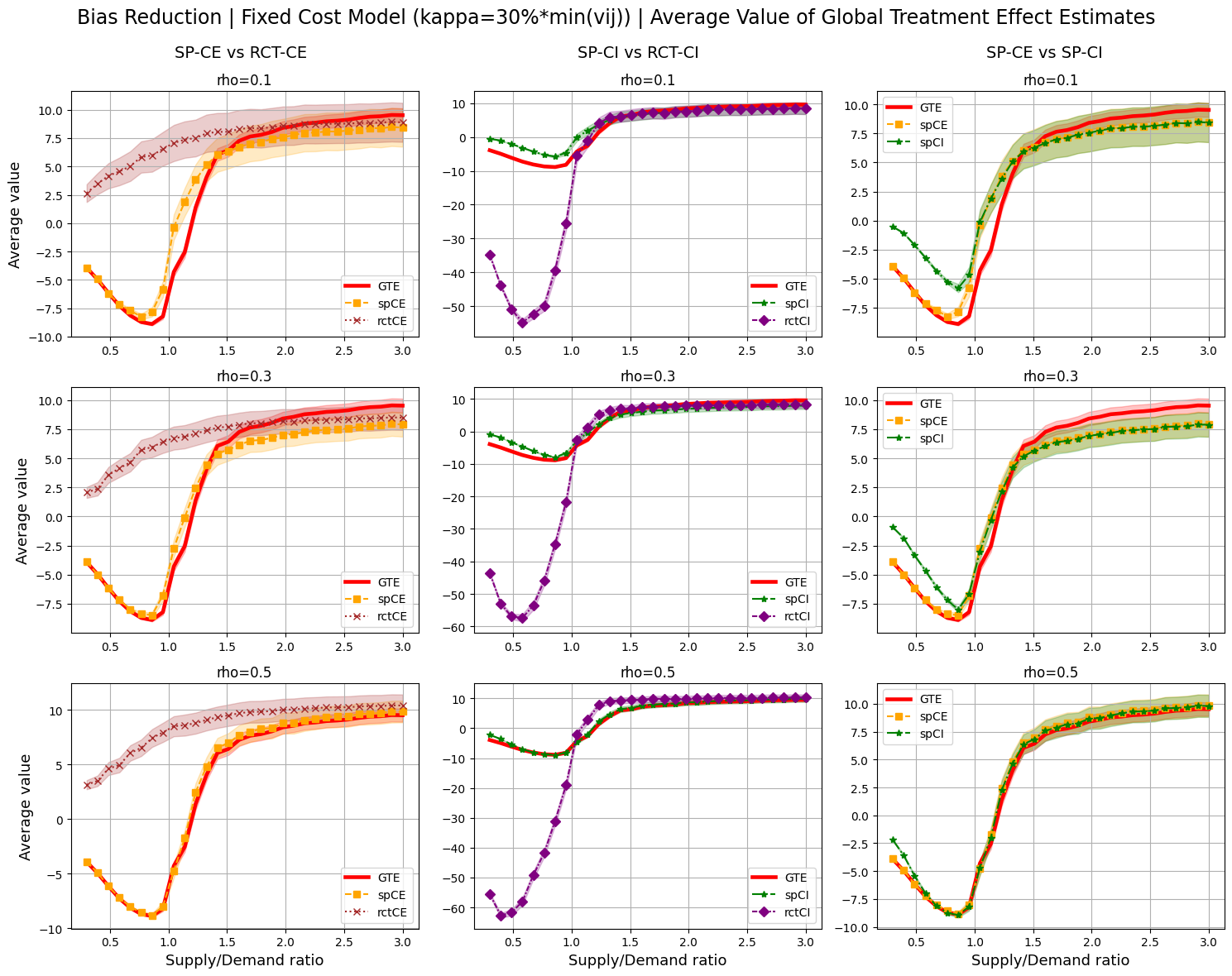}
\label{fig:simulations-bias-reduction-fixed-cost}
\end{center}
{\footnotesize\textit{Note.} Pairwise comparisons of shadow price and standard estimators across designs. The shadow price estimator in the cost-excluded design outperforms all others in bias reduction across all experiment configurations, particularly in the low-supply regime.}
\end{figure}

\subsection{Theoretical Guarantees in Practice}

This section examines the tightness of Theorem~\ref{thm:RCT-CE-vs-SP-CE} in finite-sample settings. Theorem~\ref{thm:RCT-CE-vs-SP-CE} provides sufficient conditions on the treatment assignment probability $\rho$ for the shadow price estimator to reduce bias compared to the standard estimator. Figure~\ref{fig:simulations-thm3-fixed-cost} shows that the theoretical guarantee is generally robust in undersupply and balanced settings, while being especially useful in oversupply scenarios. Importantly, the critical value of $\rho$ remains smaller than the actual turning point, indicating that shadow price estimator reduces bias over the standard estimator in the vast majority of cases. Only in highly asymmetric cases (with large $\rho$) in the oversupply scenarios may RCT show less bias than SP.

\begin{figure}[t]
\caption{Comparison of SP-CE bias vs RCT-CE bias.} 
\begin{center}
\includegraphics[height=2in]{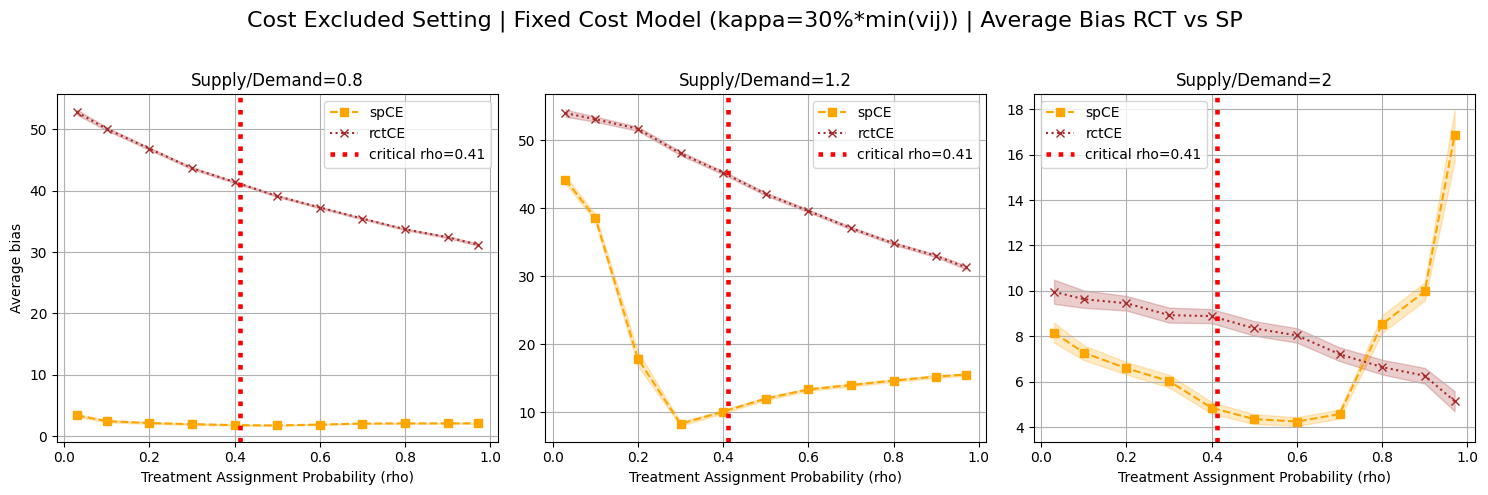}
\label{fig:simulations-thm3-fixed-cost}
\end{center}
\end{figure}

Figure~\ref{fig:simulations-thm4-fixed-cost} further explores the tightness of the upper bound on the bias ratio of SP versus RCT estimators in the cost-excluding setting. As expected, the bound is most meaningful in undersupply regimes, where the denominator of the right-hand side remains stable. In high-supply regimes, the denominator becomes very small, leading to high values on the right-hand side. While relatively loose in the low-supply regime, the bound provides a reasonable order of magnitude for the bias ratio between SP and RCT in the cost-excluding setting within finite-sample scenarios.

\begin{figure}[t]
\caption{Upper Bound of SP-CE vs RCT-CE bias ratio.} 
\begin{center}
\includegraphics[height=2.9in]{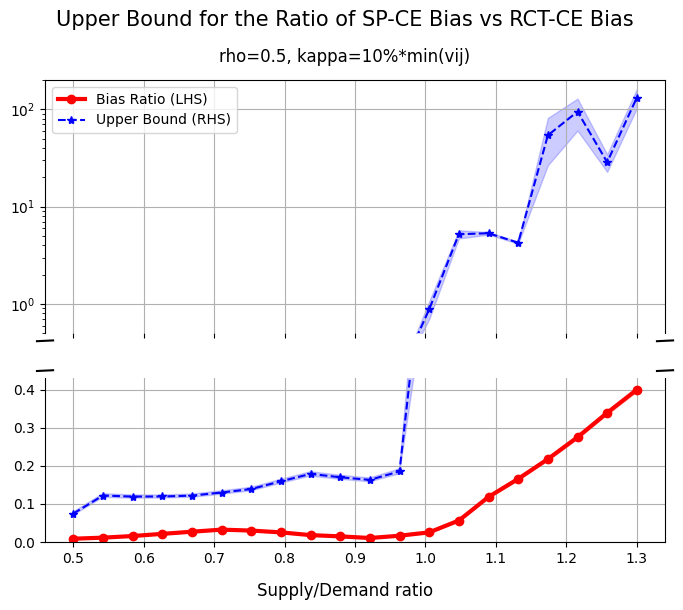}
\label{fig:simulations-thm4-fixed-cost}
\end{center}
{\footnotesize\textit{Note.} Tightness of the upper bound of the bias ratio between SP and RCT in the cost-excluding setting. The bound is most meaningful in undersupply regimes providing a reasonable order of magnitude, while it becomes large in high-supply settings.}
\end{figure}
\end{APPENDICES}

\end{document}